\documentclass[11pt,leqno]{article}

\usepackage{amsmath}
\usepackage{amssymb}
\usepackage{dsfont}
\usepackage{mathrsfs}
\usepackage{epsfig}
\usepackage{float}

\numberwithin{equation}{section}

\newtheorem{Theorem}{Theorem}
\newtheorem{Corollary}[Theorem]{Corollary}
\newtheorem{Lemma}[Theorem]{Lemma}
\newtheorem{Proposition}[Theorem]{Proposition}
\newtheorem{Definition}[Theorem]{Definition}
\newtheorem{remark}[Theorem]{Remark}

\newtheorem{example}[Theorem]{Example}

\numberwithin{Theorem}{section}

\bibstyle{plain}

\newcommand{\al}{\alpha}
\newcommand{\R}{{\mathbf R}}
\newcommand{\ds}{\displaystyle}
\newcommand{\e}{\varepsilon}
\newcommand{\Om}{\Omega}
\newcommand{\lra}{\longrightarrow}
\newcommand{\ra}{\rightarrow}
\newcommand{\p}{\partial}
\newcommand{\la}{\lambda}
\newcommand{\g}{\gamma}
\newcommand{\ov}{\overline}
\newcommand{\C}{{\mathbf C}}
\newcommand{\T}{{\mathbf T}}

\newcommand{\N}{{\mathbf N}}
\newcommand{\Z}{{\mathbf Z}}

\newcommand{\1}{{\mathds 1}}

\newcommand{\si}{\sigma}

\newcommand{\de}{\delta}

\newcommand{\ph}{\varphi}

\newcommand{\Xo}{{\mathcal X}}

\newcommand{\Lo}{{\mathcal L}}
\newcommand{\Do}{{\mathcal D}}

\newcommand{\Eo}{{\mathcal E}}

\newcommand{\Er}{{\mathcal E}(\R)}
\newcommand{\Ep}{{\mathcal E}}

\newcommand{\ii}{\int}
\newcommand{\ges}{\geqslant}
\newcommand{\les}{\leqslant}

\newcommand{\pr}[1]{\lfloor #1 \rfloor}
\newcommand{\PR}{\pr{P}}
\newcommand{\QR}{\pr{Q}}

\newcommand{\vs}{v_s}

\newcommand{\beq}{\begin{equation}}
\newcommand{\eeq}{\end{equation}}

\newcommand{\references}[1]{\theinstitutions 
\begin{thebibliography}{#1}}
\newcommand{\Endrefs}{\end{thebibliography}}

\title{Periodic traveling waves for nonlinear Schr\"odinger equations with non-zero conditions at infinity in $ \R ^2 $ } 

\date{ }

\author{ 
Mihai MARI\c S\footnote{Institut de Math\'ematiques de 
Toulouse UMR 5219, Universit\'e de Toulouse, 118 route de Narbonne, 31062 Toulouse 
Cedex 9, France. {\sf e-mail}: mihai.maris@math.univ-toulouse.fr. Corresponding author.}
$\quad \quad $ and $ \quad \quad  $
Anthony MUR\footnote{Institut de Math\'ematiques de 
Toulouse UMR 5219, Universit\'e  de Toulouse, 118 route de Narbonne, 31062 Toulouse 
Cedex 9, France. {\sf e-mail}: anthony.mur@math.univ-toulouse.fr.}
}

\begin{document}

\maketitle

\begin{abstract}
We consider the nonlinear Schr\"odinger equation with nonzero conditions at infinity in $\R^2$.  
We investigate the existence of traveling waves that are periodic in the direction transverse to the direction of propagation and minimize the energy when the momentum is kept fixed. 
We show that for any given value of the momentum, there is a critical value of the period such that traveling waves with period smaller than the critical value are one-dimensional, and those with larger periods depend on two  variables. 
\end{abstract}

\section{Introduction} 
\label{intro}

We consider  the 
nonlinear Schr\"odinger equation
\beq
\label{1.1}
i \frac{ \p \Phi}{\p t} + \Delta \Phi + F( | \Phi |^2) \Phi = 0
\qquad \mbox{ in } \R^2\times \R,  
\eeq
where $F$ is a real-valued function on $[0, \infty ) $ such that  $ F(r_0^2) = 0$ for some $r_0 >0$.
For a given $ \Lambda > 0$, we are interested in solutions $ \Phi : \R^2 \times \R \lra \C$ of (\ref{1.1})  that  are $ \Lambda -$periodic  with respect to the second variable, namely 
 $ \Phi ( x, y + \Lambda, t) = \Phi ( x, y , t )$ for all $ (x, y, t )$ 
  and satisfy the "boundary condition" 
$|\Phi (x ,y , t )  | \lra r_0 $ as $ x \lra  \pm \infty $.

If $ F'( r_0^2) <0 $ (which means that (\ref{1.1}) is defocusing), 
a simple scaling enables us to assume that $ r_0 = 1 $ and $ F'(r_0 ^2) = -1$ (see \cite{M10}, p. 108); 
we will do so throughout  this paper. The sound velocity at infinity associated to (\ref{1.1}) is then 
$ v_s = r_0 \sqrt{ -2 F' ( r_0^2) } = \sqrt{2}$ (\cite{M8}, p. 1077).

Equation (\ref{1.1}) has a  Hamiltonian structure.  
Indeed, let  $ V( s) =  \ii_s^{1} F( \tau ) \, d \tau$. 
It is then easy to see that if $ \Phi $ is a solution to (\ref{1.1}) and $ \Phi $ is $ \Lambda -$periodic with respect to the second variable, then at least formally, the "energy"
\beq
\label{1.2}
E(\Phi ) = \int_{\R \times [0, \Lambda]} \Big| \frac{\p \Phi}{\p x }  \Big|^2  + \Big| \frac{\p \Phi}{\p y }  \Big|^2 +  V( |\Phi |^2) \, dx \, dy
\eeq
does not depend on $t$.

In the sequel it is more convenient to  normalize the period to $1$: instead of working with functions $ \Phi $ 
that are $\Lambda-$periodic with respect to the second variable, we will consider the function 
$ \tilde{\Phi} ( x, y, t ) = \Phi ( x, \Lambda y, t )$ which is  $1-$periodic with respect to the second variable. 
It is clear that  $ \Phi $ is a solution of (\ref{1.1}) if an only if $ \tilde{\Phi}$  
satisfies the equation 
\beq
\label{1.3}
i \frac{ \p \tilde{\Phi}}{\p t} + \frac{\p ^2 \tilde{\Phi}}{\p x^2 } + \frac{1}{\Lambda ^2}  \frac{\p ^2 \tilde{\Phi}}{\p y ^2 } + F( |\tilde{\Phi} |^2) \tilde{\Phi} = 0
\qquad \mbox{ in } \R^2  \times \R.
\eeq
We denote $ \la = \frac{1}{\Lambda}$ and we consider the renormalized energy  
\beq
\label{Ela}
 E_{\la} ( \Psi ) =  \int_{\R \times [0, 1]} \Big| \frac{\p \Psi}{\p x }  \Big|^2  +{\lambda ^2} \Big| \frac{\p \Psi}{\p y }  \Big|^2 +  V( |\Psi |^2) \, dx \, dy. 
\eeq
If $ \Phi $ and $ \tilde{\Phi}$ are as above, we have $ E( \Phi ) = \frac{1}{\la } E_{\la } ( \tilde{\Phi}). $

We are interested in traveling waves for (\ref{1.1}), which are solutions of the form $ \Phi ( x, y , t ) = \psi ( x + ct, y )$. 
If $ \psi $ is a traveling wave and is $\Lambda-$periodic with respect to the variable $y$, then 
$ \tilde{\psi}(x, y ) = \psi( x, \Lambda y ) = \psi ( x, \frac{y}{\la })$ is $1-$periodic with respect to  $y$ and satisfies the equation 
\beq
\label{twp}
i c \frac{\p \tilde{\psi}}{\p x} + \frac{\p ^2 \tilde{\psi}}{\p x^2 } + \frac{1}{\Lambda ^2}  \frac{\p ^2 \tilde{\psi}}{\p y ^2 } + F( |\tilde{\psi} |^2) \tilde{\psi} = 0
\qquad \mbox{ in } \R^2.
\eeq

\medskip

{\bf Assumptions and some comments on the assumptions. }
We  work with general nonlinearities $F$. We will  consider the  set of assumptions {\bf (A1), (A2), (B1), (B2)}   described  below. 
We will assume throughout the paper that {\bf (A1)} holds. 
The other assumptions are not needed all the time and sometimes they can be slightly relaxed. 
For each result we will indicate the precise conditions that we use. 
For instance, assumption {\bf (A2)} is not necessary when we deal with one-dimensional traveling waves to (\ref{1.1}). 
Our aim is to consider sufficiently general nonlinearities and in the meantime to focus on ideas, avoiding irrelevant technicalities.

\medskip

{\bf (A1) } The function $F$ is continuous on $[0, \infty)$, 
$C^1$ in a neighborhood of $1$, $F(1) = 0$ and $F'(1) < 0$.

\smallskip

{\bf (A2) } 
There exist $C > 0$ and $1 < p_0 < \infty $ 
such that $|F(s) | \leq C(1 + s^{p_0}) $ for any $ s \geq 0$.

\smallskip

 We denote 
$$
V( s ) = \ii_{s }^1 F( \tau ) \, d \tau, 
$$
so that $ V(1 ) = 0 $ and $ V'(s ) = - F (s) $. If assumption {\bf(A1)} holds, we have 
\beq
\label{A1-behaviour}
V(s ) = \frac 12 ( s-1 )^2 + o \left( (s - 1)^2 \right)  \qquad \mbox{ as } s \lra 1 .
\eeq
If {\bf(A2)} holds, there is $ C'> 0 $ such that $ | V( s ) | \les C' s^{ p_0 + 1} $ for all $ s \ges 2 $.

The natural function space associated to (\ref{1.3}) is 
$$
\begin{array}{l}
\Ep  =  \ds \{ \psi \in L_{loc}^1( \R ^2 ) \; \mid \; \psi \mbox{ is $1-$periodic with respect to the second variable,  }
\\
\qquad \qquad \ds \nabla \psi \in L^2 ( \R \times [0,1]) \mbox{ and } V( |\psi |^2) \in L^1 ( \R \times [0,1]) \}.
\end{array}
$$
We will also consider the one-dimensional variant of $ \Ep$, namely 
$$
\Er = \{ \psi \in L_{loc}^1( \R ^2 ) \; \mid \; \psi ' \in L^2 ( \R) \mbox{ and } V( |\psi |^2 ) \in L^1 ( \R) \}
$$
and the associated $1-$dimensional energy 
$$
E^1( \psi ) = \ii_{\R} |\psi ' (x)|^2 + V( |\psi |^2 ) ( x) \, dx .
$$

As we can see, assumption {\bf (A1)} determines the behaviour of the nonlinear potential $V$ in a neighbourhood of $1$,   
and {\bf (A2)} gives upper bounds on $V(s ) $ for large $s$.  
In view of  (\ref{A1-behaviour}), the  renormalized Ginzburg-Landau energy 
$$
E_{GL, \la } ( \psi ) = \int_{\R \times [0, 1]} \Big| \frac{\p \psi}{\p x }  \Big|^2  +{\lambda ^2} \Big| \frac{\p \psi}{\p y }  \Big|^2 + 
\frac 12 \left(  1 -  |\psi |^2 \right)^2 \, dx \, dy,  
$$
together with its $1-$dimensional variant 
$$
E_{GL}^1 ( \psi ) = \ii_{\R} |\psi ' (x)|^2 +\frac 12 \left(  1 -  |\psi (x) |^2 \right)^2 \, dx \qquad \mbox{ for } \psi \in \Er  
$$
will be relevant throughout the article. 
Notice that the renormalized Ginzburg-Landau energy is simply the energy associated to the Gross-Pitaevskii nonlinear potential
 $ V( s ) =\frac 12( 1 - s)^2 $ corresponding to the nonlinearity $ F( s ) = 1-s$.

Whenever minimization of energy at fixed momentum is considered, we need to  assume that $ V \ges 0 $ on $[0, \infty )$
(for otherwise, the infimum is $ - \infty$). 
If $ \psi $ is any one-dimensional finite-energy traveling wave of speed $ c \neq 0$, equation (\ref{4.3}) implies that for any $ x \in \R$ we have either $ |\psi ( x ) | = 1$, or $ V( |\psi (x) |^2 ) > 0$.   
For simplicity, we will  assume that $ V (s ) > 0 $  for $ s \neq 1$. 
We discuss in Example \ref{Example3} what happens if there exists $ s_0 \in (0,1)$ such that  $ V (s_0) = 0$,  and $ V > 0 $ on $(1, \infty)$. 

The next assumptions give lower bounds on 
$V$ and will be useful to estimate the Ginzburg-Landau energy $E_{GL}(\psi)$ in terms of $ E(\psi)$; see Lemmas \ref{L2.1} and  \ref{MutualBounds} in the next section.

\medskip

{\bf (B1) } We have $ V > 0 $ on $ [0, \infty ) \setminus \{ 1 \} $ and, denoting $ H(s ) = \ds \ii_1^s \big| V(\tau ^2) \big|^{\frac 12} d \tau$, we have 
$$
H(s) \lra \infty \qquad \mbox{ as } s \lra \infty. 
$$

\smallskip

{\bf (B2) }  $ V > 0 $ on $ [0, \infty ) \setminus \{ 1 \} $ and there exists $ \g > 0 $ and $ s_ 0 \ges 1 $ such that $ V( s ) \ges s^{\g}$ for all $ s \ges s_0$. 

\medskip

Clearly, {\bf (B2)} is stronger than {\bf (B1)}.
If $ V > 0 $ on $ [0, \infty ) \setminus \{ 1 \} $, then the function $H$ is strictly increasing on $[0, \infty )$ and $H(1) = 0 $. 
If we do not assume that $H(s) \lra \infty $  as $ s \lra \infty$, it is possible to construct sequences of functions $\psi_n : \R \lra \C$ 
such that $ E^1 ( \psi _n )$ is bounded, but the Ginzburg-Landau energy $ E_{GL} ^1( \psi _n )$ is unbounded, a situation that we would like to avoid. 
Assumption {\bf (A1)} is enough to study  $1-$dimensional traveling-waves of (\ref{1.1}) by using ODE arguments, as we do in section 4. 
Assumptions  {\bf (A1)} and {\bf (B1)} are sufficient to obtain $1-$dimensional traveling-waves by minimizing the energy at fixed momentum, as in section 5. 

When we 
consider functions defined on the whole space $ \R^N$, $ N \ges 2$, assumptions {\bf (A1)} and {\bf (A2)} (with $p_0 < \infty$ if $N=2$, and $ p_0 < \frac{2}{N-2} $ if $N \geqslant 3$) and 
the fact that $ V \ges 0 $ 
are sufficient to prove that $E_{GL}( \psi )$ is bounded whenever $ E( \psi)$ is bounded, and vice-versa; see Lemma 4.8 p. 177 in \cite{CM}. 
This is no longer true in the one-dimensional space $ \R$ or in a strip $ \R \times(0,1)$.

\medskip

{\bf Momentum. } 
There is another important quantity conserved by equations (\ref{1.1}) and  (\ref{1.3}), namely {\it the momentum.}
It carries some topological information and 
defining rigorously the momentum is a difficulty in itself. We will address this issue  in Section 3. 
Roughly speaking, the momentum is a functional $Q$ whose G\^ateaux differential is $ Q'( \psi ) = 2 i \frac{ \p \psi}{\p x}$. 
If $ \psi \in \Ep  $ is a function  such that there exists $ A \ges 0 $ such that  $ \psi ( x, y ) = 1 $ if $ x \les - A$ and if $ x \ges A$,
 we have $ Q( \psi ) =  \ii_{\R \times [0, 1] } \langle i \frac{\p \psi}{\p x} , \psi \rangle \, dx \, dy. $
Given an arbitrary function $ \psi \in \Ep$, the function $ \langle i \frac{\p \psi}{\p x} , \psi \rangle$ 
does not necessarily belong to $L^1( \R \times [0, 1])$ and 
giving a meaning to the above integral is not obvious. 
For functions $ \psi \in \Er$ there is a one-dimensional variant of the momentum that we will denote $ P(\psi)$.  
 We will see that the momentum can be defined only modulo $ 2\pi$.

\medskip

{\bf Brief description of the results. } 
In this article we  focus on traveling waves  for (\ref{1.1}) that minimize  the energy when the momentum is kept fixed. 
In view of a celebrated result by T. Cazenave and P.-L. Lions \cite{CL}, such solutions are expected to be orbitally stable by the flow associated to 
(\ref{1.1}). 

We  denote 
$$
E_{\la, \min} (p) = \inf \{ E_{\la} ( \psi ) \; \mid \; \psi \in \Ep \mbox{ and } Q( \psi ) = p \}, 
$$
and 
$$
E_{\min}^1 ( p ) = \inf \{ E_{\la} ( \psi ) \; \mid \; \psi \in \Er \mbox{ and } P( \psi ) = p \}. 
$$

The main results of this paper  can be summarized  as follows. 

\begin{Theorem}
\label{T1.1}
Assume that the conditions {\bf (A1), (A2), (B2)} above are satisfied. 
Then: 

\medskip

i) The function $E_{\min}^1 $ is nonnegative, $2\pi-$periodic,  concave on $[0, 2\pi]$,  $ E_{\min}^1 ( - p ) = E_{\min}^1 ( p )$,
$E_{\min}^1 ( p ) \les \sqrt{2} |p| $ and the right-derivative of $ E_{ \min}^1 $ at the origin is $ \sqrt{2}$. 
If $F$ is $C^2$ near $1 $ and $ F''(1) < \frac 94$ we have $ E_{ \min}^1 ( p ) < \sqrt{2} |p| $ for any $p \neq 0 $.

For any $ p \in (0,  \pi]$ satisfying $ E_{ \min}^1 ( p ) < \sqrt{2} p$ there exist minimizers for $ E_{ \min}^1 ( p )  $ in $ \Er$ and any minimizing sequence 
contains a convergent 
\footnote{This statement is vague because we did not introduce a distance on $ \Er$. Please see Theorem \ref{T5.2} for a precise statement.}  subsequence modulo translations. 

\medskip

ii) For any fixed $ \la > 0 $ the function $p \longmapsto E_{\la, \min} (p)$ is $2\pi-$periodic, concave on $[0, 2\pi]$, 
$E_{\la, \min} ( p ) \les \sqrt{2} |p| $, the right-derivative of $ E_{\la, \min} $ at the origin is $ \sqrt{2}$, and 
$E_{\la, \min} ( p ) \les E_{ \min}^1 ( p )$. 

\medskip

iii) Assume that $E_{\min}^1  $ is twice differentiable at some $ p \in (0,  \pi ) $ 
and  $E_{\min}^1  ( p  ) < \sqrt{2} p$.

Then there exists $ \la_s ( p ) \in (0, \infty)  $ such that the mapping $ \la \longmapsto E_{\la, \min}(p) $ is strictly increasing on 
$(0, \la_s (p) ] $, and  $ E_{\la, \min}(p) = E_{\min}^1 ( p ) $ for all $ \la \ges \la_s ( p ) $. 
Therefore  minimizers for $ E_{\la, \min}(p) $ - whose existence is guaranteed by the next statement - depend on both variables $ x $ and $ y $ if $ \la < \la_s ( p ) $, and they depend only on $x$ if $ \la \ges \la_s ( p )$. 

\medskip

iv) For any $ \la > 0 $ and any $ p \in (0,  \pi]$ satisfying $ E_{\la,  \min} ( p ) < \sqrt{2} p$ there exist minimizers for $ E_{\la,  \min} ( p )  $ in $ \Ep$ and any minimizing sequence  has a convergent subsequence \footnote{See Theorem \ref{T6.6} for a precise statement. } after translations. 

\end{Theorem}

\medskip

{\bf Some related results in the literature.} 
In the case of the Gross-Pitaevskii equation, similar results have been obtained in \cite{dLGS1}  (corresponding to part  (iii) in Theorem \ref{T1.1}) and \cite{dLGS2} (parts (ii) and (iv) in Theorem \ref{T1.1}). 
Our approach is different. 
Our proofs are simpler, shorter and do not depend on the algebraic properties of the 
Gross-Pitaevskii nonlineaity $F(s) = 1- s$.

It has been proved in \cite{BGS-survey} that one-dimensional traveling-waves for  the Gross-Pitaevskii equation are unique (up to the invariances of the problem) and that any such solution is a minimizer of the energy at fixed momentum. 
This is no longer true for general nonlinearities; see examples \ref{Example1} and 
\ref{Example2} below. 
Since the construction of periodic traveling waves in dimension 2 relies on the properties of one-dimensional solutions that minimize the energy at fixed momentum,  
we need to develop the one-dimensional theory in sections 4 and 5. 
We mention that in the case of  general nonlinearities, various other properties of one-dimensional traveling-waves (e.g. transonic limit, 
stability or instability) have been studied in \cite{C1d} and \cite{C1d-stab}, and
nonexistence for supersonic speeds has been established in \cite{M8}.

The problem of minimizing the energy 
$$ \mathscr{E}  ( u ) = \int_{M_y ^k } \int_{\R^n} | \nabla _x u |^2 + \lambda |\nabla_y u |^2 - \frac{2}{ 2 + \al } |u |^{ 2 + \al} dx \, d\mbox{vol}_{M_y^k} 
$$
subject to the constraint $ \| u \|_{L^2 ( \R^n \times M_y^k ) }= 1 $, 
where $M_y^k$ is a $k-$dimensional compact Riemannian manifold, has been studied in \cite{TTV} in connection with the Schr\"odinger equation 
$ i \p_t u - \Delta_{x, y } u - |u |^{\al } u = 0 $ on $ \R_x^n \times M_y ^k $. 
Theorem 3.1 p. 78 in \cite{TTV} states that there exists $ \la ^* \in (0, \infty)$ 
such that minimizers depend only on $ x \in \R^n $ if $ \la > \la ^*$ and depend both on $ x $ and on $ y \in M_y^k$ if $ \la < \la ^*$.  
There is an obvious analogy between this result and Theorem \ref{T1.1} (iii) and (iv) above. 
However, the problem considered in \cite{TTV} is very different from the problem that we study in this article, and the ideas in proofs are completely different, too.

\medskip

{\bf Outline of the paper. } 
In the next section we present the functional setting and we show that one can bound the energies $E^1 $ and $E_{\la}$ in terms of the associated Ginzburg-Landau energies, and vice-versa. 
In Section 3 we give  a rigorous definition of the momentum, firstly for functions in $ \Er$ and then for functions  in $\Ep$. 
Section 4 is devoted to the study of one-dimensional traveling waves by using ODE arguments. 
In the case of the Gross-Pitaevskii equation $(F(s ) = 1-s)$ those results were already known, see \cite{BGS-survey} and references therein.
Different behaviour may occur when we consider other nonlinearities. 
In section 5 we study the $1-$dimensional minimization problem associated to $E_{\min}^1$ and we prove part (i) in Theorem \ref{T1.1}. 
Section 6 is devoted to the $2-$dimensional minimization problem for $E_{\la, \min}$, and to  the rest of the proof of Theorem \ref{T1.1}.

Proofs for the statements in  Theorem \ref{T1.1}  can be found as follows. 
The first part of (i) is contained in Lemma \ref{E1min}, and the last statement is Theorem \ref{T5.2}. 
Parts (ii) and (iii) are contained in Lemma \ref{Emin}, and (iv) is Theorem \ref{T6.6}. 

\medskip

{\bf Notation. } 
We denote by $\langle  \cdot, \cdot  \rangle $ the usual scalar product in $ \C \simeq \R^2$, namely 
$ \langle a + i b , c + id \rangle = ac + bd$, and by $ \Lo ^N$ the Lebesgue measure in $ \R^N$. 
We denote by $ C_{a_1, \dots, a_{\ell}}$ or by $ C( a_1, \dots, a_{\ell}) $ a positive constant that may change from line to line, 
and that  depends only on the parameters $ a_1, \dots, a_{\ell}$. 

\section{Energy and function spaces}

The precise representative of a given  function $ f \in L_{loc}^1( \R^N)$ is the function 
$ f^*$ defined on $ \R^N$ by $ f^*(x) = \ds \lim_{r \ra 0 } \frac{1}{| B( x, r )| } \ii_{B(x, r )} f( y ) \, dy $ if the limit exists, and 0 otherwise. 
It is well-known that for any $ f \in L_{loc}^1( \R^N)$ we have $ f = f^*$ almost everywhere
(see, e.g., Corollary 1 p. 44 in \cite{EG}). 
In the sequel we will always replace functions in $ L_{loc}^1( \R^N) $ by their precise representatives. 
We  denote 
$$
\Er  =  \{ \psi \in L_{loc}^1 ( \R) \; \mid \; \psi '\in L^2( \R), \, V( |\psi |^2) \in L^1( \R) \}. 
$$

For any given $ \psi \in \Er$ and any interval  $  I\subset \R$ we denote
$$
E^{1, I}( \psi ) = \int_{I} |\psi '|^2 + V( |\psi |^2) \, dx 
\qquad \mbox{ and } \qquad 
E_{GL}^{1, I}( \psi ) = \int_{I} |\psi '|^2 + \frac 12( |\psi |^2 - 1)^2 \, dx . 
$$
When $ I = \R$ we write simply $ E^{1}( \psi ) $ and $ E_{GL}^{1}( \psi ) $. 

\begin{Lemma}
\label{L2.1}
Assume that  conditions {\bf (A1)}  and {\bf (B1)} hold.  
Then: 

i) Any function $ \psi \in \Er$ is bounded and $\frac 12-$H\"older continuous on $ \R$. 
There exists a function $ b : [0, \infty ) \lra [0, \infty )$ satisfying $\ds \lim_{\tau \ra 0 } b( \tau ) = 0 $  
 such that 
$$
\| \, | \psi | - 1 \|_{L^{\infty}( \R)} \les b ( E^1( \psi )) \qquad \mbox{ for any } \psi \in \Er. 
$$

ii) For any $ \psi \in \Er $ we have $ |\psi | - 1 \in H^1( \R)$ and 
$$
\begin{array}{rcl}
\Er & = & \{ \psi \in H_{loc}^1( \R ) \; \mid \; \psi ' \in L^2( \R) \mbox{ and } |\psi |^2 - 1 \in L^2 ( \R) \} 
\\
\\
& = & \{ \psi \in H_{loc}^1( \R ) \; \mid \; \psi ' \in L^2( \R) \mbox{ and } |\psi | - 1 \in L^2 ( \R) \} . 
\end{array}
$$
Furthermore, there exist functions $ b _1, b_2 : [0, \infty ) \lra [0, \infty )$ such that  $\ds \lim_{\tau \ra 0 } b_i( \tau ) = 0, $ 
$\ds \lim_{\tau \ra \infty } b_i( \tau ) = \infty $ for $ i = 1, 2$,   and we have
$$
E^1( \psi )  \les b_1 ( E_{GL}^1( \psi )) \quad \mbox{ and } \quad   E_{GL}^1( \psi )  \les b_2 ( E^1( \psi )) 
\quad \mbox{ for all } \psi \in \Er. 
$$
For any  $\psi \in \Er $ and any $ v \in H^1 ( \R)$ we have  $ \psi + v \in \Er$. 
\end{Lemma}

{\it Proof. } (i) Let $ \psi \in \Er$. Since $ \psi \in L_{loc}^1( \R)$ and $ \psi ' \in L_{loc}^1( \R)$, it follows from Theorem 8.2 p. 204 in  
\cite{brezis} that  $ \psi $ is equal almost everywhere to a continuous function and 
$$
\psi ( b ) - \psi ( a ) = \ii_a^b \psi '( s ) \, ds. 
$$
The Cauchy-Schwarz inequality gives 
\beq
\label{holder}
|\psi ( b ) - \psi ( a ) | \les | b - a | ^{\frac 12} \| \psi '\|_{L^2 ([a, b ])} \qquad \mbox{ for any } a, b \in \R. 
\eeq
It is well-known that $ |\psi | \in L_{loc}^1( \R)$ and $ \big| \, |\psi |'\, \big| \les |\psi '| $ almost everywhere. 
Using the Cauchy-Schwarz inequality and denoting $ h( s ) = \sqrt{ V( s^2)}$ we get for any $ a, b \in \R$, $ a < b $, 
\beq
\label{escape}
\begin{array}{l}
\ds \ii_a ^b |\psi '|^2 + V ( |\psi |^2 ) \, dx \ges \ii_a^b  \big| \, |\psi |'\, \big| ^2 + h^2 ( |\psi | ) \, dx 
\\
\\
\ds \ges 2 \bigg| \ii_a ^b h ( |\psi |) \cdot |\psi |' \, dx \bigg| 
= 2 \Big| H ( |\psi ( b ) |) - H ( |\psi (a) | ) \Big|. 
\end{array}
\eeq

Using (\ref{escape}) and  {\bf (B1)} it is easily seen that any function $ \psi \in \Er $ is bounded. 
Since $ V( |\psi |^2 ) \in L^1 ( \R)$, there are sequences $ a_n \lra - \infty $ and $ b_n \lra \infty $ such that $ |\psi ( a_n ) | \lra 1 $ and $ |\psi ( b_n ) | \lra 1$. 
For any $ x \in \R$ we use (\ref{escape}) on $ [a_n, x]$ and on $[x, b_n ]$ and we let $ n \lra \infty $ to get 
$ 4 \big| H ( |\psi ( x ) | ) \big | \les E^1 ( \psi )$. 
We infer that 
$$ \big| \, |\psi ( x ) |  -1 \big| \les \max \left( H^{-1}\left( \frac 14   E^1 ( \psi )\right) - 1 ,  1 - H^{-1}\left(-  \frac 14   E^1 ( \psi )\right) \right) \qquad \mbox{ if } 0 \les |\psi(x ) | \les 2 , 
$$
respectively $ | \psi ( x ) | -1 \les  H^{-1}\left( \frac 14   E^1 ( \psi )\right) - 1 $ if $ |\psi ( x )| >2$, and (i) is proven. 

\medskip

(ii) Assume that $ \psi ' \in L^2( \R)$ and $ V ( |\psi |^2 ) \in L^1( \R)$. 
By {\bf (A1)} and {\bf (B1)} there exists $C > 0 $ such that $ ( 1 - s ) ^2 \les C V(s) $ for all $ s \in [0,4]$. Then 
$
\left( 1 - |\psi |^2\right)^2  \1_{\{ |\psi | \les 2 \}} \les C V( |\psi |^2 ) \in L^1( \R). 
$

The set $ A = \{ x \in \R \; \mid \; |\psi (x ) | > 2 \} $ has finite Lebesgue measure in $ \R$. 
Indeed, by (i) we know that $ \psi $ is bounded. 
If $ m := \sup |\psi | > 2$, let $ i(m ) = \ds \inf_{s \in [4, m ^2]} V( s)$. 
Since $ V > 0 $ on $(1, \infty)$, we have  $ i(m ) > 0 $ and we obtain  the rough estimate
$ 
\ii_{\R} V(| \psi |^2 ) \, dx \ges \ii_{A} V( |\psi |^2) \, dx \ges  i(m ) \Lo ^1( A),  
$
hence $ \Lo ^1 ( A) \les \frac{1}{i ( m )}\|  V(| \psi |^2 )\|_{L^1( \R)}   $ and 
$$ \ii_{A} \left( 1 - |\psi |^2\right)^2 \, dx
\les ( m^2 - 1) ^2 \Lo ^1 ( A) \les \frac{( m^2 - 1)^2}{i(m) } \|  V(| \psi |^2 )\|_{L^1( \R)} .
$$ 

It follows from the above estimates that $1 - |\psi |^2 \in L^2( \R)$ and there exists  a function $ b_2 $ with the desired properties. 

Conversely, if $1 - |\psi |^2 \in L^2( \R)$ and $ \psi '\in L^2( \R)$, a similar argument shows that $ \psi $ is bounded, $ V( |\psi|^2 ) \in L^1( \R)$ and $E^1( \psi )$ 
can be estimated in terms of $E_{GL}^1( \psi )$.

We have $ \big| \, |\psi | - 1 \big| = \frac{ \big| \, |\psi |^2 - 1 \big|}{|\psi | + 1 } \les  \big| \, |\psi |^2 - 1 \big|$, hence 
$  |\psi | - 1 \in L^2( \R)$ whenever $  |\psi |^2 - 1  \in L^2( \R)$. 
Conversely, if $ \psi '\in L^2( \R)$ and $  |\psi | - 1 \in L^2( \R)$ we see as above that $ \psi $ is bounded and therefore 
$  |\psi |^2 - 1  \in L^2( \R)$.  

The fact that $ \psi + v \in \Er$ whenever   $\psi \in \Er $ and $ v \in H^1 ( \R)$ is proven  as  Lemma \ref{PG-bis} (i) below. 
\hfill
$\Box $

\medskip

The following result is contained in Theorem 1.8 p. 134 in \cite{PG2}: 

\begin{Lemma} \label{PG}
(\cite{PG2}) Let $ \psi \in \Er$. Then: 

i) There exist a real-valued function $ \ph \in \dot{H}^1(\R)$ 
 and $ w \in H^1( \R)$ such that $ \psi = e^{i \ph } + w$.

ii) If $(\ph_1, w_1)$ and $ (\ph_2, w_2)$ are as above, there exist $ k_-, k_+ \in \Z$ such that 
$ \ph_1 - \ph_2 - 2 \pi k_{\pm } \in L^2( \R_{\pm})$. 

iii) Moreover, the function $ \ph $ can be chosen such that $ \ph \in C^{\infty} ( \R)$  and
$ \ph ^{(k)} \in L^2( \R)$ for any $ k \in \N^*$.

\end{Lemma}

As already mentioned in the introduction, 
the natural "energy space" for the study of (\ref{1.3}) is 
\beq
\label{Ep}
\begin{array}{rcl}
\Ep & = & \{ \psi \in L_{loc}^1( \R ^2 ) \; \mid \; \psi \mbox{ is $1-$periodic with respect to the second variable and }\\
& & 
\nabla \psi \in L^2( \R \times [0,1] ) \mbox{ and } V( |\psi|^2 ) \in L^1 ( \R\times [0,1]) \} .
\end{array}
\eeq
Obviously, for any $ \psi \in \Er$ the function $ \psi^{\sharp} $ defined by $ \psi^{\sharp}(x, y ) = \psi ( x)$ belongs to $ \Ep$. 
We will also denote 
$$
H_{per}^1 = \{ v \in H_{loc}^1( \R^2) \; \mid \; v \mbox{ is $1-$periodic in the second variable and } v \in H^1( \R \times (0,1) ) \} .
$$

It is clear  that $E_{\la } ( \psi )  $ is well-defined for any $ \psi \in \Ep $ and for any $ \la > 0$, where $E_{\la}$ is as in (\ref{Ela}). 
We will show in Lemma \ref{MutualBounds} below that for any $ \psi \in \Ep $ we have  $\left( |\psi |^2 - 1 \right)^2  \in L^1 ( \R\times [0,1])$. 
For any $ \psi \in \Ep$, $ \la > 0 $ and for any interval $ I \subset \R$ we denote 
$$
\begin{array}{l}
\ds E_{\la }^I ( \psi ) = \ii_{I \times [0,1] } \Big| \frac{ \p \psi }{\p x} \Big|^2 + \la ^2  \Big| \frac{ \p \psi }{\p y} \Big|^2 + V( |\psi |^2 ) \, dx \, dy 
\quad \mbox{ and } 
\\
\\
\ds E_{GL, \la }^I ( \psi ) = \ii_{I \times [0,1] } \Big| \frac{ \p \psi }{\p x} \Big|^2 + \la ^2  \Big| \frac{ \p \psi }{\p y} \Big|^2 
+ \frac 12 \left( |\psi |^2 - 1 \right)^2 \, dx \, dy.
\end{array}
$$
We will simply write  $E_{\la }(\psi)$, respectively $ E_{GL, \la } ( \psi ) $ when $ I = \R$. 
We will write $ E( \psi ) $ and $E_{GL}( \psi )$ when $ \la = 1$.

\begin{Lemma}
\label{MutualBounds}
Assume that the conditions {\bf (A1), (A2)} and {\bf (B2)} in the introduction are satisfied. Let $ \Ep $ be as in (\ref{Ep}). Then we have
$$
\begin{array}{rcl}
\Ep & = & \{ \psi \in H_{loc}^1( \R ^2 ) \; \mid \; \psi \mbox{ is $1-$periodic with respect to the second variable and }\\
& & 
\nabla \psi \in L^2( \R \times [0,1] ) \mbox{ and } |\psi | - 1 \in L^2 ( \R\times [0,1]) \} 
\\
\\
& = & \{ \psi \in H_{loc}^1( \R ^2 ) \; \mid \; \psi \mbox{ is $1-$periodic with respect to the second variable and }\\
& & 
\nabla \psi \in L^2( \R \times [0,1] ) \mbox{ and } |\psi |^2 - 1 \in L^2 ( \R\times [0,1]) \} .
\end{array}
$$
In particular, for any $ \psi \in \Ep $ we have $ |\psi | - 1 \in H_{per}^1$. 

Moreover, for any $ \la > 0 $ there exist $ a, b, c, d > 0 $ such that for all $ \psi \in \Ep $ and for any interval $ I \subset \R$ of length at least $1$ we have 
\beq
\label{Mutual-Bounds}
E_{\la } ^I ( \psi ) \les  a E_{GL, \la } ^I ( \psi )  + b E_{GL, \la } ^I ( \psi ) ^{ p_0 + 1} 
\qquad \mbox{ and } \qquad 
E_{GL, \la } ^I ( \psi )  \les c E_{\la } ^I ( \psi )  + d  E_{\la } ^I ( \psi ) ^{ \frac{2}{\g}}, 
\eeq
where $ p_0 $ is as in {\bf(A2)} and $ \g \in (0, 1]$ is the exponent appearing in {\bf (B2)}.
\end{Lemma}

{\it Proof. } We have $ \big| \, |\psi |^2 - 1 \big| = \big| \,  |\psi | - 1\big| \cdot  \big| \,  |\psi | + 1 \big| \ges   \big| \,  |\psi | - 1\big|$. 
If $ |\psi |^2 - 1 \in L^2 ( \R\times [0,1]) , $ it is obvious that $ |\psi | - 1 \in L^2 ( \R\times [0,1]) . $

Conversely, assume  that $ |\psi |- 1 \in L^2 ( \R\times [0,1])  $  and $ \nabla \psi \in L^2 ( \R\times [0,1])  $. 
Since $ \big| \nabla |\psi | \, \big| \les |\nabla \psi | $ almost everywhere, we infer that $ |\psi |- 1 \in H^1 ( \R\times [0,1])  $.
Then the Sobolev embedding implies that $ |\psi |- 1 \in L^p  ( \R\times [0,1])  $ for any $ p \in [2, \infty)$. 
We have 
$$
\big| \, |\psi |^2 - 1 \big| \les 5 \big| \, |\psi | - 1 \big| \1_{ \{|\psi | \les 4 \} } + \frac 53 \big| \, |\psi | - 1 \big| ^2 \1_{ \{|\psi | > 4 \} }
$$
and we infer that $  |\psi |^2 - 1 \in L^2 ( \R \times [0,1]). $

We will repeatedly use the following simple observation. 
Let  $ I \subset \R$ be an interval of length greater than or equal to $1$. 
Proceeding as in \cite{brezis}, Section 9.2 we may use four successive "mirror symmetries" to extend  any function 
$ u \in H^1( I \times (0, 1)) $ to a function $ \tilde{u} \in H^1(   \Om_I)$, where $ \Om _I $ is a domain containing 
$  ( I + ( -1, 1) ) \times ( -1, 2))$.  
Then we choose a cut-off function $ \chi \in C_c ^{\infty} ( \Om _I)$ such that $ 0 \les \chi \les 1$, $ \chi = 1 $ on $ I \times [0, 1]$ and 
$ \nabla \chi $ is bounded independently of $I$. Denoting $ P(u) = \chi \tilde{u}$, we see that 
$ P(u) \in H^1 ( \R^2)$ and $ \| Pu \|_{H^1( \R^2 ) } \les C \| u \|_{H^1 ( I \times (0, 1))}$, where $C$ is independent of $I$. 
Using the Sobolev embedding in $ \R^2$ we see that for any $ p \in [2, \infty ) $ there exists $ C_p > 0 $ depending only on $p$ such that 
$$
\| u \|_{L^p ( I \times (0, 1) )} \les \| P(u) \|_{L^p ( \R^2 ) } \les C_p \| P(u) \|_{H^1( \R^2 )} \les C C_p \| u \|_{H^1( I \times (0 , 1))}. 
$$

Assume that $ \psi \in H_{loc}^1 ( \R \times (0, 1) ) $ and $ E_{GL } ( \psi ) < \infty$. 
From the above arguments it follows that $ |\psi | - 1 \in H^1 ( \R \times (0, 1)) $ and 
$ \| |\psi | - 1 \|_{H^1( I \times (0,1 )) }^2 \les C E_{GL} ^I  ( \psi ) $ for any interval $I $ of length at least $1$. 
By {\bf (A1)} and {\bf (A2)} there exist $ C_1, C_2 > 0 $ such that 
$$
V( s^2 ) \les C_1 \left( s^2-1 \right)^2 \quad \mbox{ if } 0 \les s \les 2, 
\qquad \mbox{ and } \qquad
V( s^2 ) \les C_2 (s-1)^{ 2 p_0 + 2 } \quad \mbox{ if } s > 2. 
$$
Using the Sobolev embedding we have 
$$
\begin{array}{l}
\ds \ii_ { I \times (0, 1) } V ( |\psi |^2 ) \, dx 
\les \ii_ { I \times (0, 1) } C_1 \left( \, |\psi |^2 - 1 \right)^2 + C_2 \big| \, |\psi | - 1 \big|^{ 2 p_0 + 2} \, dx 
\\
\\
\ds \les C_1 E_{GL }^I ( \psi ) + C_2 C_{2 p_0 + 2} \| \,  |\psi | - 1 \|_{H^1 ( I \times (0, 1))} ^{ 2 p_0 + 2 }
\les C_1 E_{GL }^I ( \psi ) + C_2 '  E_{GL }^I ( \psi )^{ p_0 + 1}. 
\end{array}
$$
The first estimate in (\ref{Mutual-Bounds}) is thus proven. 
If $I = \R$ we see that $ V( |\psi |^2) \in L^1 ( \R  \times (0, 1))$. 

Conversely, assume that {\bf(B2)} holds,  $ \nabla \psi \in L^2 ( I \times [0,1]) $ and $ V ( |\psi |^2 ) \in L^1( I \times [0,1])$. 
Without loss of generality we may assume that $ 0 < \g <  1$. 
Take an increasing,    concave, $C^1$  function $ G : \R \lra \R$ such that  $ G(s ) = s $ if $ s \les 4$ and $ G( s ) = s^{\g}$ for  $s \ges s_1$, where $ s_1 > 4$, and $ 0 < G'\les 1$. 
By assumptions {\bf (A1)} and {\bf (B2)} there is $C > 0 $ such that $ G( |\psi | - 1) ^2 \les C V( |\psi |^2) $, 
hence $ G( |\psi | - 1) \in L^2 (  I \times [0,1])$. 

We claim that $ G ( |\psi | - 1) \in H^1 ( I \times (0,1))$ and $ \nabla \left( G ( |\psi| -1 ) \right) = G'( |\psi | - 1) \nabla (|\psi|)$ almost everywhere. 
Indeed, let $ u_n = \min( |\psi | -1, n)$. We have $ |\nabla u _n | \les |\nabla \psi |$ a.e., 
and there is $ C_n > 0$ such that $ |u_n|^2 \les C_n V( |\psi|^2)$, hence $ u_n \in H^1( I \times (0,1))$. 
By Proposition 9.5 p. 270 in \cite{brezis} we have $ G( u_n ) \in H^1 ( I \times (0,1)) $ and 
$  \nabla  (G( u_n ) ) = G'(u_n ) \nabla u_n = G'( u_n ) \nabla( |\psi |) \1_{ \{|\psi | \les n+1 \} }  $ a.e. 
The claim follows by letting $ n \lra \infty$ and using the Dominated Convergence Theorem.

It is obvious that $ \big| \nabla ( G( |\psi| - 1 )) \big| \les | \nabla  \psi |$ a.e. and we conclude that 
$ \| G( |\psi| - 1 ) \|_{H^1 ( I \times (0,1)) } ^2 \les C_{\la}  E_{\la } ^{I} ( \psi )$. 
By the Sobolev embedding we have $ \| G( |\psi| - 1 ) \|_{L^p ( I \times (0,1)) } \les C(p, \la )  E_{\la } ^{I} ( \psi )^{\frac{p}{2}} $ for any $ p \in [2, \infty)$. 
We have 
$$ (|\psi |^2 - 1 )^2 \les 25 G( |\psi | -1 ) ^2 \1_{\{ |\psi | \les 4 \} } + C G( |\psi | - 1) ^{\frac{4}{\g} }\1_{\{ |\psi | > 4 \} } $$ 
and we infer that 
$$
\ii_{I \times (0, 1)} (|\psi |^2 - 1 )^2 \les 25 \| G( |\psi | -1 )\|_{L^2 ( I \times (0,1))}^2 + C 
\| G( |\psi | -1 )\|_{L^{\frac{4}{\g}} ( I \times (0,1))}^{\frac{4}{\g}}
\les C_1 E_{\la } ^{I} ( \psi ) + C_2 E_{\la } ^{I} ( \psi )^{\frac{2}{\g}}. 
$$
This gives the second estimate in (\ref{Mutual-Bounds}). 
\hfill
$\Box$

\medskip

\begin{Lemma}
\label{PG-bis} Assume that {\bf (A1), (A2)} and {\bf (B2)} hold.  Then: 

\medskip

i) For any $ \psi \in \Ep $ and any $ v \in H_{per}^1 $ we have $ \psi + v \in \Ep$. 

\medskip

ii) Let $ \psi \in \Ep$. Then for almost all $ y \in \R$ the mapping $ \psi (\cdot, y )$ belongs to $ \Er$ and for almost all $ x \in \R$ 
the mapping $ \psi ( x, \cdot )$ belongs to $H^1( (0,1))$. 

\medskip

iii) For a given $ \psi \in \Ep,$  define $ \breve{\psi}(x) = \int_0^1 \psi( x, y ) \, dy $ and 
$ v_{\psi } ( x, y ) = \psi ( x, y ) - \breve{\psi}(x)$. Then we have 
$ \breve{\psi} \in \Er, $ $ v_{\psi } \in H_{per}^1 $ and 
$$
(\breve{\psi}) '(x) = \int_0^1 \frac{\p \psi }{\p x } (x, y ) \, dy, \qquad 
\frac{ \p v_{\psi}}{\p y } = \frac{\p \psi}{\p y } 
\quad \mbox{ almost everywhere.} 
$$
For any interval $ I \subset \R$ of length  greater than or equal to  $1$ we have
\beq
\label{2.5}
\| v_{\psi }\|_{H^1 ( I \times (0,1) )} ^2 \les 2 \| \nabla \psi \|_{L^2 ( I \times (0, 1))}^2 
\quad 
\mbox{and} 
\quad 
E_{GL}^{1, I} ( \breve{\psi} ) \les C \! \left( E_{GL}^I ( \psi ) + E_{GL}^I ( \psi )^2 \right) .
\eeq

\medskip

iv) For any $ \psi \in \Ep $ there exist a real-valued function $ \ph \in C^{\infty} ( \R)$ satisfying
$ \ph ^{(k)} \in L^2( \R)$ for any $ k \in \N^*$ and $ w \in H_{per}^1$ such that $ \psi (x, y) = e^{i \ph (x) } + w(x,y) . $

\end{Lemma}

{\it Proof. } 
(i) It is clear that $ \psi + v \in H_{loc}^1 ( \R^2 ) $ and $ \nabla ( \psi + v ) \in L^2 ( \R \times [0,1])$. 
We only need to show that $| \psi + v |^2 - 1 \in L^2 ( \R \times [0,1])$. 
Recall that $|\psi | - 1 , v \in H_{per} ^1 \subset L^p ( \R \times [0,1])$ for any $ p \in [2, \infty)$ by the Sobolev embedding. 
We have
\beq
\label{2.6}
| \psi + v |^2 - 1 = \left( |\psi |^2 -1 \right) + 2 \langle \psi, v \rangle + | v|^2. 
\eeq
It is clear that $ |\psi |^2 -1 \in  L^2 ( \R \times [0,1])$ because $ \psi \in \Ep $ and 
$|v|^2 \in  L^2 ( \R \times [0,1])$ because $ v \in  L^4 ( \R \times [0,1])$.
We have also 
\beq
\label{2.7}
|   \langle \psi, v \rangle  |\les | \psi | \cdot | v | \les \big| \, |\psi | -1 \big| \cdot |v| + |v|. 
\eeq
The  last function is in  $ L^2 ( \R \times [0,1])$ because $  |\psi | -1 , v \in L^4  ( \R \times [0,1])$
and $ v \in L^2( \R \times [0,1])$.

\medskip

(ii) is a consequence of Theorem 2 p. 164 in \cite{EG} and of Fubini's Theorem. 

\medskip

(iii) By Fubini's Theorem, $ \breve{\psi} $ is measurable. By the Cauchy-Schwarz inequality we have 
$ \big| \breve{\psi} ( x ) \big| ^2 \les \ii_0^1 |\psi ( x, y )|^2 \, dy $ and we infer that $ \breve{\psi } \in L_{loc}^2 ( \R)$. 

Let $ g ( x ) =  \ii_0 ^1 \frac{\p \psi }{\p x } ( x, y ) \, dy . $
As above, using  the Cauchy-Schwarz inequality we find
$ |g ( x )|^2 \les  \ii_0^1 \big| \frac{\p \psi }{\p x } ( x, y ) \big|^2 \, dy $ and we infer that 
$ g \in L^2 ( \R)$ and $ \| g \|_{L^2( I) } \les \big\|  \frac{\p \psi }{\p x } \big\|_{L^2 ( I \times (0,1) )} $ for any interval $ I \subset \R$. 
For any $ \phi \in C_c^{\infty}( \R)$ we have 
$$
\begin{array}{l}
\ds \ii_{\R}  \breve{\psi} ( x ) \phi '(x) \, dx = \ii_{\R} \left( \ii_0^1 \psi ( x, y ) \phi ' ( x ) \, dy \right) dx 
\\
\\
\ds = \ii_0^1 \left( \ii_{\R} \psi ( x, y ) \phi '( x) \, dx \right) dy \quad \mbox{ by Fubini because } \psi ( x,y) \phi ' (x) \in L^1 ( \R \times [0,1])
\\
\\
\ds = -  \ii_0^1 \left( \ii_{\R} \frac{\p \psi}{\p x} ( x, y ) \phi ( x) \, dx \right) dy \quad \mbox{ because } \psi ( \cdot, y ) \in H_{loc}^1 ( \R) \mbox{ for a.e. } y \in [0,1]
\\
\\
\ds = - \ii_{\R} \left( \ii_0^1 \frac{\p \psi}{\p x} ( x, y ) \, dy \right)  \ph ( x) \, dx \quad \mbox{ by Fubini again because }
\frac{\p \psi}{\p x } \phi \in L^1 ( \R \times [0,1])
\\
\\
= - \ds \ii_{\R} g( x ) \phi ( x ) \, dx. 
\end{array}
$$
We conclude that $ \breve{\psi } \in H_{loc}^1 ( \R)$ and $ \left( \breve{\psi } \right) ' = g$. 

It is clear that $ v_{\psi}$ is $1-$ periodic with respect to the second variable and $ \frac{\p v_{\psi}}{\p y } = \frac{ \p \psi }{\p y }$. 
For almost every $ x \in \R$ we have $ v_{\psi }( x, \cdot ) = \psi ( x, \cdot ) - \breve{\psi} ( x ) \in H^1 ((0,1))$ and 
$ \ii_0^1  v_{\psi }( x, y ) \, dy = 0 $.
For any such $x$, using the Poincar\'e-Wirtinger inequality (see \cite{brezis} p. 233) we get  
$$
\ii_0^1 \big| v_{\psi } ( x, y ) \big|^2 \, dy \les \| v_{\psi} ( x, \cdot ) \|_{L^{\infty}(0,1)}^2 \les 
\Big\| \frac{\p v_{\psi }}{\p y } ( x, \cdot ) \Big\|_{L^1 ( (0,1))}^2 \les \ii_0^1 \Big| \frac{\p \psi }{\p y } ( x, y ) \Big|^2 \, dy.
$$
Integrating with respect to $x$ we infer that $ v_{\psi } \in L^2 ( \R \times [0,1])$ and we have 
$ \| v_{\psi } \|_{L^2 ( I \times [0,1] )} \les \big\| \frac{\p \psi }{\p y }  \big\| _{L^2 ( I \times [0,1] )}$
for any interval $ I \subset \R$.

We have 
$$ \frac{\p v_{\psi}}{\p x } (x, y ) = \frac{ \p \psi }{\p x } (x, y) -  \left( \breve{\psi } \right) ' (x) 
= \frac{ \p \psi }{\p x } (x, y) - \ii_0^1  \frac{ \p \psi }{\p x } (x, y ) \, dy $$
and we see that for almost every  $ x \in \R$ there holds 
$$
\ii_0^1 \Big|  \frac{\p v_{\psi}}{\p x } (x, y )  \Big| ^2 \, dy = \ii_0^1 \Big|  \frac{\p \psi}{\p x } (x, y )  \Big| ^2 \, dy 
- \left(\ii_0 ^1 \frac{\p \psi}{\p x } (x, y ) \, dy \right)^2
\les  \ii_0^1 \Big|  \frac{\p \psi}{\p x } (x, y )  \Big| ^2 \, dy  . 
$$
The above estimates imply that  $ v_{\psi } \in H_{per}^1$ and the first estimate in (\ref{2.5}) holds. 

We have $ \breve{\psi } = \psi - v_{\psi }$. Then using (\ref{2.6}), (\ref{2.7})
and the Sobolev inequality  
$$ \| \, |\psi | - 1 \| _{L^4 ( I \times [0,1]) } \les C \| \, |\psi | - 1 \| _{H^1 ( I \times [0,1]) } \les C E_{GL}^I ( \psi )^{\frac 12} $$ 
as well as the similar estimate for $ v_{\psi}$ we get the second estimate in (\ref{2.5}). 

\medskip

(iv) By Lemma \ref{PG} there exist $ \ph \in C^{\infty } ( \R)$ such that $ \ph ^{(k)} \in L^2 ( \R)$ for any $ k \ges 1 $
 and $ w_1 \in H^1( \R)$ such that $ \breve{\psi } = e^{ i \ph } + w_1$. 
 Letting $ w ( x , y ) = v_{\psi }( x, y ) + w_1 ( x)$ we see that $ w \in H_{per}^1 $ and (iv) holds.
 \hfill
 $\Box$

\section{The momentum} 

\subsection{Definition of the momentum on $ \Er$}

From a mathematical point of view, the momentum should be a functional defined on $ \Er$ such that for any $ \psi \in \Er $ and for  any $ v \in H^1( \R)$, 
\beq
\label{3.1}
\lim_{t \ra 0 } \frac{P( \psi + t v ) - P(\psi) }{t } = 2 \int_{\R} \langle i \psi ', v \rangle \, dx. 
\eeq. 

Notice that functions in $ \Er$ may oscillate at infinity. 
One can introduce a distance and define a manifold structure on $ \Er$, see \cite{PG2}. 
The tangent space of $ \Er $ at $ \psi $ contains $H^1( \R)$, but is larger than $ H^1(\R)$ (see \cite{PG2} p. 140). 
We require (\ref{3.1}) to hold only for $ v \in H^1( \R)$, hence condition (\ref{3.1}) is weaker than G\^ateaux differentiability and this allows  some flexibility in the choice of the  definition of the momentum. 

We were inspired by the definitions of the momentum in higher space dimensions given in \cite{M10, CM}. 
We recall that the energy space associated to eq. (\ref{1.1}) in $ \R^N $ with "boundary condition" $|\psi | \lra 1 $ as $ |x |\lra \infty $  is 
$$ \Eo ( \R^N) = \{ \psi \in H_{loc}^1( \R^N) \; \mid \; \nabla \psi \in L^2( \R^N) \mbox{ and } |\psi | - 1 \in L^2( \R^N) \}$$
(see (1.11) p. 154 in \cite{CM}). 
The momentum with respect to the $x_1-$direction is a functional $ P : \Eo ( \R^N) \lra \R$ satisfying 
$$
\lim_{t \ra 0 } \frac{P( \psi + t v ) - P(\psi) }{t } = 2\int_{\R^N} \langle i \frac{\p \psi }{\p x_1}, v \rangle \, dx
\quad \mbox{ for any } v \in H^1( \R^N) \mbox{ with compact support.}
$$
Formally we should take $ P(\psi ) = \ds \int_{\R^N} \langle  i \frac{\p \psi }{\p x_1}, \psi \rangle \, dx$, 
except that for an arbitrary function $ \psi \in  \Eo ( \R^N) $,   the function $ \langle  i \frac{\p \psi }{\p x_1}, \psi \rangle $ 
is not necessarily in $L^1( \R^N)$. 
However, it has been shown in \cite{M10, CM} that for any $ \psi \in \Eo ( \R^N)$ there exist 
$ f \in L^1( \R^N)$ and $ g \in \dot{H}^1( \R^N)$ such that 
$$
\langle  i \frac{\p \psi }{\p x_1}, \psi \rangle  = f + \p_{x_1} g. 
$$
Obviously, $f$ and $g$ are not unique. 
However, if $ h \in \dot{H}^1( \R^N)$ and $ \p_{x_1} h \in L^1( \R^N)$, then  Lemma 2.3 p. 122 in \cite{M10} implies that 
 necessarily $ \int_{\R^N} \p_{x_1} h \, dx = 0 $. 
This allows to define unambiguously the momentum by $ P( \psi ) = \int_{\R^N} f(x) \, dx. $

The situation is different in space dimension one. 
If $ h \in \dot{H}^1( \R)$ and $ h'\in L^1( \R)$, the integral $ \int_{- \infty } ^{\infty} h'(x) \, dx $ can take any value. 
This is due to the fact that functions in $\dot{H}^1( \R)$   may have different limits or  may oscillate at $\pm \infty$. 
To give an example, let $ \al \in (\frac 12, 1)$ and consider $ \chi \in C^{\infty}( \R, \R)$ such that $ \chi' ( x ) = \frac{1}{| x|^{\al} }$ on $ (-\infty, -1] \cup [ 1 , \infty)$. 
Then $ \chi( x ) = \frac{1}{1 - \al } x^{ 1 - \al } + C_1 $ on $[1, \infty)$ and a similar formula holds on $ (-\infty, -1]$. 
We have $ \chi \in \dot{H}^1( \R)$,  $ e^{ i \chi } \in \Er $ and $ \chi ' \not\in L^1( \R)$. 
For $a, b \in \R$, $ a < b $, let 
$$
 \chi_{a, b }( x ) = \left\{ \begin{array}{ll} \chi( a ) & \mbox{ if } x < a, \\
\chi( x ) & \mbox{ if } x \in [a, b], \\
\chi( b ) & \mbox{ if } x > b.
\end{array}
\right. 
$$
We have $ \chi_{a, b } \in \dot{H}^1( \R)$,  $ e^{ i \chi_{a,b} } \in \Er $, $ \chi_{a,b} ' \in L^1( \R)$
 and $ \int_{- \infty }^{\infty} \chi_{a,b} ' (x) \, dx $ may take any  value in $ \R$ as $a $ and $ b $ vary.

Assume that $ \psi \in H_{loc}^1 ( \R^N)$ can be written in the form  $ \psi = e^{ i \ph } + w$, where  $ \ph $ is real-valued and 
 $ \ph ,  w \in H_{loc}^1( \R^N)$. 
A simple computation gives 
\beq
\label{analog}
\langle i \frac{ \p \psi}{\p x_1 } , \psi \rangle = - \frac{\p \ph }{\p x_1 } + \frac{\p }{ \p x_1 } \left( \langle iw, e^{i \ph } \rangle \right) 
- 2 \langle \frac{ \p \ph }{\p x_1 } e^{ i \ph } , w \rangle  + \langle i \frac{ \p w }{ \p x_1 }, w \rangle. 
\eeq
In space dimension $ N = 2$, a variant of Lemma \ref{PG}  asserts that for any $ \psi \in \Eo ( \R^2)$ there exist 
 a real-valued function $ \ph \in \dot{H}^1( \R^2)$ and $w \in H^1( \R^2)$ such that $ \psi = e^{ i \ph } + w$ (see Theorem 1.8 p. 134 in \cite{PG2}). 
 Then we have $\langle iw, e^{i \ph } \rangle \in H^1( \R^2) $ (see the proof of Lemma 2.1  p. 158 in \cite{CM}). 
According to Lemma 2.3 p. 122 in \cite{M10} we must have 
$ \int_{\R^2}   - \frac{\p \ph }{\p x } + \frac{\p }{ \p x } \left( \langle iw, e^{i \ph } \rangle \right) \, dx \, dy = 0 $ whenever 
$  - \frac{\p \ph }{\p x } + \frac{\p }{ \p x } \left( \langle iw, e^{i \ph } \rangle \right) \in L^1( \R^2)$. 
This observation enables to define unambiguously the momentum on $ \Eo ( \R^2) $ by 
$$
Q( \psi ) = \int_{\R^2} - 2 \langle \frac{ \p \ph }{\p x } e^{ i \ph } , w \rangle  + \langle i \frac{ \p w }{ \p x }, w \rangle \, dx \, dy
$$
for any function $ \psi = e^{ i \ph } + w$, where  $ \ph $ and $ w $ are as above. 
The integrand belongs to $L^1( \R^2)$ by the Cauchy-Schwarz inequality, and the value of the integral does not depend on the choice of the functions  $ \ph $ and $ w $ satisfying the above properties. 

The situation is more complicated in space dimension  $ N = 1$ because the integral of a derivative of a function in $ \dot{H}^1( \R)$, when it exists,  does not necessarily vanish. 
Nevertheless, we will use an analogy to the two-dimensional case. 
More precisely, for any real-valued function $ \ph \in \dot{H}^1( \R)$ and for any $ w \in H^1( \R)$ we define 
\beq
\label{p}
p(\ph, w) = \int_{\R} - 2 \langle \ph ' e^{ i \ph } , w \rangle  + \langle i  w ', w \rangle \, dx. 
\eeq
Notice that  the integrand is in $ L^1( \R)$ by the Cauchy-Schwarz inequality  and   consequently 
$p( \ph, w)$ is well-defined. 

Assume that $ \psi \in \Er $ can be written as $ \psi = e^{ i \ph _1 } + w_1 = e^{ i \ph _2 } + w_2$, where $ (\ph_1, w_1)$ and $ (\ph_2, w_2)$ are as in Lemma \ref{PG}. 
By (\ref{analog}) we have 
$$
- 2 \langle \ph_j ' e^{ i \ph_j} , w_ j \rangle  + \langle i  w_j ', w_j \rangle = \langle i \psi ' , \psi \rangle + \ph _j ' 
- \left(\langle  i w_j, e^{i \ph _j } \rangle \right) '
$$ 
for $ j = 1, 2 $, 
therefore 
\beq
\label{3.4}
\begin{array}{l}
(- 2 \langle \ph _2' e^{ i \ph_2 } , w_2 \rangle  + \langle i  w_2', w_2 \rangle ) - 
( - 2 \langle \ph_1 ' e^{ i \ph_1 } , w_1\rangle  + \langle i  w_1 ', w_1 \rangle  ) 
\\
\\
= \ph_2 '- \ph _1 '  - \left(\langle  i w_2, e^{i \ph _2 } \rangle \right) ' + \left(\langle  i w_1, e^{i \ph _1 } \rangle \right) '.
\end{array}
\eeq
By Lemma \ref{PG} there exist $ k_+, k_- \in \Z$ such that $ \ph_2 - \ph _1 - 2  k_{-}\pi \in L^2 ( ( - \infty, 0 ]) $ and
 $ \ph_2 - \ph _1 - 2  k_{+}\pi \in L^2 ( [0, \infty )) $. Then we have  $ \ph_2 - \ph _1 - 2  k_{-}\pi \in H^1 ( ( - \infty, 0 )) $ and
 $ \ph_2 - \ph _1 - 2  k_{+}\pi \in H^1 ( (0, \infty )) $, hence 
 $  \ph_2 - \ph _1 \lra 2  k_{-}\pi  $ as $ x \lra - \infty $ and  $  \ph_2 - \ph _1 \lra  2  k_{+}\pi   $ as $ x \lra  \infty $.

It is easy to see that $ \langle  i w_j, e^{i \ph _j } \rangle  \in H^1( \R)$ for $ j = 1, 2$. 
For any function $ f \in H^1( \R)$ we have $ f(s) \lra 0 $ as $ s \lra \pm \infty $ and $ \int_{-R}^R f'(t ) \, dt = f( R) - f( - R) \lra 0 $ as $ R \lra \infty$. 
Using (\ref{p}), Lebesgue's dominated convergence theorem, then (\ref{3.4}) we get 
\beq
\label{quantif}
\begin{array}{l}
p( \ph_2, w_2 ) - p( \ph_1, w_1) 
\\
\\
= \ds  \lim_{R \ra \infty } \int_{-R}^R (- 2 \langle \ph _2' e^{ i \ph_2 } , w_2 \rangle  + \langle i  w_2', w_2 \rangle ) - 
( - 2 \langle \ph_1 ' e^{ i \ph_1 } , w_1\rangle  + \langle i  w_1 ', w_1 \rangle  ) \, dx 
\\
\\
= \ds \lim_{R \ra \infty } \int_{-R}^R   \ph_2 '- \ph _1 '  - \left(\langle  i w_2, e^{i \ph _2 } \rangle \right) ' + \left(\langle  i w_1, e^{i \ph _1 } \rangle \right) ' \, dx 
= 2 \pi ( k_+ - k_-). 
\end{array}
\eeq

Let  $ \psi = e^{ i \ph  } + w \in \Er$, where $ \ph $ and $  w $ are as in Lemma \ref{PG}. 
Let  $ k \in \Z$. 
Consider a real-valued function $ \chi \in C^{\infty } ( \R)$ such that $ \chi = 0 $ on $ (- \infty, 0 ]$ and $ \chi = 1 $ on $[1, \infty)$. 
Define $ \tilde{\ph } = \ph  + 2 k \chi $ and $ \tilde{w} = w + e^{ i \ph } - e^{ i \tilde{\ph}}$. 
It is easily seen that $ \psi =  e^{ i \tilde{\ph } } + \tilde{w} $,  where $\tilde{ \ph }$ and $\tilde{  w }$ also satisfy the conclusion of  Lemma \ref{PG}, and the above computation shows that 
$ p( \tilde{\ph}, \tilde{w}) - p( \ph, w ) = 2 k \pi$. 

We conclude given any $ \psi \in \Er$, it can be written as $ \psi = e^{ i \ph  } + w \in \Er$, where $ \ph $ and $  w $ are as in Lemma \ref{PG}, but  the quantity $ p( \ph, w)$ is well-defined only modulo $ 2 \pi \Z$. 
We denote by $ \pr{\cdot }$ the projection of $ \R$ onto $ \R / 2 \pi \Z$, namely 
$ \pr{ x } = \{ x + 2 k \pi \; | \; k \in \Z\}$. 

For $ \pr{p} \in {\R} / {2 \pi \Z}$  we denote $ | \pr{p} | = \inf \{ |p'| \; \mid \; p'\in \pr{p} \}$. 
Notice that $ | \pr{p} | \in [0, \pi]$ and the infimum is achieved. 
It is obvious that $|\pr{ p_1 + p_2 } | \les | \pr{p_1}| + | \pr{p_2} |.$
The distance between the classes $\pr{p_1} $ and $\pr{p_2}$ is 
$$ dist( \pr{p_1},  \pr{p_2}) :=  | \pr{p_1} - \pr{p_2} | = |\pr{ p _ 1 - p_2} |
= \min \{ | p_1 ' - p_2 '| \, \mid \;  p_1 ' \in  \pr{p_1}, p_2 ' \in  \pr{p_2} \}.$$
Endowed with this distance, ${\R} / {2 \pi \Z}$ is a compact group. 
We have $\pr{p_n } \lra \pr{p}$ if and only if there exists a sequence of integers $ (\ell_n)_{n \ges 1} \subset \Z$ such that $ p_n + 2 \pi \ell _ n \lra p$. 

\begin{Definition}
\label{mom1}
Given any $ \psi \in \Er$, the momentum of $ \psi $ is 
$$
\PR(\psi) = \pr{p ( \ph, w )}, 
$$
where $ \psi =  e^{ i \ph  } + w$ and $ \ph $ and $ w $ are as in Lemma \ref{PG} (i). 
We call a {\rm valuation} of the momentum of $ \psi $ any number in  the set $ \PR(\psi)$, 
and the {\rm canonical valuation} the only number in the set $ [0, 2 \pi) \cap \PR(\psi) $.
\end{Definition}

It follows from the above discussion that a number $ p \in \R$ is a valuation of the momentum of a mapping $ \psi \in \Er $ if and only if there exist  $ \ph \in \dot{H}^1 ( \R, \R) $ and $  w \in H^1( \R, \C)$ 
 such that $  \psi = e^{ i \ph  } + w $ and
$ p( \ph, w ) = p $. 

For any $ \psi = e^{  \ph } + w \in \Er$ and any $ \al \in \R$ we have $ e^{ i \al } \psi = e^{ i ( \ph + \al )} + e^{ i \al } w$. It is clear that 
$p ( \ph, w ) = p ( \ph + \al, e^{ i \al } w )$, and consequently we have $ \PR(e^{ i \al } \psi) = \PR(\psi).$

For all  $ \ph \in \dot{H}^1 ( \R, \R )$ we have $ e^{ i \ph } \in \Er $ and $ E^1 ( e^{ i \ph }) = \| \ph '\|_{L^2( \R)}^2$. 
Definition \ref{mom1} gives  $ \PR ( e^{ i \ph }) = \pr{ p ( \ph, 0 ) } = \pr{0 }$. 
However,  we have $ \langle i \left( e ^{ i \ph } \right) ', e^{ i \ph } \rangle = - \ph '$ and $ \ph '$ does not necessarily belong to $ L^1 ( \R)$; when it does, $ \ii_ {\R} \ph '( x ) \, dx $ can take any value. 

Assume that $ \psi \in \Er $ is constant outside a bounded interval $[a, b]$, say $ \psi (x ) = e^{ i \al _1 }$ on $ ( - \infty, a ] $ and 
$ \psi (x ) = e^{ i \al _2 }$ on $[b, \infty)$, where $ \al_1, \al _2 \in \R$. 
Consider any function $ \ph \in C^{ \infty }( \R)$ such that $ \ph = \al _1 $ on $ ( - \infty, a] $ and $ \ph  = \al _2 $ on $[b, \infty)$. 
Let $ w = \psi - e^{ i \ph}$. 
Then $ w \in H^1( \R)$, $ \mbox{supp}(w) \subset [a, b ]$ and using (\ref{analog}) and (\ref{p}) we get 
$$
p ( \ph, w ) = \al _2 - \al _ 1 + \ii_a ^b \langle i \psi ', \psi \rangle \, dx. 
$$
In particular, if $ \al _1 = \al _2 $ we see that a valuation of the momentum of $ \psi $ is $ \ii_ {\R}   \langle i \psi ', \psi \rangle \, dx. $

Given any $ \psi = e^{ i \ph  } + w \in \Er $ and  any $ v \in H^1( \R)$, we have $ \psi + v = e^{i \ph } + ( w + v ) $, hence a valuation of the momentum of $ \psi + v $ is $ p ( \ph, w + v )$. It is obvious that 
$$
 p ( \ph, w + v ) - p ( \ph, w ) = \ii_{\R}  - 2 \langle \ph 'e^{ i \ph }, v \rangle + 2 \langle i w', v \rangle + \langle i v', v \rangle \, dx 
 = \ii_{\R} \langle i \left( 2 \psi ' + v ' \right) , v \rangle \, dx
$$
and
$$
\lim_{t \ra 0 } \frac{ p( \ph, w + tv ) - p ( \ph, w) }{t} = 2 \int_{\R} \langle - \ph ' e^{ i \ph}, v \rangle + \langle i w', v \rangle\, dx  =  2 \int_{\R} \langle i \psi ', v \rangle \, dx. 
$$

\begin{remark}
\label{lift1}
\rm
Assume that $ \psi \in \Er $ admits a lifting $ \psi (x) = \rho (x) e^{ i \theta (x)}$ where $ \rho = | \psi |$ and $ \theta \in \dot{H}^1( \R)$. 
By Lemma \ref{MutualBounds} we have $ 1 - \rho \in H^1 ( \R)$.   Let $ w = ( \rho - 1) e^{ i \theta }$. It is easy to see that $ w \in H^1 ( \R)$, 
$ \psi = e^{ i \theta } + w $ and 
$$
- 2 \langle \theta 'e^{ i \theta  } , w \rangle + \langle i w ', w \rangle = ( 1 - \rho ^2 ) \theta '. 
$$
Therefore a valuation of the momentum of $ \psi $ is $ \ds p ( \theta, w ) = \int_{\R} ( 1 - \rho ^2 ) \theta ' \, dx. $
This is in perfect agreement with formula (2.12) p. 123 in \cite{M10} and with formula (2.7) p. 159 in \cite{CM}. 
\end{remark}

\subsection{Definition of the momentum on $ \Ep$}
\label{defmomp}

For any given $ \ph \in \dot{H}^1( \R)$ and any $ w \in H^1 ( \R  \times (0, 1))$ we define 
$$
\ds d [ \ph, w ] ( x, y ) =  - 2 \langle \ph ' (x) e^{ i \ph ( x) } , w (x, y ) \rangle  + \langle i  \frac{ \p w }{ \p x } (x, y ), w (x,y)\rangle .
$$
The Cauchy-Schwarz inequality implies that $ d [ \ph, w ] \in L^1 ( \R  \times (0, 1))$, thus we may define
$$
q ( \ph, w ) = \ds \ii_{\R \times (0, 1)} d [\ph, w ] \, dx \, dy. 
$$

Assume that $ \ph_1, \ph _2\in \dot{H}^1( \R)$,  $ w_1, w_2  \in H^1 ( \R  \times (0, 1)) $  and 
$ e^{ i \ph _1 } + w_ 1 = e^{ i \ph _ 2 } + w _ 2 $ a.e. on $ \R  \times (0, 1)$. 
Let $ h ( x ) = e^{ i \ph _ 2 ( x ) } - e^{ i \ph _ 1 ( x ) }.$ 
Then we have $ h \in \dot{H}^1( \R)$ and $ h = w_1 - w_2 \in H^1 ( \R \times ( 0, 1 ))$. 
By Fubini's theorem it follows that $ h \in L^2 ( \R)$, hence $ h \in H^1 ( \R)$. 

By Theorem 2 p. 164 in \cite{EG}, there is a set $ A \subset (0, 1)$ such that 
$(0, 1) \setminus A$ has zero Lebesgue measure and 
for any $ y \in A $ the mappings $ w_j ( \cdot, y ) $ belong to $ H^1( \R)$ for $ j = 1, 2 $. 
For all $ y\in A $ we have $ \psi (\cdot, y ) = e^{i \ph _j } + w _ j ( \cdot, y ) \in \Er $. 
Using Lemma \ref{PG} (ii), there exist $ k_+, k_- \in \Z$ such that $ \ds \lim_{ x \ra \pm  \infty } \ph_2 ( x ) - \ph _1 ( x ) = 2 \pi k_{\pm}.$
Using (\ref{analog}) we get 
$$
d [ \ph _2, w_2 ] - d[ \ph _1, w_1] = \ph _2 ' - \ph _1 '  - \frac{ \p }{\p x } 
\left( \langle i w_2, e^{ i \ph_2 } \rangle -  \langle i w_1, e^{ i \ph_1 } \rangle    \right) 
= \ph _2 ' - \ph _1 ' + \frac{ \p }{\p x } \left( \langle i h, e ^{ i \ph _1 } + w_2 \rangle \right). 
$$
Proceeding exactly as in (\ref{quantif}) we see that for any $ y \in A$ there holds 
$$
\int _{\R}  d [ \ph_2, w _2] ( x, y ) - d [ \ph_1 , w _1] ( x, y ) \, dx = 2 \pi ( k_+ - k_-) 
$$
and then integrating with respect to $y$ and using Fubini's theorem we get 
$$   q ( \ph_2, w_2 ) -  q ( \ph_1, w_1 )  = 2 \pi ( k_+ - k_-). $$ 

Let   $ \ph \in \dot{H}^1( \R)$, $ w \in H^1 ( \R  \times (0, 1))$
and  $ k \in \Z$  be arbitrary. 
Take a real-valued function $ \chi \in C^{\infty } ( \R)$ such that $ \chi = 0 $ on $ (- \infty, 0 ]$ and $ \chi = 1 $ on $[1, \infty)$
and define $ \tilde{\ph } = \ph  + 2 k \chi $ and $ \tilde{w} = w + e^{ i \ph } - e^{ i \tilde{\ph}}$. 
It is easily seen that $ e^{ i \ph } + w =  e^{ i \tilde{\ph } } + \tilde{w} $ and
$ q( \tilde{\ph}, \tilde{w}) - q( \ph, w ) = 2 k \pi$. 

Given any $ \psi \in \Ep $, by Lemma \ref{PG-bis} (iv) there exist $ \ph \in \dot{H}^1( \R)$ and $ w \in H^1 ( \R  \times (0, 1))$
such that $ \psi (x, y ) = e^{ i \ph (x ) } + w ( x, y )$. The previous discussion shows that the quantity 
 $ q( \ph, w)$ is well-defined  modulo $ 2 \pi \Z$. 
This enables us to give the following 

\begin{Definition}
\label{mom}
Given any $ \psi \in \Ep$, the momentum of $ \psi $ is 
$$
\QR(\psi) = \pr{q ( \ph, w )}, 
$$
where $ \ph \in \dot{H}^1( \R , \R)$ and $ w \in H^1 ( \R  \times (0, 1), \C)$ are such that $ \psi =  e^{ i \ph  } + w$.
A {\rm valuation} of the momentum of $ \psi $ is any number in  the set $ \QR(\psi)$, 
and the {\rm canonical valuation} is the only number in the set $ \QR(\psi) \cap  [0, 2 \pi)$.
\end{Definition}

\begin{remark} \rm 
\label{basic}
i) As in the one-dimensional case, a number $ q \in \R$ is a valuation of the momentum of a mapping $ \psi \in \Ep $ if and only if there exist  $ \ph \in \dot{H}^1 ( \R , \R) $ and $  w\in H_{per}^1 $  such that $  \psi = e^{ i \ph  } + w $ and
$ q( \ph, w ) = q $. 
If $ \psi,$   $ \ph $ and $ w $ are as above, then for almost any $ y \in \R$ we have $ \psi ( \cdot, y ) \in \Er $ and $ w ( \cdot, y ) \in H^1( \R)$,  and then  $ p ( \ph, w( \cdot, y )) $ is a valuation of the momentum of $ \psi(\cdot, y )\in \Er$. 
By Fubini's Theorem we have 
$$
q ( \ph, w ) = \int_0 ^1 p ( \ph, w( \cdot, y ) ) \, dy . 
$$

ii) {\bf Warning!} If $ Q_0 ( \psi )$ and $ P_0 ( \psi ( \cdot, y)) $ are the canonical valuations of the momenta of $ \psi \in \Ep $ and of $ \psi ( \cdot, y ) \in \Er$, respectively, we may have 
$$
Q_0 ( \psi ) \neq \int_0 ^1 P_0 ( \psi(\cdot, y ) ) \, dy . 
$$

iii) For  $ \psi = e^{ i \ph  } + w \in \Ep $ and   $ v \in H^1( \R \times (0, 1))$, 
a valuation of the momentum of $ \psi + v $ is $ q ( \ph, w + v )$ and we have 
\beq
\label{diff-mom}
\begin{array}{l}
\ds q ( \ph , w  + v ) - q ( \ph, w ) = \ii_{\R \times [0,1]} - 2 \langle \ph ' e^{ i \ph } , v \rangle + 2 \langle  i \frac{ \p w }{\p x } , v \rangle + 
  \langle  i \frac{ \p v }{\p x } , v \rangle \, dx \, dy 
\\
\\
= 
\ds \ii_{\R \times [0,1]}  2 \langle i \frac{ \p \psi}{\p x }   , v \rangle +  \langle  i \frac{ \p v }{\p x } , v \rangle \, dx \, dy 
= \ii_{\R \times [0,1]}   \langle i \frac{ \p \psi}{\p x }  +i \frac{ \p ( \psi + v )}{\p x}  , v \rangle  \, dx \, dy .
\end{array}
 \eeq
Using the Cauchy-Schwarz inequality we get 
\beq
\label{diff-mom2}
| q ( \ph , w  + v ) - q ( \ph, w )  | \les 
\| v \|_{L^2 ( \R \times [0,1])} \left( \Big\| \frac{ \p \psi }{\p x } \Big\|_{L^2 ( \R \times [0,1])}  + \Big\| \frac{ \p ( \psi + v )}{\p x } \Big\|_{L^2 ( \R \times [0,1])} \right). 
\eeq
From  (\ref{diff-mom}) we  obtain
\beq
\label{diff-mom3}
\lim_{t \ra 0 } \frac{ q( \ph, w + tv ) - q ( \ph, w) }{t} 
=  2 \int_{\R \times (0,1)} \langle i \frac{ \p \psi}{\p x} , v \rangle \, dx \, dy. 
\eeq
Notice that (\ref{diff-mom})-(\ref{diff-mom3})  are analogous to Lemma 2.5 and Corollary 2.6 p. 123-124 in \cite{M10} and to Lemma 2.3 and Corollary 2.4 p. 159 in \cite{CM}. 
\end{remark}

\begin{Lemma}
\label{lift}
Assume that  $ \psi \in \Ep $  satisfies $ \rho _0 := \ds \inf_{(x, y) \in \R^2 } |\psi ( x, y ) | > 0 .$ Let $ \rho = |\psi |$. 
There exists a real-valued function $ \theta \in H_{per}^1 $ such that $ \psi = \rho e^{ i \theta }$ and a valuation of the momentum of $ \psi $ is $$  \int_{\R \times (0, 1 ) } ( 1 - \rho ^2 )\frac{ \p \theta}{\p x } \, dx \, dy. $$
\end{Lemma}

{\it Proof. } It is well-known that $ \rho \in H_{loc}^1 ( \R^2) $ and $ | \nabla \rho | \les |\nabla \psi |$ almost everywhere. 
The mapping $ \frac{ \psi }{\rho } $ belongs to $ H_{loc}^1 ( \R^2, \mathbb{S}^1)$, hence it admits a lifting, 
in other words there exists 
$ \theta \in H_{loc}^1 ( \R^2, \R)$ such that $  \frac{ \psi }{\rho } = e^{ i \theta}$, or equivalently 
$ \psi = \rho e^{ i \theta}$. 
We have 
$$
| \nabla \psi |^2 = |\nabla \rho |^2  + \rho ^2 |\nabla \theta |^2 \quad \mbox{ almost everywhere}
$$
and we infer that $ \nabla \theta \in L^2 ( \R \times [0, 1])$. 

We claim that $ \theta $ is $1-$periodic with respect to the second variable $y$.
Indeed, since $ \psi $ and $ \rho $ are $1-$periodic with respect to $y$, we infer that for any $ ( x, y ) \in \R^2$ there holds 
$ 1 = \frac{ \psi ( x, y + 1) }{\psi ( x, y )} = e^{i ( \theta ( x, y + 1 ) - \theta ( x, y ))}$. 
The mapping $ ( x, y ) \longmapsto  \theta ( x, y + 1 ) - \theta ( x, y ) $ belongs to $ H_{loc}^1 ( \R^2, \R)$ and takes values in $ 2 \pi \Z$, hence it must be constant. 
We infer that there exists $ k _{\theta } \in \Z$ such that $ \theta ( x, y + 1 ) = \theta ( x, y ) + 2 \pi k_{\theta }$ for any $(x, y ) \in \R^2$. 
For almost any $ x \in \R$ we have $ \theta ( x , \cdot ) \in H^1 (( -1, 2))$ and then using the Cauchy-Schwarz inequality we get
$$
2 \pi |k_{\theta }| = | \theta ( x, 1 ) - \theta ( x, 0 ) | = \bigg| \int_0^1 \frac{ \p \theta}{\p  y } (x, y ) \, dy  \bigg| \les \int_0 ^1 \Big| \frac{ \p \theta}{\p  y } (x, y ) \Big| ^2 \, dy . 
$$
Integrating the above inequality with respect to $x$ and using the fact that  $ \frac{ \p \theta}{\p  y }  \in L^2 ( \R \times (0, 1))$ we see that necessarily  $ k _{\theta } = 0 $ and the claim is proven.

Denote $ \breve{\theta} ( x ) = \int_0 ^1  \theta ( x, y ) \, dy $ and 
$ \theta^{\sharp} (x, y ) = \theta( x, y ) - \breve{\theta} ( x)$. 
Proceeding as in the proof of Lemma \ref{PG-bis} (iii) we see that $ \breve{\theta} '\in L^2 ( \R)$ and $ \theta ^{\sharp} \in H^1( \R \times (0, 1))$.
Let $ w(x,y) = \psi (x, y ) - e^{ i \breve{\theta}(x)}$. 
 We have 
$$
w = ( \rho -1 ) e^{ i \theta } + e^{ i \breve{\theta}} \left( e^{ i \theta^{\sharp} } - 1 \right) .
$$
It is easy to see that $ w \in H^1 ( \R \times (0, 1))$. 
Using (\ref{analog}) we have 
\beq
\begin{array}{l}
\label{3.6}
\ds d[ \breve{\theta}, w ] = \langle i \frac{ \p \psi }{\p x } , \psi \rangle + \breve{\theta} ' 
- \frac{ \p }{\p  x } \left( \langle iw, e^{i \breve{\theta}} \rangle \right) 
\\
\\
\ds = - \rho ^2 \frac{ \p \theta}{\p x } + \breve{\theta} ' - \frac{ \p }{\p  x } \left( \langle i \rho e^{ i \theta} -  i e^{i \breve{\theta}}, e^{i \breve{\theta}} \rangle \right) 
=  - \rho ^2 \frac{ \p \theta}{\p x } + \breve{\theta} '  + \frac{ \p }{\p  x } \left( \rho \sin ( \theta - \breve{\theta} ) \right) 
\\
\\
\ds = ( 1 - \rho ^2 )\frac{ \p \theta}{\p x } + \frac{ \p }{\p  x } \left( - \theta ^{\sharp} +  \rho \sin ( \theta^{\sharp}  ) \right) .
\end{array}
\eeq
By the Cauchy-Schwarz inequality we have $ d[ \breve{\theta}, w ] \in L^1 ( \R \times (0, 1 ) ) $ and
 $ ( 1 - \rho ^2 )\frac{ \p \theta}{\p x } \in L^1 ( \R \times (0, 1 ) ) $, and  then (\ref{3.6}) 
 gives  $ \frac{ \p }{\p  x } \left( - \theta ^{\sharp} +  \rho \sin ( \theta^{\sharp}  ) \right)  \in L^1 ( \R \times (0, 1 ) ) $. 
 
Next we use the following simple observation: whenever $ f \in H^1 ( \R \times (0, 1 ) ) $ satisfies  $ \p _ 1 f \in L^1 ( \R \times (0, 1 ))$ we must have $ \ii_{\R \times (0, 1 ) } \p_1 f ( x, y ) \, dx \, dy = 0 $. 
Indeed, using  Theorem 2 p. 164 in \cite{EG} and Fubini's Theorem we infer that for almost every $ y \in (0, 1 ) $ we have 
$ f ( \cdot, y ) \in H^1 ( \R)$ and
$ \frac{d}{dx } [ f ( \cdot, y ) ]  = \p _1 f ( \cdot, y ) \in L^1 \cap L^2 ( \R)$. For any such $ y $ we have
$$
\int_{\R} \p _ 1 f ( x, y ) \, dx = \lim_{R \ra \infty } \int_{- R }^R \p _1 f ( s, y ) \, ds =  \lim_{R \ra \infty } ( f( R, y ) - f ( - R, y ) ) = 0. 
$$
Integrating with respect to $y$ we get the desired result. 

Since $ \theta ^{\sharp} \in H^1 ( \R \times (0, 1 ) )$ and $ \rho - 1  \in H^1 ( \R \times (0, 1 ) )$ 
 and $ ( \rho -1 ) \nabla \theta \in L^2(\R \times (0, 1)) $,  it is easily seen that 
$  - \theta ^{\sharp} +  \rho \sin ( \theta^{\sharp}  ) =  
- \theta ^{\sharp} + \sin (   \theta ^{\sharp} ) + ( \rho - 1 ) \sin (   \theta ^{\sharp} )  \in H^1 ( \R \times (0, 1 ) )$ and 
then using (\ref{3.6}) and the above observation we infer that 
$$
 q (  \breve{\theta}, w )  = \int_{\R \times (0, 1 ) }  d[ \breve{\theta}, w ]  \, dx \, dy =  \int_{\R \times (0, 1 ) } ( 1 - \rho ^2 )\frac{ \p \theta}{\p x } \, dx \, dy. 
 $$
This completes the proof of Lemma \ref{lift}. 
Notice that this is in agreement with 
 formula (2.12) p. 123 in \cite{M10} and with formula (2.7) p. 159 in \cite{CM}. 
\hfill
$ \Box $

\begin{Lemma}
\label{approx}
i) For any $ \psi \in \Er $ 
 there exist sequences $ (\ph _n)_{n \ges 1} \subset C^{\infty} ( \R, \R) $ and  $ ( w_n )_{n \ges 1} \subset C_c^{\infty}( \R ) $ such that 
 $$  \PR \left( e^{ i \ph _n } + w_n\right) = \PR ( \psi ) \qquad \mbox{ for all } n,  $$ 
 $$  \left( e^{ i \ph _n } + w_n\right) ' \lra  \psi ' \qquad \mbox{ in } L^2 ( \R ), $$ 
 $$ e^{ i \ph _n } + w_n \lra \psi \quad \mbox{ in } H_{loc}^1 (\R ), $$ 
 $$ V \left( \big| e^{ i \ph _n } + w_n \big| ^2 \right) \lra V ( |\psi |^2 )  \;
\mbox{  and  } \; \left( 1 - \big| e^{ i \ph _n } + w_n \big|^2 \right)^2 \lra \left( 1 - |\psi |^2 \right)^2 \qquad \mbox{ in }
 L^1 ( \R ), $$
and there  exist sequences of real numbers $ ( A_n )_{n \ges 1 }$,  $  ( \al_n )_{n \ges 1 }$, $ ( \beta_n )_{n \ges 1 }$
such that $ A_n \lra  \infty$, $ \mbox{ supp}( w_n ) \subset [-A_n, A_n ]$,   $ \ph _n (x) = \al _n $ for $ x \in ( - \infty, -A_n]$, and $ \ph _n (x) = \beta _n $ for $ x \in [A_n , \infty)$. 

\medskip

ii)
For any $ \psi \in \Ep $ 
 there exist sequences $ (\ph _n)_{n \ges 1} \subset C^{\infty} ( \R, \R) $ and  $ ( w_n )_{n \ges 1} \subset C_c^{\infty}( \R \times (0, 1)) $ satisfying
 $$  \QR \left( e^{ i \ph _n } + w_n\right) = \QR ( \psi ) \qquad \mbox{ for all } n,  $$ 
 $$ \nabla \left( e^{ i \ph _n } + w_n\right) \lra \nabla \psi \qquad \mbox{ in } L^2 ( \R \times (0, 1)), $$ 
 $$ e^{ i \ph _n } + w_n \lra \psi \quad \mbox{ in } H_{loc}^1 (\R \times (0, 1)), $$ 
 $$ V \left( \big| e^{ i \ph _n } + w_n \big| ^2 \right) \lra V ( |\psi |^2 )  \;
\mbox{  and  } \; \left( 1 - \big| e^{ i \ph _n } + w_n \big|^2 \right)^2 \lra \left( 1 - |\psi |^2 \right)^2 \quad \mbox{ in }
 L^1 ( \R \times (0, 1)), $$
and there  exist sequences  $ ( A_n )_{n \ges 1 }$, 
 $  ( \al_n )_{n \ges 1 }$, $ ( \beta_n )_{n \ges 1 } \subset \R$
such that $ A_n \lra  \infty$, 
$ \mbox{ supp}( w_n ) \subset [- A_n, A_n ] \times(0, 1)$,   $ \ph _n (x) = \al _n $ for $ x \in ( - \infty, -A_n]$, and $ \ph _n (x) = \beta _n $ for $ x \in [A_n , \infty)$.

\end{Lemma}

{\it Proof. }  We only prove (ii). The proof of (i) is similar.  

If $ \frac{ \p \psi}{\p x } = 0 $ in $ L^2( \R \times (0,1))$, then 
for almost all $ y \in (0, 1)$ the mapping $ \psi  ( \cdot , y )$ belongs to  $ H_{loc}^1 ( \R)$ and $ \psi ( \cdot, y )$ is  constant.
We infer that there exists $ h \in H^1 ( 0, 1 ) $ such that $ \psi ( x, y ) = h( y )$ almost everywhere in $ \R \times (0, 1 )$, and then 
$ \frac{ \p \psi}{\p y } = h ' (y)$. Since $  \frac{ \p \psi}{\p y } \in L^2 ( \R \times (0, 1 ))$, using Fubini's theorem we see that $ h'= 0 $ in $L^2 ( 0, 1)$ and consequently  $ h$ is constant, hence $ \psi $ is constant. 
 In this case it suffices to take $  \ph _n $ a constant function such that $ e^{ i \ph _n } = \psi $, and $ w _n = 0 $. 

If $ \frac{ \p \psi}{\p x } \neq 0 $ in $ L^2( \R \times (0,1))$, there exists $ v \in C_c^{\infty } ( \R \times (0, 1)) $ such that 
$$
\int_{\R \times (0, 1 ) } \langle i \frac{ \p \psi }{\p x }  , v \rangle \, dx \, dy \neq 0 . 
$$
By Lemma \ref{PG-bis} (iv) there exist $ \ph \in C^{\infty} ( \R,\R) $ such that 
$ \ph ^{(k)} \in L^2 ( \R)$ for any $ k \in \N^*$ and $ w \in H_{per}^1 $ such that 
$ \psi = e^{ i \ph } + w $. 
There exists a sequence $ (\tilde{w}_n)_{n \ges 1} \subset C_c^{\infty} ( \R \times (0, 1 ) ) $  such that 
$ \tilde{w}_n \lra w $ in $ H^1 ( \R \times ( 0 , 1 ))$. 
We have 
$$
q ( \ph, \tilde{ w}_n + t v ) = q ( \ph, \tilde{w}_n ) + 2 t \int_{\R \times (0, 1 ) } \langle i \frac{ \p }{\p x } \left( e^{ i \ph  } + \tilde{w}_n \right) , v \rangle \, dx \, dy   + t^2 \int_{\R \times (0, 1 ) } \langle i \frac{ \p v}{\p x } , v \rangle \, dx \, dy .
$$
Since $  q ( \ph, \tilde{w}_n ) \lra q ( \ph, w) $ and 
$$
\int_{\R \times (0, 1 ) } \langle i \frac{ \p }{\p x } \left( e^{ i \ph  } + \tilde{w}_n \right) , v \rangle \, dx \, dy
\lra \int_{\R \times (0, 1 ) } \langle i \frac{ \p \psi  }{\p x } , v \rangle \, dx \, dy \neq 0, 
$$
we infer that there is  a sequence $ t_n \lra 0 $ such that $ q ( \ph, \tilde{ w}_n + t_n v ) = q ( \ph , w )$ for all $ n $ sufficiently large. 
Then we take $ w_n = \tilde{ w}_n + t_n v $. 
If $ \mbox{supp}( w _n ) \subset [a_n, b_n ] \times (0, 1 ) $, we choose $ A_n > 2\max( |a_n|, | b_n |)+ 1 $ such that $ A_n \lra \infty$.
Take $ \chi \in C^{\infty} ( \R, \R)$ such that $ \mbox{ supp}( \chi ) \subset [ - 2, 2]$ and $ \chi = 1 $ on $ [ -1, 1 ]$.  
Take $ \ph _n ( x ) = \ii_0 ^x \chi\left( \frac{2s}{A_n} \right) \ph '( s)\, ds$, so that $ \ph _n = \ph $ on $ [ - A_n /2, A_n/2]$, 
and $ \ph _n $ is constant on  $  (- \infty, - A_n ]$ and on $ [A_n, \infty)$. 
It is easy to see that $ ( \ph_n, w_n)_{n \ges 1}$ satisfy all desired properties. 
\hfill
$\Box$

\begin{Lemma}
\label{continuous}
Let $ \ph \in \dot{H}^1 ( \R)$ and let $ w \in C^2 ( \R^2 ) $ such that there exists $ a > 0 $ satisfying 
$\mbox{supp} ( w ) \subset [-a, a ] \times \R$. Let $ \psi ( x , y ) = e^{ i \ph ( x ) } + w ( x, y )$. 
Then for any $ y_1, y_2 \in \R$, $ y_ 1 < y _2 $ we have 
$$
 p ( \ph , w ( \cdot, y_2)) - p ( \ph , w ( \cdot, y_1)) =  2 \ii_{\R \times [y_1, y_2] }   \langle i \frac{ \p \psi}{\p x }, \frac{ \p \psi}{\p y } \rangle \,  dx \, dy
$$
and consequently 
$$
|  p ( \ph , w ( \cdot, y_2)) - p ( \ph , w ( \cdot, y_1))  | \les 2 \Big\| \frac{ \p \psi }{\p x } \Big\|_{L^2(\R \times [y_1, y_2])}
\Big\| \frac{ \p \psi }{\p y } \Big\|_{L^2(\R \times [y_1, y_2])} .
$$
\end{Lemma}

{\it Proof. } We have $ d[\ph, w ],  \langle i \frac{ \p \psi}{\p x }, \frac{ \p \psi}{\p y } \rangle  \in L^1(\R \times [y_1, y_2])$
and a standard computation gives 
$$
\begin{array}{l} 
\ds p ( \ph , w ( \cdot, y_2)) - p ( \ph , w ( \cdot, y_1)) 
= \ii _{\R} d [ \ph, w ] ( x, y_2 ) -  d [ \ph, w ] ( x, y_1 ) \, dx 
\\
\\
\ds = \ii _{\R} \ii _{y _1}^{y _2} \frac{\p }{ \p y } \left( d [ \ph, w ] ( x, y ) \right) \, dy \, dx 
\\
\\
\ds = \ii_{\R} \ii_{y_1}^{ y_2} \langle - 2 \ph ' ( x ) e^{ i \ph (x)}, \frac{ \p w }{\p y } (x, y )  \rangle 
+ \langle i \frac{ \p^2 w}{\p y \p x } , w \rangle + \langle i \frac{ \p w }{ \p x }, \frac{ \p w }{\p y } \rangle \, dy \, dx 
\\
\\
\ds = \ii_{y_1}^{y_2}  \ii_{\R}  \langle - 2 \ph ' ( x ) e^{ i \ph (x)}, \frac{ \p w }{\p y } (x, y )  \rangle 
+ \langle i \frac{ \p^2 w}{\p y \p x } , w \rangle + \langle i \frac{ \p w }{ \p x }, \frac{ \p w }{\p y } \rangle \, dx \, dy 
\quad \mbox{(Fubini)}
\\
\\
\ds = 2 \ii_{y_1}^{y_2}  \ii_{\R}  \langle -  \ph ' ( x ) e^{ i \ph (x)}, \frac{ \p w }{\p y } (x, y )  \rangle 
 + \langle i \frac{ \p w }{ \p x }, \frac{ \p w }{\p y } \rangle \, dx \, dy 
\quad \mbox{(integration by parts)}
\\
\\
\ds =  2 \ii_{y_1}^{y_2}  \ii_{\R}  \langle i \frac{ \p }{\p x } ( e^{ i \ph } + w ), \frac{ \p w }{ \p y } \rangle \,  dx \, dy 
= 2 \ii_{y_1}^{y_2}  \ii_{\R}  \langle i \frac{ \p \psi}{\p x }, \frac{ \p \psi}{\p y } \rangle \,  dx \, dy .
\end{array}
$$
The second statement follows from the first one and the Cauchy-Schwarz inequality. 
\hfill
$\Box$

\section{One-dimensional traveling waves for (\ref{1.1})}
\label{sect1D}

We consider (\ref{1.1}) in $ \R \times \R$ and we look for traveling waves, namely solutions of the form $ \Phi ( x, t ) = \psi ( x + ct )$. 
The traveling wave profile $ \psi $  satisfies the ordinary differential equation 
\beq
\label{4.1}
i c \psi '+ \psi '' + F ( |\psi |^2 ) \psi = 0 \qquad \mbox{ in } \R. 
\eeq
We will only consider solutions of (\ref{4.1}) in $ \Er$. 
In the case of the Gross-Pitaevskii equation ($F(s ) = 1- s$), 
an extensive study of solutions to (\ref{4.1}) has been carried out in \cite{BGS-survey}, Section 2. 
For more general nolinearities we refer to \cite{C1d} and to \cite{M8} (the latter focuses mainly on the non-existence of supersonic traveling waves).

If assumption {\bf{(A1)}} is satisfied, it has been shown in Theorem 5.1 p. 1099 in \cite{M8} that the only solutions of (\ref{1.1})
with $ c^2 > \vs ^2 = 2 $ are constants. 
It follows from the proof of Theorem 5.1  in \cite{M8} that all traveling waves in $ \Er$ are $ C^2$ functions on $ \R$. 
Let $ \psi \in \Er $ be a solution of (\ref{4.1}) and let $ \varrho = |\psi |^2$. 
Then $ \varrho -1 \in H^1( \R)$ by Lemma \ref{L2.1}  and  it can be shown that $ \varrho $ satisfies the    equation 
\beq 
\label{4.2}
\varrho '' + c^2 ( \varrho - 1) - 2 V ( \varrho) + 2 \varrho F( \varrho) = 0 \qquad \mbox{ in } \R  
\eeq
(for the proof see (5.10) p. 1100 in \cite{M8}).  Multiplying (\ref{4.2}) by $ 2 \varrho '$ and integrating we get 
\beq
\label{4.3} 
(\varrho ')^2 + c ^2 ( \varrho - 1)^2 - 4 \varrho V( \varrho ) = 0 . 
\eeq
Denote $$ g( s, c ) = 4 s V( s ) - c^2 ( s-1)^2. $$
By (\ref{4.3}), for any $ x \in \R$ we must have $ g( \varrho ( x ) , c ) \ges 0 $ and $ \varrho ' ( x ) = \pm \sqrt{ g( \varrho (x), c)}$. 
Since $ g( 0, c ) = - c^2$, we see that for any $ c \neq 0$, solutions of (\ref{4.3}) must stay away from zero. 
This implies that for $ c \neq 0$, any solution $ \psi \in \Er $ of (\ref{4.1}) does not vanish, and therefore has a lifting 
$ \psi (x)= \sqrt{\varrho (x)} e^{i \theta ( x)} $, where the function $ \theta $ is $ C^2 $ on $ \R$. 
Taking the scalar product of (\ref{4.1}) with $ i \psi $ we get 
$$
\frac c2 \varrho ' + ( \varrho \theta ')'= 0 . 
$$
We infer that there is a constant $ k_1 \in \R$ such that $ \frac c2 \varrho  +  \varrho \theta ' =  k_1 $ in $ \R$.  
We have $ \varrho ( x ) \lra 1 $ as $ x \lra \pm \infty $ because $ \varrho - 1 \in H^1 ( \R)$. 
Since $ |\psi '|^2 = \frac{|\varrho '|^2}{4 \varrho } + \varrho |\theta '|^2 \in L^1( \R)$, 
we deduce that $ \theta ' \in L^2( \R)$, therefore we must have $ k_1 = \frac c2 $ and consequently
\beq
\label{4.4}
\theta ' = \frac c2 \frac{1 - \varrho}{\varrho}.
\eeq
Let $ c \neq 0 $.
If we are able to solve (\ref{4.2}) and if we get a solution $ \varrho $ such that $\varrho - 1 \in H^1( \R)$, 
it follows from (\ref{4.3}) that $\ds \inf_{x \in \R} \varrho ( x )  > 0 $, 
 and  then from (\ref{4.4}) we obtain $ \theta $ up to a constant. 
Then by (\ref{4.4}) we have $ \theta ' \in H^1 ( \R)$, $ \sqrt{\varrho } e^{ i \theta } \in \Er $ and it is straightforward to  see that $ \sqrt{\varrho } e^{ i \theta }$ is a solution of (\ref{4.1}). 
Moreover, all solutions of (\ref{4.1}) in $ \Er$ are obtained in this way. 
Notice that (\ref{4.1}) is invariant by translations and by multiplication by complex numbers of modulus  $1$, 
so the phase $\theta $ can be determined only up to a constant. 

If assumption {\bf (A1)} is satisfied, 
we have $ V( s ) = \frac 12 ( s - 1 )^2 + o ( (s-1)^2) $ and $ g( s, c ) = ( 2 - c^2 ) ( s - 1 )^2 +  o ( (s-1)^2) $  as $ s \lra 1$. 
If  $ c^2 < \vs^2 = 2$ we have $ g( s, c ) > 0 $ whenever  $ s $ is sufficiently close to  $1$ and $ s \neq 1$. 
Since $ g( 0, c ) = - c^2 \les 0$, we infer that there exists $ \zeta \in [0, 1 ) $ such that $ g( \zeta, c ) = 0 $ and we denote 
$ \zeta ( c ) = \ds \sup \{ \zeta \in [0, 1 ) \; | \; g( \zeta , c ) = 0 \}$. 
It is clear that $ g( \zeta (c), c ) = 0 $, $ g > 0 $ on $ ( \zeta(c), 1 ) $ and the mapping $ c \longmapsto \zeta (c ) $ is even on $ ( - \sqrt{2}, \sqrt{2}) $  and is strictly increasing on $(0, \sqrt{2})$. 
We denote 
$ D = \{ ( s, c ) \in  (0, 1 ) \times ( - \sqrt{2}, \sqrt{2}) \; | \; \zeta ( c ) < s < 1 \}. $
The set $D$ is connected, but not necessarily open. 
We consider a continuous function $ G : D \lra \R$ 
such that $\frac{ \p G}{\p s }$ exists and $\frac{ \p G}{\p s } ( s, c ) = \frac{1}{\sqrt{g(s, c ) }} $ for any $ (s, c ) \in D$, and
$G$ is $ C^2 $ in $\mathring{D}$.  
(Such a function exists: it suffices to take a smooth curve $ c \longmapsto a(c) $ defined on $ ( - \sqrt{2}, \sqrt{2}) $ 
such that $ \zeta ( c ) < a(c ) < 1 $ for all $ c $, then put $ G( s, c ) = \ii_{a(c) }^s \frac{1}{\sqrt{ g( \tau, c )}} \, d \tau .$)
For any fixed $c \in ( - \sqrt{2}, \sqrt{2})$, the mapping $s \longmapsto G(s, c)$ is strictly increasing on $ (\zeta(c), 1)$ and tends to $ \infty $ as $ s \lra 1 $ because $ \frac{1 - s}{\sqrt{g( s, c )}} \lra \frac{1}{\sqrt{2 - c^2}}$. 
Let $ L(c ) =\ds  \lim_{s \searrow \zeta(c) } G( s, c ) $. Then $ G( \cdot, c) $ is a $C^2-$diffeomorphism between $(\zeta(c), 1 )$ and $(L(c), \infty)$.

\begin{Proposition}
\label{P4.1}
Assume that condition {\bf (A1)} in the introduction  is satisfied and  $ c \in ( - \sqrt{2}, \sqrt{2} )$. 
Then: 

\medskip

i) Equation (\ref{4.1}) admits a solution $ \psi \in \Er $ satisfying $ \ds \inf_{x \in \R} |\psi ( x ) | < 1 $ if and only if 
$L(c) := \ds \lim_{s \searrow \zeta(c)} G( s, c ) $ is finite. 

Whenever  $ L(c) $ is finite, let 
$$ 
\varrho _c ( x ) = \left\{ \begin{array}{ll}  G(\cdot, c )^{-1} ( L(c ) - x ) & \mbox{if } x < 0 \\
\zeta(c ) & \mbox{if } x = 0 , \\
 G(\cdot, c )^{-1} ( L(c ) + x ) & \mbox{if } x > 0 
 \end{array}
 \right. 
\quad \mbox{ and } \quad 
\theta_c  ( x ) = \frac c2 \int_0 ^x \frac{ 1 - \varrho_c (s)}{\varrho _c ( s ) } \, ds \; \mbox{ if } c \neq 0. 
$$
If $ c \neq 0$ we define $ \psi_c ( x ) = \sqrt{ \varrho _c ( x )} e^{ i \theta_c (x)}$. \\ If $ c = 0 $ there are three subcases: 

- either  there exists $ s \in (0, 1)$ such that  $V(s) = 0$, then we have $ \zeta ( 0 ) > 0$  and we put  $\psi _0 (x ) = \sqrt{ {\varrho}_0 ( x )} $; 

- or  $ V > 0 $ on $ [0, 1)$, so that $ \zeta ( 0 ) = 0 $ and we put $ \psi _0 (x ) = \mbox{\rm sgn}(x) \sqrt{ {\varrho}_0 ( x )} $;  

 - or  $ V > 0 $ on $ (0, 1)$  and $ V(0) = 0$, and then $ V(s )  \les C s $ for $ s $ sufficiently small and consequently $ L( 0 ) = - \infty$; in 
 this subcase we do not define $ \psi _0$. 

Whenever $ \psi _c $ is defined as above, we have $ \psi _c \in \Er $ and $ \psi _ c $ is a solution of (\ref{4.1}).

\medskip

ii)  Equation (\ref{4.1}) admits a solution $ \psi \in \Er $ satisfying $ \ds \sup_{x \in \R} |\psi ( x ) | > 1 $ if and
only if 

$\bullet $ the mapping $ g( \cdot , c ) $ admits zeroes in $(1, \infty)$, and 

$\bullet $ denoting $ \tilde{\zeta} (c ) = \inf \{ \zeta > 1 \; | \; g( \zeta, c ) = 0 \}$ and by $ \tilde{G} (\cdot, c ) $ a primitive of 
$\frac{ 1}{\sqrt{ g( \cdot, c )}} $ on the interval $(1, \tilde{\zeta}_c ) $, the limit 
$\tilde{L}(c) := \ds \lim_{s \nearrow \tilde{\zeta}(c)} \tilde{G} ( s, c ) $ is finite. 

 In this case, define 
$$ 
\tilde{\varrho } _c ( x ) = \left\{ \begin{array}{ll}  \tilde{G}(\cdot, c )^{-1} ( \tilde{L}(c )  + x ) & \mbox{if } x < 0 ,\\
\tilde{\zeta}(c ) & \mbox{if } x = 0 , \\
\tilde{G}(\cdot, c )^{-1} ( \tilde{L}(c ) - x ) & \mbox{if } x > 0, 
 \end{array}
 \right. 
\quad \mbox{ and } \quad 
\tilde{\theta}_c  ( x ) = \frac c2 \int_0 ^x \frac{ 1 - \tilde{\varrho}_c (s)}{\tilde{\varrho} _c ( s ) } \, ds \; \mbox{ if } c \neq 0 .
$$
Let $ \tilde{\psi}_c ( x ) = \sqrt{ \tilde{\varrho}_c ( x )} e^{ i \tilde{\theta}_c (x)}$ if $ c \neq 0 $, 
respectively   $ \tilde{\psi}_0 ( x ) = \sqrt{ \tilde{\varrho}_0 ( x )} $ if $ c = 0 $.

Then $ \tilde{\psi} _c \in \Er $ and $ \tilde{\psi} _ c $ is a solution of (\ref{4.1}). 

\medskip

iii) Let $ \psi \in \Er $ be a nonconstant solution of (\ref{4.1}). 
We have either $ |\psi (x) | < 1 $ for any $ x \in \R$, or $ |\psi | > 1 $  for any $ x \in \R$.
For any $ \e \in ( 0, \sqrt{ 2 - c^2} ) $ there exist $ A_{\e}, B_{\e } $ (depending on $ \e $ and on $ \psi $) such that 
\beq
\label{4.5}
e^{- |x | \sqrt{ 2 - c^2 + \e ^2 } } \les \big| |\psi |^2 - 1 \big| \les   e^{- |x | \sqrt{ 2 - c^2 - \e ^2 } } 
 \qquad \mbox{ on } ( - \infty, A_{\e} ] \cup [B_{\e} , \infty ). 
\eeq

\medskip

iv) Any solution $ \psi \in \Er $ of (\ref{4.1}) satisfies $ |\psi ' (x)  |^2 = V( |\psi (x) |^2 ) $ for any $ x \in \R$. 
In particular, we have $ E ^1( \psi ) = 2 \ii_{\R} |\psi '|^2 \, dx = 2 \ii_{\R} V( |\psi |^2 ) \, dx $. 

\medskip

v) If the nonlinearity $F$ is locally Lipschitz, then any solution $ \psi \in \Er$ of (\ref{4.1}) is either a constant of modulus one, or is of the form 
$ e^{ i \al } \psi _c ( \cdot - x_0 ) $ or $ e^{ i \al } \tilde{\psi} _c ( \cdot - x_0 ) $ for some $ \al , x_0 \in \R$, where 
$ \psi_c $ and $ \tilde{\psi}_c $ are as in (i) and (ii), respectively. 

\medskip

vi) If $ c_0 \neq 0 $ and the mapping $ c \longmapsto \zeta ( c ) $ is differentiable at $ c_0$, then $ L(c_0)$ is finite. \\
Consequently, equation (\ref{4.1}) admits solutions in $ \Er $ for almost every $ c \in ( - \sqrt{2},  \sqrt{2})$. 

\end{Proposition}

{\it Proof. } 
(i) "$\Longrightarrow$"
Let $ \psi \in \Er $ be a solution of (\ref{4.1}) and let $ \varrho = |\psi |^2$. 
Then $\varrho \lra 1 $ at $ \pm \infty$, $ \varrho $ is a $ C^2 $ function on $ \R$ and satisfies (\ref{4.2}) and (\ref{4.3}). 
Assume that $ \inf \varrho < 1 $. 
Take $ x_1 \in \R$ such that $ \varrho ( x_1 ) \in ( \zeta(c), 1)$. 
Let $I = (a, b) $ be the maximal interval containing $ x_1 $ such that $ \varrho ( x ) \in (\zeta(c), 1) $ for any $ x \in I$. 
It follows from (\ref{4.3}) that $ \varrho ' \neq 0 $ on $I$, hence $ \varrho '$ has constant sign and $ \varrho $ is strictly monotonic on $I$. 
If $ \varrho '> 0 $ on $I$, from (\ref{4.3}) we get $\ds \frac{ \varrho '(x)}{ g( \varrho (x), c )} = 1 $ on $I$, and integrating we obtain 
$ G( \varrho (x ) , c ) = x + k_2$ on $I$, where $ k_2 $ is a constant, hence $ \varrho ( x ) = G(\cdot, c )^{-1} ( x + k_2) $ on $I$. 
Then necessarily $ b = \infty$. Indeed, if $ b $ is finite we have  $\zeta ( c ) <\varrho ( x ) < G(\cdot, c )^{-1} ( b + k_2) < 1 $ on $I$, hence there is some $ \e > 0 $ such that $\zeta ( c ) <\varrho ( x ) < 1 $ on $[b, b+ \e )$, contradicting the maximality of $I$. 
If $ L(c ) = - \infty$ we must have $ a = - \infty $ (for otherwise, $ \zeta ( c ) < G( \cdot, c )^{-1} (a  + k_2) = \varrho ( a) < \varrho ( x ) < 1$ for all $ x \in I$, and we would have   $ \zeta ( c ) < \varrho ( x ) $ for $ x \in ( a - \e, a ]$ for some positive $ \e$, contradicting again  the maximality of $I$). 
But if $ a = - \infty $ we have $ \ds \lim_{x \ra - \infty } \varrho (x ) = \lim_{x \ra - \infty }  G(\cdot, c )^{-1} ( x + k_2) = \zeta( c)$, 
impossible because $ \varrho (x ) \lra 1 $ as $ x \lra  - \infty$. 
Thus necessarily $ L(c) $ is finite. This implies that $ a $ is finite and $ \varrho ( a  ) = \zeta ( c )$, and we find 
$ a + k_2 = L(c)$. In conclusion, if $I$ is a maximal interval such that $  \zeta(c) < \varrho < 1 $ and $ \varrho ' > 0 $ on $I$ then necessarily $ I$ is of the form $(a, \infty )$ for some $ a \in \R$ and we have $ \varrho ( x ) = G(\cdot, c )^{-1} ( x - a + L(c)) $ on $I$. 
Similarly, if $I$ is a maximal interval such that $  \zeta(c) < \varrho < 1 $ and $ \varrho ' < 0 $ on $I$
we show that $ L(c) $ must be finite (for otherwise $I = \R$ and $ \varrho $ would not tend to $1$  at $  \infty$)
and $ I$ is of the form $(- \infty, b )$ for some $ b \in \R$ and 
  $ \varrho ( x ) = G(\cdot, c )^{-1} ( - x + b + L(c)) $ on $I$. 

"$\Longleftarrow$"
One easily proves that $ \varrho _c $ satisfies (\ref{4.3}) and (\ref{4.2}), and is $ C^2 $ in $ \R$.
It is obvious that  $ \theta _c $ satisfies (\ref{4.4}) and then en easy computation shows that $\psi _c $ solves (\ref{4.1}). 
By (\ref{4.5}) proven below we have $ \varrho _c - 1 \in L^2 ( \R)$, and then  (\ref{4.3}), the boundedness of $ \varrho _c $ and {\bf (A1)} imply that $ \varrho '\in L^2 ( \R)$. 
By (\ref{4.4}) we get $ \theta '\in L^2 ( \R)$ and then we infer that $ \psi _c \in \Er$. 
Notice that $ \varrho$  may vanish only if $ c = 0 $,  and  in this case we have  $ \theta _ 0 '= 0 $ in $ \R^*$ by (\ref{4.4}).

\medskip

(ii) The proof of (ii) is similar to the proof of part (i), so we omit it. 

\medskip

(iii) Let $ \psi $ be a solution of (\ref{4.1}) and let $ \varrho = |\psi |^2$. 
If $\varrho = 1$, we may write $ \psi = e^{i \theta}$ and (\ref{4.4}) implies that $ \theta '= 0 $, hence $ \psi $ is constant. 
Assume that there exist $ x_1, x_2 \in \R$ such that $ \varrho ( x _1 ) < 1 $ and $ \varrho ( x_2 ) \ges 1 $. 
If $ x_ 1 < x_2$, there exists $ x_3 \in ( x_1, x_2)$ such that $ \zeta ( c ) < \varrho ( x _3 ) < 1 $ and $ \varrho '( x_3 ) > 0 $. 
The argument in the proof of part (i) "$\Longrightarrow$" shows that the maximal interval $I$ containing $ x_3$ such that 
$ \zeta(c) < \varrho < 1 $ and $ \varrho '> 0 $ on $I$ is of the form $(b, \infty)$, contradicting the fact that $ \varrho ( x_2 ) \ges 1 $. 
If $ x_1 > x_2$, there exists $ x_3 \in ( x_2, x_1)$ such that $ \zeta ( c ) < \varrho ( x _3 ) < 1 $ and $ \varrho '( x_3 ) < 0 $. 
As above, we have then $ \zeta ( c ) < \varrho ( x _3 ) < 1 $ and $ \varrho '( x_3 ) < 0 $ on $( - \infty, x_3]$, contradicting the fact that $ \varrho ( x_2 ) \ges 1 $. 

A similar argument leads to a contradiction  if we assume that there exist  $ x_1, x_2 \in \R$ such that $ \varrho ( x _1 ) > 1 $ and $ \varrho ( x_2 ) \les 1 $.

Fix $ \e \in ( 0 , \sqrt{ 2 - c^2}) $. 
There is $ \de _{\e } > 0 $ such that 
$$
( 2 - c^2 - \e^2 ) ( s- 1)^2 \les g( s, c) \les ( 2 - c^2 + \e^2 )(s-1)^2 \qquad \mbox{ for any } s \in (1 - \de_{\e}, 1 + {\de_{\e}}). 
$$
and integrating we see that there exist constants $ C_1, C_2 \in \R$ such that 
\beq
\label{4.6}
C_1 - \frac{ \ln  | 1 -s | }{\sqrt{2 - c^2 + \e ^2 }} \les G( s, c ) \les C_2 - \frac{ \ln | 1 -s  |}{\sqrt{2 - c^2 - \e ^2 }}
\qquad \mbox{ for } s \in [ 1 - \de_{\e}, 1) \cup   (1, 1 +   \de_{\e}). 
\eeq
If $ \varrho < 1$, there exist $ a, b, k_1, k_2 \in \R $ such that $ a < b $, 
$ \varrho ( x ) \in ( 1 - \de_{\e}, 1) $ for all $ x \in ( - \infty, a ) \cup ( b , \infty )$, 
 and $ G( \varrho (x) , c ) =  - x + k_1 $ on $ ( - \infty , a )$, respectively $ G ( \varrho (x), c ) = x + k_2 $ on $( b, \infty)$. 
Then using (\ref{4.6}) we see that there are constants $ C_3, C_4 \in \R$ such that 
$$
C_3 e^{- |x | \sqrt{ 2 - c^2 + \e ^2 } } \les  \big| 1 - \varrho ( x ) \big| \les C_4 e^{ - |x | \sqrt{ 2 - c^2 - \e ^2}}
\qquad 
\mbox{ on } ( - \infty, a ) \cup ( b , \infty ).
$$
We obtain (\ref{4.5})  for any $ 0 < \e ' < \e $ 
by choosing conveniently $ A_{\e'} < a $ and $ B_{\e'} > b$. 
The proof of (\ref{4.5}) is similar if $ \varrho > 1 $. 

\medskip 

(iv) Taking the scalar product of  (\ref{4.1})  with $2 \psi '$ we get $ (|\psi ' |^2 ) ' - ( V(|\psi |^2 ) ) '= 0 $, hence 
$ |\psi ' |^2  -  V(|\psi |^2 )  $ is constant. Since $ |\psi '|^2 $ and $ V ( |\psi |^2) $ belong to $ L^1( \R)$, the constant must be zero.

\medskip 

(v)  Let $ \psi \in \Er $ be a traveling wave of speed $c$. Let $ \varrho = |\psi |^2$. Assume that $ \varrho (x) < 1 $ for some $x\in \R$.  
By (i) we know that $ L(c) $ is finite.
There is some $ x_1 \in \R $ such that $ \zeta ( c) < \varrho ( x_1 ) < 1 $ and $ \varrho ' ( x_1 ) < 0 $. 
As in the proof of part (i), there exists $ x_0  > x_1 $ such that 
$ \varrho ( x ) = G(\cdot, c)^{-1} ( L(c) + x_0 - x ) $ for all $ x \in ( - \infty, x_0)$ and $ \varrho ( x_0 ) = \zeta( c ) $. 
Using (\ref{4.3}) and the continuity of $ \varrho '$ we get 
$ \varrho '( x_0 ) = \lim_{x \uparrow x_0} \varrho '( x ) = \lim_{x \uparrow x_0} ( - \sqrt{g( \varrho (x), c) } ) = 0 $. 
Let $ \varrho _ c $ be as in (i). 
By (\ref{4.2}),  $(\varrho, \varrho ')$ and $ ( \varrho_c ( \cdot - x_0), \varrho_c '( \cdot - x_0) ) $ are both solutions of the Cauchy problem 
$$
\left \{
\begin{array}{l}
y'(x ) = z(x ) 
\\
z'(x ) = - c^2 ( y(x ) -1 ) + 2 V( y(x)) - 2 y(x ) F( y(x)) 
\\
(y(x_0), z( x_0) ) = ( \zeta(c), 0)
\end{array} 
\qquad \mbox{ in } [x_0, \infty).
\right. 
$$
Since $F$ is locally Lipschitz, the solution of the above Cauchy problem is unique and we infer that 
$ \varrho = \varrho _c ( \cdot - x_0) $ on  $[x_0, \infty)$, thus on $ \R$. 
Then using (\ref{4.4}) we see that the difference between the phase of $ \psi $  and $ \theta_c ( \cdot - x _0)$ is a constant, say $ \al $, and therefore $ \psi = e^{ i \al } \psi _c ( \cdot - x_0)$. 

The proof is analogous if there exists $ x \in \R$ such that $ \varrho ( x ) > 1 $. 

\medskip 

(vi) We have $ g ( \zeta(c), c ) = 0 $ for any $ c \in ( - \sqrt{2}, \sqrt{2})$. 
If $ \zeta $ is differentiable at $ c_0$, differentiating this equality with respect to $c$ we get 
$$
\p_1 g ( \zeta( c_0), c_0) \cdot \zeta ' ( c_0 ) + \p_2 g (  \zeta( c_0), c_0) = 0 ,
$$ 
that is $ \p_1 g ( \zeta( c_0), c_0) \cdot \zeta ' ( c_0 ) = 2 c_0 ( \zeta( c_0 ) - 1 )^2 > 0 $. 
We infer that $ \p_1 g ( \zeta( c_0), c_0) \neq 0 $. Since $ g ( \zeta( c_0), c_0) = 0 $ and 
$ g( s, c_0) > 0 $ on $ ( \zeta( c_0), 1)$, we have  $ \p_1 g ( \zeta( c_0), c_0) > 0 $.
Then 
$ \frac{\sqrt{ s - \zeta( c_0 )}}{\sqrt{g( s, c_0) }} \lra \frac{1}{\sqrt{ \p_1 g ( \zeta( c_0), c_0) }} > 0 $ as $ s \searrow \zeta (c_0)$, hence $ s \longmapsto \frac{1}{\sqrt{g( s, c_0)}}$ is integrable on an interval $ (\zeta(c_0), \zeta(c_0) + \e)$ and 
consequently $ L(c_0)$ is finite. 

It is well-known that a monotonic function is differentiable almost everywhere. 
\hfill
$\Box$

\medskip

The following question arises naturally: is it true that (\ref{4.1}) admits solutions in $ \Er$ for all but countably many $c$'s ? 

\begin{remark}
\label{R4.2} \rm 
One may compute the energy and the momentum of traveling waves provided by Proposition \ref{P4.1}. 

If $ c \in ( - \sqrt{2}, \sqrt{2}) $ and $ c \neq 0$, the function $ \psi _c $ admits a lifting $ \psi _c = \sqrt{ \varrho _ c } e^{ i \theta _c}$. Using Remark \ref{lift1}, 
 (\ref{4.4}), then the fact that $ \varrho _c $ is even and satisfies $ \varrho _c '( x ) = \sqrt{ g ( \varrho _c ( x ), c )}$ on $ ( 0, \infty)$ and performing the change of variable $ s = \varrho _c ( x)$ we find that a valuation of the momentum of $ \psi _ c $ is
$$
\int_{\R} ( 1 - |\psi_c |^2 ) \theta_c ' \, dx 
= \int_{\R} ( 1 - \varrho _c) \cdot \left( \frac c2 \frac{ 1 - \varrho_c}{\varrho _ c} \right) \, dx 
= c \int_0 ^{\infty } \frac{(1 - \varrho _ c )^2 }{\varrho_c } \, dx = c \int_{\zeta(c ) }^1 \frac{ ( 1 - s )^2}{s \sqrt{ g( s, c )}} \, ds. 
$$

If $ V $ vanishes on $(0,1)$ we have  $ \zeta( 0 ) > 0 $. If $ L(0)$ is finite, 
the function $ \psi_0 $ is real-valued and satisfies $ 0 < \zeta(0) \les \psi _0 < 1$. 
By Remark \ref{lift1}, a valuation of the momentum of $ \psi _0 $ is $0$. 

If $ V > 0 $ on $[0, 1)$,  we have $ \zeta (  0 ) = 0$. 
Take $ \de > 0 $ and $ \ph _0 \in C^{\infty }( \R)$ such that $ \ph _0 = - \pi $ on $(-\infty, - \de ]$ and $ \ph _0 = 0 $ on $ [\de, \infty)$. 
It is easily seen that $ w_0 := \psi _0 - e^{ i \ph _0 } \in H^1( \R)$. 
Using (\ref{analog}), (\ref{p}) and the fact that $ \langle i \psi _0 ', \psi_0 \rangle = 0 $ and 
$ \langle i w_0 , e^{ i \ph _0 } \rangle \lra 0 $ as $ x \lra \pm \infty$, an immediate computation shows that a valuation of the momentum of $ \psi _0 $ is $ p ( \ph_0, w_0 ) = \pi$. 

Using Proposition \ref{P4.1} (iv) and the change of variable $ s = \varrho _c ( x)$ we obtain 
$$
E^1( \psi_c ) = 2 \int_{\R} V( |\psi _c|^2) \, dx = 4 \int_0^{\infty } V( \varrho_c ) \, dx 
= 4 \int_{\zeta(c)}^1 \frac{ V(s)}{\sqrt{ g( s, c ) } } \, ds. 
$$
In particular, if $ V > 0 $ on $[0, 1)$ we have $ \zeta (0 ) = 0$ and 
$$ E^1( \psi _0 ) = 4 \int_0^1 \frac{V(s)}{\sqrt{4 s V(s)}} \, ds = 4 \int_0 ^1 \sqrt{V( \tau ^2 )} \, d \tau.$$

\end{remark}

\begin{example} 
\label{GP}
\rm
Consider the particular case of the Gross-Pitaevskii equation, namely $ F( s ) = 1- s. $
We have 
$$ V( s ) = \frac 12 ( 1 - s )^2 , \qquad  
 g ( s, c ) = (2s- c^2)( 1 - s)^2 , \qquad 
  \zeta ( c ) = \frac {c^2}{2} \quad \mbox{  and } \quad  g( s, c) > 0 \mbox{  if } s > 1.
$$
  One can use the change of variable $ t = \sqrt{ 2s - c^2} $ to compute a primitive of $\frac{1}{\sqrt{g ( \cdot, c )}} $ 
and it is easily seen that we may take 
\beq
\label{4.7}
G(s, c ) = \frac{1}{\sqrt{2 - c^2}}\ln \Big| \frac{\sqrt{ 2 - c^2} + \sqrt{ 2s - c^2}}{\sqrt{ 2 - c^2} - \sqrt{ 2s - c^2}} \Big| .
\eeq
With this choice of $G$ we have 
$ L( c ) = \ds \lim_{s \searrow \zeta ( c ) } G( s, c ) = \ds \lim_{s \searrow \frac{c^2}{2} } G( s, c ) = 0$ for all $ c \in ( - \sqrt{2}, \sqrt{2} ). $
For any fixed $ c$, the function $ G( \cdot, c ) $ is an increasing  diffeomorphism between 
$ ( \frac{c^2}{2}, 1 ) $ and $ (0, \infty)$. 
We find $ G( \cdot, c)^{-1} ( x ) = \frac{ c^2}{2} + \frac{ 2 - c ^2}{2} \left[ \tanh \left( \frac{ \sqrt{2 - c^2}}{2} x \right) \right]^2  . $
Proceeding as in Proposition \ref{P4.1} (i) we get
\beq
\label{4.8}
\varrho _ c ( x ) =  \frac{ c^2}{2} + \frac{ 2 - c ^2}{2} \left[ \tanh \left( \frac{ \sqrt{2 - c^2}}{2} x \right) \right]^2 \qquad \mbox{ for any } x \in \R. 
\eeq
Using (\ref{4.4}) and the change of variable $ t = \tanh \left( \frac{ \sqrt{ 2 - c^2}}{2} s \right) $ we compute $ \theta _c $ and we find
\beq
\label{4.9}
\theta _c ( x ) = \frac c2 \int_0 ^x \frac{1 - \varrho _c ( s )}{\varrho _c ( s ) } \, ds
= \arctan \left[ \frac{\sqrt{ 2 - c^2}}{c} \tanh \left( \frac{ \sqrt{2 - c^2}}{2} x \right) \right]. 
\eeq
Since $ \cos ( \arctan (z)) = \frac{1}{\sqrt{1 + z^2}} $ and $ \sin ( \arctan (z)) = \frac{z}{\sqrt{1 + z^2}} $, we finally obtain 
\beq
\label{4.10}
\begin{array}{rcl}
\psi_c ( x  ) = \sqrt{\varrho (x ) } e^{ i \theta (x)}  & = & \frac{c}{\sqrt{2}}
+ i \sqrt{1 - \frac{c^2}{2}} \tanh \left( \frac{ \sqrt{2 - c^2}}{2} x \right)
\\
\\
&  = &  i \left[ \sqrt{1 - \frac{c^2}{2}} \tanh \left( \frac{ \sqrt{2 - c^2}}{2} x \right) - i  \frac{c}{\sqrt{2}} \right].
\end{array}
\eeq
It follows from Proposition \ref{P4.1} (i), (ii), (v) that all traveling waves for the Gross-Pitaevskii equation in $ \Er$ are either constants of modulus one, or are of the form $ e^{ i \al } \psi _c ( \cdot - x_0)$ for some $ \al , x_0 \in \R$. 

Using Proposition \ref{P4.1} (iv), then (\ref{4.8}) and the change of variable $ t = \tanh \left( \frac{ \sqrt{ 2 - c^2}}{2} x \right) $ we get 
\beq
\label{4.11}
\begin{array}{rcl}
E^1( \psi _ c )  = E_{GL}^1( \psi _c) & = & \ds 2 \ii_{\R} |\psi_c '|^2 \, dx = 2 \ii_{\R}  V( |\psi _c  |^2 ) \, dx 
= \ii_{\R} ( 1 - \varrho _c (x) )^2 \, dx 
\\
\\
& = & \ds \left( \frac{ 2 - c^2}{2} \right)^2 \frac{2}{\sqrt{ 2 - c^2}} \ii_{-1}^1 1 - t^2 \, dt = \frac 23 ( 2 - c^2 )^{\frac 32}. 
\end{array}
\eeq

If $ c \neq 0$, it follows from Remark \ref{lift1}, identity  (\ref{4.4})  and the change of variable 
$ t = \tanh \left( \frac{ \sqrt{ 2 - c^2}}{2} x \right) $
that a valuation of the momentum of $ \psi _c $ is 
\beq
\label{4.12}
\begin{array}{l}
\ds \ii_{\R} ( 1 - \varrho _c ) \theta_c '\, dx  = \frac c2 \ii_{\R} \frac{ (1 - \varrho _c )^2 }{\varrho _c } \, dx
=
\frac c4 ( 2 - c^2 )^{\frac 32} \ii_{-1}^1 \frac{ (1 - t ^2)^2 }{\frac{ c^2}{2 } + \frac{ 2 - c^2}{2} t^2} \frac{1}{1 - t^2} \, dt 
\\
\\
\ds = \frac c4 ( 2 - c^2 )^{\frac 32}  \ii_{-1}^1 \frac{ 1 - t ^2 }{\frac{ c^2}{2 } + \frac{ 2 - c^2}{2} t^2} \, dt 
= \frac c2 ( 2 - c^2 )^{\frac 12}  \ii_{-1}^1 \frac{ 1 - t ^2 }{\frac{ c^2}{2 - c^2} + t^2} \, dt 
\\
\\
\ds = 2 \arctan \left( \frac{ \sqrt{ 2 - c ^2}}{c} \right) - c \sqrt{ 2 - c ^2}. 
\end{array}
\eeq
For more information on traveling waves for the Gross-Pitaevskii equation and for further references we refer to \cite{BGS-survey}, Section 2. 
It has been shown in \cite{BGS-survey} that the functions $ \psi _c $ minimize the energy when the momentum is kept fixed. 
Formulae (\ref{4.11}) and (\ref{4.12}) here above correspond to formulae (2.23) and (2.24) p. 63 in \cite{BGS-survey}, respectively. 

\end{example}

\begin{example}
\label{E4.3}
\rm 
Consider a function $V \in C^{\infty}( \R, \R)$ having the following properties: 

$\bullet $ There exists $ \de _1> 0 $ such that $ V(s ) = \frac 12 ( 1 - s )^2 $ for all $ s \in [ 1 - \de_1, 1 + \de_1] $, 

$ \bullet$ There exist $ c_0 \in (0, \sqrt{2})$, $ s_0 \in (0, 1 ) $ and $a > 0 $, $ \de _2 > 0 $ such that 
$$
V(s ) = \frac{1}{4s} \left( c_0 ^2 ( 1 - s )^2 + a^2 ( s - s_0 )^3 \right) 
\; \mbox{ on } ( s_0 - \de _2, s_0 + \de _2 ],  
\quad \mbox{ and } V( s ) > \frac{c_0 ^2 ( 1 - s)^2}{4s} \; \mbox{ on } (s_0, 1), 
$$

$ \bullet$ $ V( s ) \ges \frac{( s - 1)^2}{2 s} $ on $ (1, \infty)$. 

Let $ F(s ) = - V'(s)$. It is obvious that $F$ satisfies the assumption {\bf (A1)}.
We consider equation (\ref{1.1}) with nonlinearity $F$. 
Using  Proposition \ref{P4.1} (i), (ii) and (v) it is easily seen that 
for $c $  sufficiently close to $ \sqrt{2} $ (more precisely, for $ c \in [\sqrt{ 2 ( 1 - \de _1)}, \sqrt{2})$), 
traveling waves of speed $c$ for this equation in $ \Er $ are the same as traveling waves for the Gross-Pitaevskii equation 
(and they are either constant, or are equal to functions $ \psi _c $ in Example \ref{GP} up to a translation and a phase shift). 
Letting $ g( s, c ) = 4s V(s) - c^2 ( s - 1 )^2$, 
we have $ g ( s, c_0 ) = a^2 ( s - c _0 )^3 $ on $( s_0, s_0 + \de _2 ) $ and $ g( s, c_0 ) > 0 $ on $( s_0, 1 )\cup( 1, \infty)$,
hence $ \zeta( c_0 ) = s_0 $ and $ L( c_0 ) = - \infty$.
Then Proposition \ref{P4.1} (i) and (ii) implies that all traveling waves of (\ref{1.1}) in $\Er$ with speed $ c_0$ must be constant. 


 \end{example}

\begin{Lemma}
\label{L4.4}

Assume that $ V > 0 $ on $ [0, 1) $.
Then we have 
\beq
\label{4.13moins}
\inf \left\{E^1( \psi ) \; \mid \; \psi \in \Er \mbox{ and } \inf_{x \in \R } |\psi ( x ) |  = 0 \right\} = 4 \ii_0^1 \sqrt{V( s^2) } \, ds. 
\eeq
If the infimum is achieved by a function $ \psi$, then there exist $ x_0, \al _-, \al _+ \in \R$ such that 
\beq
\label{4.13}
 \psi (x ) = \left\{ \begin{array}{ll} \sqrt{\varrho_0 ( x - x_0 )} e^{ i \al _-} & \mbox{ if } x < x_0, 
\\
\\ 
\sqrt{\varrho _0 ( x + x_0 )} e^{ i \al _+} & \mbox{ if } x \ges x_0, 
\end{array}
\right.
\eeq
 where $ \varrho _0$  is as in Proposition \ref{P4.1} (i). 

\end{Lemma}

Equality  (\ref{4.13moins}) still holds if $ V> 0 $ on $ ( 0, 1) $ and $ V( 0 ) = 0 $, but the infimum is not achieved. 
In this  case we have $ L(  0 ) = - \infty $ and it follows from the proof of Proposition \ref{P4.1} (i) that there is no $ \rho _ 0 \in \Er $ satisfying (\ref{4.14}) below and $ \rho (0) = 0$. 

\medskip

{\it Proof. } For any $ \psi \in \Er$ we have  $ |\psi | - 1 \in  H^1 ( \R)$ (see Lemma \ref{L2.1} (ii)), hence $|\psi |$ is continuous and tends to $1$ at 
$ \pm \infty $.
If $ \ds \inf_{x \in \R } |\psi ( x ) | <1$, then the infimum is achieved at some $ x _ 0 \in \R$. 

Consider any $ \psi \in \Er $ such that $ \psi ( x _0 ) = 0 $ for some $ x _0 \in \R$. Take two sequences $ ( x_n ^{\pm} )_{n \ges 1 } \subset \R$ such that $ x_n^{\pm} \lra \pm \infty $. 
For $ n $ sufficiently large we have $ x_n ^- < x_0 < x_n ^+ $ 
and  using (\ref{escape}) we get 
$$
\ii_{x_n ^-} ^{x_0 } |\psi '|^2 + V( |\psi |^2 ) \, dx \ges 2 \big| H ( |\psi ( x_0 ) |) - H ( |\psi ( x_n ^-)|) \big|
$$
 and 
$$
\ii_{x_0} ^{x_n ^+}  |\psi '|^2 + V( |\psi |^2 ) \, dx \ges 2 \big|   H ( |\psi ( x_n ^-)|)  - H ( |\psi ( x_0 ) |)\big|.
$$
Summing up and letting $ n \lra \infty $ we get $ \ds E^1( \psi ) \ges 4  \big| H ( 0) \big| = 4 \ii_0 ^1 \sqrt{ V( s^2 ) } \, ds. $

If $ V( 0 ) > 0 $ it follows that any primitive of the function $ s \longmapsto \frac{1}{\sqrt{g( s, 0 )} } = \frac{1}{2\sqrt{ s V(s)} } $ has finite limit at $ 0 +$ (in other words, $L(0) $ is finite), and the function $ \varrho  _0 $ in Proposition \ref{P4.1} (i) is well-defined and satisfies
$$
\varrho _0 ' = - 2\sqrt{ \varrho _0 V( \varrho _0 ) } \mbox{ on } ( - \infty, 0 ), 
\quad \mbox{ respectively } \quad 
\varrho _0 ' = 2 \sqrt{ \varrho _0 V( \varrho _0 ) } \mbox{ on } ( 0,  \infty ). 
$$
Let $ \rho_0  = \sqrt{ \varrho _0}$. It is easily seen that 
\beq
\label{4.14}
\rho _0 '= - \sqrt{ V ( \rho _0^2) }  \; \mbox{ on } ( - \infty, 0 ), 
\qquad \mbox{ respectively } \qquad 
\rho _0 ' =  \sqrt{  V( \rho _0 ^2) } \; \mbox{ on } ( 0,  \infty ). 
\eeq
Using (\ref{4.14}) and the change of variable $ \tau = \rho _0 (x) $ we obtain
$$
\ii_{- \infty }^0 |\rho _ 0 '|^2 + V ( \rho _0 ^2) \, dx = - 2 \ii _ {- \infty }^0  \sqrt{V ( \rho _0 ^2 (x) )}  \rho _ 0 ' (x) \, dx  = 
2 \ii_ 0 ^1 \sqrt{V ( \tau ^2 )} \, d \tau
$$ 
and a similar computation holds on $[0, \infty)$, hence $ \ds E^1 ( \rho _ 0 ) = 4  \ii_ 0 ^1 \sqrt{V ( \tau ^2 ) } \, d \tau. $

Consider any $ \psi \in \Er $ such that $ E^1 ( \psi ) = 4  \ii_ 0 ^1 \sqrt{V ( \tau ^2 ) } \, d \tau $ 
and there exists $ x _0 \in \R$ such that $ \psi ( x_0 ) = 0 $. 
As above, using (\ref{escape}) we see that 
\beq
\label{4.15}
\ii_{- \infty }^{x_0 } \!   |\psi '|^2 + V( |\psi |^2 ) \, dx \ges 2 \!  \ii_ 0 ^1 \! \! \sqrt{V ( \tau ^2 ) } \, d \tau
\mbox{ and } 
\ii_{x_0 }^{- \infty }\!  |\psi '|^2 + V( |\psi |^2 ) \, dx \ges 2 \!  \ii_ 0 ^1 \! \! \sqrt{V ( \tau ^2 ) } \, d \tau. 
\eeq
Hence we must have equality in both inequalities in (\ref{4.15}) and we infer that  the point $ x_0 $ must be unique, 
and $ \rho := |\psi | $ must satisfy $ \rho '= \pm \sqrt{V( \rho ^2)} $ on each of the intervals 
$ ( - \infty, x_0)$ and $ ( x_0, \infty)$. Then it follows easily that $ \rho ( \cdot + x_0 ) $ satisfies (\ref{4.14}), and finally that  
$ \rho ( \cdot + x_0 ) = \rho _0 $. The function $ \frac{ \psi }{\rho }$ must be constant on each of the intervals $ ( - \infty, x_0 )$ and $ (x_0, \infty)$ (for otherwise, we would have $ \psi = \rho e^{ i \theta }$ on those intervals for some $ \theta \in H_{loc}^1$, and then 
$ |\psi ' | ^2 = |\rho '|^2 + \rho ^2 |\theta '|^2  \ges |\rho '|^2 $. 
If $ \theta '\not\equiv 0 $  we would get $\ii_{\R} |\psi '|^2 \, dx > \ii_{\R} |\rho '|^2 \, dx$, hence 
$ E^1 ( \psi ) > E^1( \rho )$, a contradiction). 
\hfill
$\Box$

\begin{Corollary}
\label{C4.5}
For any $ \psi \in \Er $ satisfying $ \ds E^1 ( \psi ) <  4  \ii_ 0 ^1 \sqrt{V ( \tau ^2 ) } \, d \tau $ we have 
$ \ds \inf_{ x \in \R} | \psi ( x ) | > 0 $ and  there exists a lifting $ \psi = \rho e^{ i \theta} $, where $ 1 - \rho \in H^1 ( \R)$ and 
$ \theta \in \dot{H}^1 ( \R)$. 

The same conclusion holds if  $ \ds E^1 ( \psi )  =   4  \ii_ 0 ^1 \sqrt{V ( \tau ^2 ) } \, d \tau $
and $ \psi $ is not one of the functions in (\ref{4.13}). 

\end{Corollary}

\begin{remark} \rm 
It can be shown that for any $ a \in ( 0 , 1) $ we have 
$$
\inf \left\{E^1( \psi ) \; | \; \psi \in \Er \mbox{ and } \inf_{x \in \R } |\psi ( x ) |  = a \right\} = 4 \ii_a^1 \sqrt{V( s^2)  } \, ds.
$$
Moreover, if $ V > 0 $ on $[0, 1)$, the only minimizers, up to translations in $ \R$ and multiplication by complex numbers of modulus $1$, are the functions
$ 
\rho _ a ( x ) = \left\{ \begin{array}{ll}
\rho_0 ( \cdot - b) & \mbox{ if } x < 0, \\
\rho_0 ( \cdot + b) & \mbox{ if } x \ges 0, 
\end{array}
\right.
$
where $ \rho _0 $ is as in (\ref{4.14}) and $b > 0 $ is chosen so that $ \rho_0 ( b ) = a $. 
\end{remark}

\section{Minimizing the energy at fixed momentum in $ \Er$} 

We define 
\beq
\label{E1-min}
E_{\min}^1 ( p ) = \inf \{ E^1( \psi ) \; \big| \; \psi \in \Er \mbox{ and } \PR (\psi ) = \pr{p} \}. 
\eeq

We collect in the next Lemma the main properties of the function $E_{\min}^1$.

\begin{Lemma}
\label{E1min}
Assume that $V$ satisfies {\bf (A1)}. 
The function $ E_{\min}^1 $ has the following properties: 

\medskip

i) $ E_{\min}^1 $ is non-negative, $2 \pi-$periodic, $ E_{\min}^1 ( -p ) = E_{\min}^1 ( p )$ for all $ p \in \R$ and 
$$
E_{\min}^1 ( p ) = \inf \{ E^1( e^{ i \ph } + w ) \; \big| \; \ph \in \dot{H}^1( \R, \R), 
w \in H^1 ( \R, \C) 
 \mbox{ and } p (\ph , w ) = p \}.
$$
In other words, we have $ E_{\min}^1 ( p ) = E_{\min}^1 ( |\pr{p} | )$ for any $ p \in \R$. 

\medskip

ii) $ E_{\min}^1( p ) \les \sqrt{2} p $. 

\medskip

iii) For any $ \e > 0 $ there exists $ p_{\e } > 0 $ such that $ E_{\min}^1( p ) \ges (1 - \e) \sqrt{2} p $
 for any $ p \in (0, p_{\e})$. 
 
 \medskip

 iv) $ E_{\min}^1$ is sub-additive: for any $ p_1, p_2 \in \R$ there holds 
 $$ E_{\min}^1( p_1 + p_2 ) \les E_{\min}^1 ( p_1 ) + E_{\min}^1 ( p_2).$$
 
v) $ E_{\min}^1$  is Lipschitz on $ \R$ and its best Lipschitz constant is $ \vs = \sqrt{2}$. 

\medskip

vi) Assume that there exists $ \de > 0 $ such that 
$ V( s ) \les \frac 12 ( 1 - s )^2 + \frac 38 ( 1 - s )^3 $ for any $ s \in [1 - \de, 1).$     
(This condition is fulfilled, for instance, if $F$ is $C^2$ near $ 1$ and $ F ''(1) < \frac 94$.) 
Then we have $E_{\min}^1 ( p ) < \sqrt{ 2 } p $ for any $ p > 0 $. 

\medskip

Assume, in addition, that $ V > 0 $ on $[0, 1)$. Then: 

\medskip

vii)  For any $ p \in \R$ we have $ \ds E_{\min}^1 ( p ) \les  4  \ii_ 0 ^1 \sqrt{V ( \tau ^2 ) } \, d \tau$. 

\medskip

viii) $ E_{\min}^1 $ is nondecreasing on $[0, \pi]$ and is concave on $[0, 2 \pi]$.

\end{Lemma}

{\it Proof. } (i) 
For any given  $  \ph \in \dot{H}^1( \R, \R)$ and $ w \in H^1( \R)$, let $ \tilde{\ph}( x ) = \ph ( - x ) $ and 
$ \tilde{w}(x ) = w ( - x)$. Then we have $  \tilde{\ph } \in \dot{H}^1( \R, \R)$, $\tilde{ w} \in H^1( \R)$, 
$ p ( \tilde{\ph}, \tilde{w}) = - p ( \ph, w )$ and $ E^1(e^{i \tilde{\ph}} + \tilde{w} )= E^1( e^{ i \ph } + w )$. 
This implies that $ E_{\min}^1 ( -p ) = E_{\min}^1 ( p )$ for any $ p \in \R$.

Let $ p \in \R$. Let $ k \in \Z$. 
Consider $ \ph \in \dot{H}^1 ( \R) $ and $ w \in H_{per}^1 $ satisfying $ p (\ph , w ) = p $.
Let $ \chi \in C^{\infty } ( \R)$ such that $ \chi = 0 $ on $ ( - \infty, 0 ]$ and $ \chi = 1 $ on $ [1, \infty)$. 
Let $\tilde{ \ph } ( x ) = \ph ( x ) + 2 k \pi \chi (x )$ and $ \tilde{w} (x, y) = e^{ i \ph (x )  }  - e^{ i \tilde{ \ph }(x ) } + w(x,y)$.
Then $ p ( \tilde{\ph} , \tilde{w} ) = p( \ph, w ) + 2 k \pi  = p + 2 k \pi$ and $  e^{ i \ph } + w  =  e^{ i \tilde{\ph }} + \tilde{ w}$. Since for any $( \ph, w ) \in \dot{H}^1(\R, \R) \times H_{per}^1$ satisfying $ p (\ph , w ) = p $ we may construct 
$( \tilde{\ph}, \tilde{w})$ as above, 
 we conclude that $   E_{\min}^1 ( p  + 2 k \pi) = E_{\min}^1 ( p )$. 
The  rest of part (i) is obvious. 

\medskip

(ii) The proof is very similar to the proof  of Lemma 3.3 p. 604 in \cite{BGS} and of Lemma 4.5 p. 173 in \cite{CM}. 
Take $ \chi \in C_c^{\infty}( \R)$ such  that $ \ii_{\R} | \chi'(x )|^2 \,  dx = 1 $ and $ \ii_{\R}  \chi'(x )^ 3\,  dx = 0 $
(for instance, we may take $\chi $ an even function). 
Let $ A = \ii_{\R} |\chi ''( t ) |^2 \, dt $ and $ B =   \ii_{\R} |\chi ' ( t ) |^4 \, dt $.
For $ \e , \la , \si > 0 $ (to be chosen later), let 
$ \rho_{\e, \la }( x ) = 1 - \frac{\e}{\la} \chi ' \left( \frac{x}{\la } \right) $, 
$ \theta _{\la, \si } ( x ) = \si \chi  \left( \frac{x}{\la } \right) $ and 
$ \psi _{\e, \la, \si  }( x ) = \rho_{\e , \la }(x ) e^{ i \theta _{\la, \si} (x)}. $
It is clear that $ \psi_{\e , \la, \si  } \in \Er $  and a simple computation gives 
$$
\begin{array}{l}
\ds \ii_{\R} |\psi_{\e , \la, \si  } '(x ) |^2 \, dx 
= \ii_{\R } |\rho_{\e, \la } ' ( x) |^2 + |\rho_{\e, \la }|^2 |\theta _{\la, \si  } ' ( x ) |^2 \, dx 
\\
\\
\ds = \ii_{\R} \frac{ \e ^2}{\la ^3} | \chi ''(t )|^2 + \frac{\si ^2}{\la } |\chi' ( t ) |^2 \left( 1 - \frac{ 2 \e}{ \la } \chi' ( t ) + \frac{ \e ^2}{\la ^2 } |\chi' ( t ) |^2 \right) \, dt
= \frac{ \si^2}{\la } + A \frac{ \e^2}{\la ^3} + B \frac{ \e^2 \si^2}{\la ^3}, 
\end{array}
$$
$$
\ii_{\R} ( 1 - \rho_{\e, \la }^2 ) \theta_{\la, \si }' ( x ) \, dx 
= \ii_{\R } \frac{ 2 \si \e }{\la } |\chi' ( t ) |^2  - \frac{ \si \e ^2}{\la ^2} \chi ' (t ) ^3 \, dt  = \frac{ 2 \si \e}{\la}, 
$$
$$
\ii_{\R} ( 1 - \rho_{\e, \la }^2 ) ^2 \, dx = \ii_{\R} \frac{4 \e ^2}{\la } |\chi'(t) |^2 - \frac{4 \e ^3}{\la ^2} \chi'(t) ^3 + \frac{\e ^4}{\la ^3} | \chi'(t ) | ^4 \, dt = \frac{4 \e^2}{\la } + B \frac{\e ^4}{\la ^3}  . 
$$

Fix $ p > 0 $. For $ \la > 0 $ we choose $ \e = \e( \la ) = 2^{ - \frac 34} \sqrt{ p \la}$ and $ \si = \si(\la ) = 2^{ - \frac 14} \sqrt{ p \la}$.
Let $ \psi_{\la } = \psi_{\e( \la ), \la, \si(\la )}$. 
As $ \la \lra \infty $ we have $ \frac{ \e( \la)}{\la } \lra 0 $, $ \frac{ \si( \la)}{\la } \lra 0 $, and
$ \big| 1 - |\psi_{\la } | \big| = \big| 1 - \rho_{\e(\la), \la } \big| \les  \frac{ \e( \la)}{\la }  \| \chi ' \|_{L^{\infty}} \lra 0 $. 
For all $ \la $ sufficiently large we have $ \psi_{\la } \in \Er $ and the above computations show that a valuation of the momentum of $ \psi _{\la } $ is $p$, and 
$$
\ds \ii_{\R} |\psi_{ \la    } '(x ) |^2 \, dx = 2^{- \frac 12 } p + 2^{ - \frac 32}  \frac{p}{\la ^2} A + \frac{p^2}{4 \la } B
\lra \frac{p}{\sqrt{2}} \quad \mbox{ as } \la \lra \infty. 
$$
By assumption {\bf (A1)} we have $ V( s ) = \left( \frac 12 + o ( |s -1 |) \right) ( s -1 )^2 $ as $ s \lra 1 $, hence 
$$
\ii_{\R} V( | \psi_{\la} |^2) \, dx 
= \left( \frac 12 + o \left( \frac{ \e( \la)}{\la }  \| \chi ' \|_{L^{\infty}} \right) \right) \ii_{\R} ( 1 - \rho_{\e (\la), \la }^2 ) ^2 \, dx 
\lra \frac{p}{\sqrt{2}} \quad \mbox{ as } \la \lra \infty. 
$$
Since $ E_{\min}^1 ( p ) \les E^1 ( \psi_{\la})$ for all sufficiently large $ \la$ and 
$ E^1 ( \psi _{\la } ) \lra  \sqrt{2} p $ as $ \la \lra \infty$, (ii) follows. 

\medskip

(iii)  Fix $ \e > 0 $. We may assume that $ \e \les \frac 12$. By assumption {\bf (A1)} there exists $ \de = \de( \e ) > 0 $ 
such that $ ( 1 - \de )^2 \ges 1 - \e $ and 
$$
V( \rho ^2 ) > \frac 12 ( 1 - \e) ( 1 - \rho ^2 )^2 \qquad \mbox{ for any } \rho \in [1 - \de, 1 + \de].
$$
By Lemma \ref{L2.1} (i) there exists $ \kappa > 0 $ such that for any $ \psi \in \Er $ satisfying $ E ^1 ( \psi ) \les \kappa $ we have $ \| 1 - |\psi | \|_{L^{\infty} ( \R ) } \les \de$. 
Let $ p_{\e } = \min \left( \frac{ \kappa}{2 \sqrt{2}} , \frac{\pi}{4} \right)$ and let $ p \in (0, p_{\e}]$.  
Consider any $ \ph \in \dot{H}^1 ( \R)$ and any $ w \in H^1 ( \R)$ such that  $ p ( \ph, w ) = p $ and $ E ^1 ( e^{ i \ph } + w ) \les 2 \sqrt{2} p $. 
Denoting $ \psi =  e^{ i \ph } + w$, we have $ \psi \in \Er$ and $ E^1 ( \psi ) \les \kappa$, and then  Lemma \ref{L2.1} (i) implies that 
$ \| 1 - |\psi | \|_{L^{\infty} ( \R ) } \les \de$. We infer that $ \psi $ admits a lifting $ \psi = \rho e^{ i \theta }$ where $ \theta \in \dot{H}^1 ( \R)$ and $ 1 - \de \les \rho ( x ) \les 1 + \de $ on $ \R$. Then we have 
\beq
\label{minoration}
\begin{array}{l}
\ds E^1 ( \psi ) = \int_{\R} |\rho ' |^2 ( x ) + \rho^2 ( x ) |\theta ' |^2 ( x ) + V( \rho ^2 ( x)) \, dx 
\\
\\
\ds \ges \int_{\R} ( 1 - \de )^2 |\theta ' |^2 ( x )  + \frac 12 ( 1 - \e) ( 1 - \rho ^2 ( x ))^2 \, dx 
\ges ( 1 - \e ) \sqrt{2} \Big| \int_{\R} ( 1 - \rho ^2 ( x) ) \theta '(x ) \, dx \Big|. 
\end{array}
\eeq
By Definition \ref{mom1} and Remark \ref{lift1}, 
$ p = p( \ph, w )$ and $ \ds \int_{\R} ( 1 - \rho ^2 ( x) ) \theta '(x ) \, dx  $ are both valuations of the momentum of $ \psi$, hence 
$ \ds \int_{\R} ( 1 - \rho ^2 ( x) ) \theta '(x ) \, dx  = p + 2 k \pi $ for some $ k \in \Z$. 
We have $ 0 < p \les p_{\e } \les \frac{\pi}{4} $ and by (\ref{minoration}) we get 
${\ds \Big| \int_{\R} ( 1 - \rho ^2 ( x) ) \theta '(x ) \, dx \Big| } \les \frac{ E^1( \psi ) }{( 1 - \e ) \sqrt{2} }\les 4 p_{\e} <  \pi$, and we conclude that necessarily $ \ds \int_{\R} ( 1 - \rho ^2 ( x) ) \theta '(x ) \, dx  = p$. 
Then (\ref{minoration}) gives $ E^1( \psi ) = E^1( e^{i \ph } + w ) \ges ( 1 - \e ) \sqrt{2} p $. 
Since this inequality holds for any $ \ph $ and $ w$ as above, (iii) follows. 

\medskip 

(iv)  Let $ p_1, p_2 \in \R$.  Fix $ \e > 0 $. 
By Lemma \ref{approx} (i) there exist $ \ph_1, \ph _2 \in C^{\infty} ( \R) $, 
$ w_1, w_2 \in C_c ^{\infty } ( \R)$,  $ A > 0 $ and $ \al _j, \, \beta _j \in \R$ such that 
$ \mbox{ supp}( w_j ) \subset ( - A, A)$, $ \ph _j = \al _j $ on $ ( - \infty, - A]$, $ \ph _j = \beta_j $ on $[A, \infty )$,  
$$
p ( \ph_j, w_j ) \in \pr{p_j}  \qquad \mbox{ and } \qquad E^1(e^{ i \ph _j  } + w _j ) < E_{\min}^1 ( p_j ) + \frac{ \e }{2} 
\qquad \mbox{ for } j = 1,2.  
$$
Let 
$$
 \psi ( x ) = \left\{ \begin{array}{l} e^{ i \ph _1 (x ) } + w_1 (x ) \mbox{ if } x \les 2A, 
\\
\\
e^{i ( \beta_1 - \al _2)} \left(  e^{ i \ph _2 (x - 3A ) } + w_2 (x - 3A) \right) \mbox{ if } x > 2A,
\end{array}
\right.
$$
$$
\ph ( x ) = \left\{ \begin{array}{l} 
\ph_1 ( x ) \mbox{ if } x \les 2A, \\
\beta _1 - \al _2 + \ph_2 ( x - 3A) \mbox{ if } x > 2A, 
\end{array} \right.
\quad 
w( x ) = \left\{ \begin{array}{l} 
w_1 ( x ) \mbox{ if } x \les 2A, \\
e^{i ( \beta _1 - \al _2 )} w_2 ( x - 3A) \mbox{ if } x > 2A. 
\end{array} \right.
$$
Then we have $ \psi \in \Er \cap C^{\infty}( \R)$. It is easy to see that 
$ p ( \ph, w ) = p ( \ph _1, w_1) + p ( \ph _2, w_2)$, and
$$ \PR ( \psi ) = \PR( e^{ i \ph _1} + w_1 ) + \PR \left( e^{i ( \beta_1 - \al _2)} \left( e^{i ( \ph_2 ( \cdot - 3A))} + w_2 ( \cdot - 3A)\right)
 \right) = \pr{p_1} + \pr{p_2}.  
 $$
 We infer that 
$$ 
\begin{array}{rcl}
E_{\min}^1 ( p_1 + p_2 )  & \les & 
E^1 ( \psi ) = E^1 \left(  e^{ i \ph _1} + w_1  \right) + 
E^1 \left( e^{i ( \beta_1 - \al _2)} \left( e^{i ( \ph_2 ( \cdot - 3A))} + w_2 ( \cdot - 3A)\right) \right)
\\
\\
 & = &  E^1 \left(  e^{ i \ph _1} + w_1  \right) +  E^1 \left(  e^{ i \ph _2} + w_2  \right)  < E_{\min}^1 ( p_1 ) + E_{\min}^1 ( p_2 ) + \e. 
 \end{array}
$$
Since $ \e $ is arbitrary, (iv) follows. 

\medskip 

(v) The sub-additivity of $E_{\min}^1$ and part (ii) imply that 
$$
 | E_{\min}^1 ( p_2 ) -  E_{\min}^1 ( p_1 ) | \les  E_{\min}^1 ( p_2 - p_1 ) \les \sqrt{2} | p_2 - p_1| 
 \qquad \mbox{for any } p_1, p_2 \in \R.
 $$
Part (iii) implies that  $ \sqrt{2}$ is the best Lipschitz constant of $E_{\min}^1$. 

\medskip 

(vi) The sub-additivity of $E_{\min}^1 $ gives $ E_{\min}^1 ( np ) \les n E_{\min}^1 ( p )$ for any $ p > 0 $ and any $ n \in \N^*$. 
Hence it suffices to show that $ E_{\min}^1 ( p ) < \sqrt{2} p $ for sufficiently small $ p $.

We use as "test functions" the traveling-waves for the Gross-Pitaevskii equation in Example \ref{GP}.  
Proceeding as in (\ref{4.11}) and using the change of variable $ t = \tanh \left( \frac{ \sqrt{ 2 - c ^2}}{2} x \right) $ we get 
\beq
\label{5.2bis}
\ii_{\R} \! \left( 1 - |\psi _c |^2 \right) ^3  dx = \! \ii_{\R} \! ( 1 - \varrho _c (x) )^3 \,   dx 
= \! \ds \left( \frac{ 2 - c^2}{2} \right)^3 \! \frac{2}{\sqrt{ 2 - c^2}} \ii_{-1}^1 (1 - t^2 )^2 \, dt = \frac{4}{15} ( 2 - c^2 )^{\frac 52}. 
\eeq
There is some $ c _{\de} \in (0, \sqrt{2}) $ such that for all $ c \in (c_{\de}, \sqrt{2}) $ we have $ 1 - \de < |\psi _c | < 1 $. 
Using the fact that $ V( s ) \les \frac 12 ( 1 - s )^2 + \frac 38 ( 1 - s )^3 $ for  $ s \in [1 - \de, 1)$, (\ref{4.11})  and (\ref{5.2bis})
we get 
$$
E^1( \psi _c ) \les E_{GL}^1 ( \psi _ c ) + \frac 38 \ii_{\R} \left( 1 - |\psi _c |^2 \right) ^3\, dx
\les 
\frac 23  ( 2 - c^2 )^{\frac 32} + \frac{1}{10}  ( 2 - c^2 )^{\frac 52}.
$$
It is useful to denote $ \e ( c ) = \sqrt{2 - c^2}$, so that $ \e ( c ) \lra 0 $ as $ c \nearrow \sqrt{2}$. 
The above inequality can be written as 
\beq
\label{E1fe}
E^1( \psi _c ) \les f ( \e ( c)), \qquad \mbox{ where } f ( \e ) = \frac 23 e^3 + \frac{1}{10} \e ^5. 
\eeq
Recall that by (\ref{4.12}), a valuation of the momentum of $ \psi _c $ is 
$$
m(c) :=  2 \arctan \left( \frac{ \sqrt{ 2 - c ^2}}{c} \right) - c \sqrt{ 2 - c ^2}
= 2 \arctan \left( \frac{ \e (c)}{\sqrt{ 2 - \e^2 (c) }} \right) - \e (c) \sqrt{2 - \e ^2(c)} = : g( \e ( c)). 
$$
We have $ g ( 0 ) = 0 $ and an elementary computation gives $ g'(\e ) = \frac{2 \e^2}{\sqrt{ 2 - \e^2}} $ for $ \e \in (0, \sqrt{2}). $
The function $ g $ is increasing and continuous on $[0, \sqrt{2})$, 
and therefore $ c \longmapsto m(c)$ is decreasing, positive on $(c_{\de} , \sqrt{2})$, and tends to $0 $ as $ c \nearrow \sqrt{2}$. 
We have 
$$
E_{\min}^1 ( m(c) ) \les E^1( \psi _c ) \les f ( \e (c))
$$
and it suffices to show that $ f( \e ) < \sqrt{2} g( \e ) $ for all sufficiently small $ \e$. 
A straightforward computation gives
$$
g'( \e ) = \sqrt{2} \left( \e^2 + \frac 14 \e^4 + \frac{3}{32} \e^6 + o( \e^6 ) \right) 
\quad \mbox{ and }
$$
$$
g( \e ) = \sqrt{2} \left( \frac 13 \e^3 + \frac{1}{20} \e^5 + \frac{3}{224} \e^7 + o( \e^7 ) \right) 
\quad \mbox{ as } \e \lra 0.
$$
By (\ref{E1fe}) there is $ \e_0 > 0 $ such that $ f( \e ) < \sqrt{2} g( \e ) $ for any $ \e \in (0, \e _0)$, as desired.

\medskip 

(vii) Let $ p \in \R$. Fix $ \e > 0 $. Choose $ \de > 0 $ such that $ \de V( 0 ) < \e$. 
Let $ \varrho _0 $ be as in Proposition \ref{P4.1} (i). Define 
$$
\rho (x ) = \left\{ \begin{array}{ll} \sqrt{\varrho_0 ( x )}  & \mbox{ if } x <0, 
\\
0 &  \mbox{ if } 0 \les x \les  \de , 
\\ 
\sqrt{\varrho _0 ( x - \de)}  & \mbox{ if } x \ges x_0. 
\end{array}
\right.
$$
It is easily seen that $ 1 - \rho \in H^1( \R)$. 
Choose $ \theta \in C^{\infty}( \R)$ such that $ \theta $ is constant on $(- \infty, 0 ]$ and on $[\de, \infty)$, 
and $ \ii _0 ^{\de} \theta ' \, dx = \theta ( \de )  - \theta ( 0 ) = p. $
Let $ w = ( \rho - 1 ) e^{ i \theta}$ and $ \psi = \rho e^{ i \theta } = e^{ i \theta } + w$. 
We have $ w \in H^1 ( \R)$, $ \psi \in \Er $ and using  Remark \ref{lift1} we get 
$$
p ( \theta, w ) = \ii_{\R} ( 1 - \rho ^2 ) \theta ' \, dx = p . 
$$
We have $  |\rho '|^2 + V( \rho ^2 ) = V( 0 ) $ on $ (0, \de) $ and consequently 
$$
E_{\min}^1 ( p ) \les E^1( \psi ) = \ii_{( - \infty, 0 ] \cup [\de, \infty) }\left(  |\rho '|^2 + V( \rho ^2 ) \right) \, dx + \de V( 0 ) 
\les 4  \ii_ 0 ^1 \sqrt{V ( \tau ^2 ) } \, d \tau +\e. 
$$
Since $ \e $ was arbitrary, the conclusion follows. 

\medskip

(viii) We proceed in several steps. 

{\it Step 1. "Reflection" of functions in $ \Er$ that have a lifting. } 
Assume that $ \psi \in \Er $ can be written in the form $ \psi = \rho e^{ i \theta }$, where $ \rho $ and $ \theta$ are real-valued functions, 
$ 1 - \rho \in H^1( \R) $ and $ \theta \in \dot{H}^1( \R)$. 
For any $ t \in \R$ we define $ \psi_t ( x ) = e^{ - i \theta ( t )} \psi ( x) = \rho ( x ) e^{ i (\theta (x ) - \theta (t))}$. 
It is obvious that 
$$ \psi _t \in \Er,  \quad \psi _t ( t ) = \rho ( t ) \in [0, \infty), \quad 
E^1( \psi _t ) = E^1( \psi ) \quad \mbox{ and } \quad \PR (\psi_t) = \PR (\psi) . 
$$
We define 
\beq
\label{lop1}
\psi_{t, 1} ( x ) = \left\{ \begin{array}{l} 
\psi _t ( x ) \mbox{ if } x \les t , 
\\
 \ov{\psi _t} (2 t -x) \mbox{ if } x > t, 
\end{array}
\right.
\qquad
\psi_{t, 2} ( x ) = \left\{ \begin{array}{l} 
 \ov{\psi _t}( 2 t - x ) \mbox{ if } x < t, 
\\
 \psi _t(x) \mbox{ if } x \ges t, 
\end{array}
\right.
\eeq
\beq
\label{lop2}
\theta_{t, 1} ( x ) = \left\{ \begin{array}{l} 
\theta ( x ) - \theta(t)  \mbox{ if } x \les t , 
\\
 - \theta(2 t -x)  + \theta (t) \mbox{ if } x > t, 
\end{array}
\right.
\qquad
\theta_{t, 2} ( x ) = \left\{ \begin{array}{l} 
 - \theta( 2 t - x )  + \theta ( t) \mbox{ if } x < t, 
\\
 \theta(x) - \theta (t) \mbox{ if } x \ges t, 
\end{array}
\right.
\eeq
\beq
\label{lop3}
\rho_{t, 1} ( x ) = \left\{ \begin{array}{l} 
\rho ( x ) \mbox{ if } x \les t , 
\\
 \rho (2 t -x) \mbox{ if } x > t, 
\end{array}
\right.
\qquad
\rho_{t, 2} ( x ) = \left\{ \begin{array}{l} 
 \rho ( 2 t - x ) \mbox{ if } x < t, 
\\
 \rho (x) \mbox{ if } x \ges t. 
\end{array}
\right.
\eeq
It is easy to check that $ \rho_{t, j } \in H^1( \R)$, $ \theta_{t, j } \in \dot{H}^1( \R)$, $ \psi _{t, j } \in \Er$ and
$ \psi _{t, j } = \rho_{t, j } e^{ i \theta_{t, j }}  =  e^{ i \theta_{t, j }}  + w_{ t, j }$ for $ j = 1, 2$, where 
$ w_{t, j } =   (\rho_{t, j } -1 )e^{ i \theta_{t, j }} \in  H^1( \R)$. 
An immediate computation gives 
\beq
\label{lop4}
E^1(\psi_{t, 1} ) = 2 \ii_{- \infty}^t |\psi '|^2 + V( |\psi |^2 ) \, dx
\quad \mbox{ and } \quad 
E^1(\psi_{t, 2} ) = 2 \ii_t^{ \infty} |\psi '|^2 + V( |\psi |^2 ) \, dx, 
\eeq
\beq
\label{lop5}
\begin{array}{l}
\ds p( \theta_{t, 1} , w_{t, 1}) = \ii_{\R} ( 1 - \rho_{t, 1}^2 ) \theta _{t, 1} ' \, dx = 2  \ii_{- \infty}^t ( 1 - \rho^2) \theta '\, dx, 
\\
\\
\ds p( \theta_{t, 2} , w_{t, 2}) = \ii_{\R} ( 1 - \rho_{t, 2}^2 ) \theta _{t, 2} ' \, dx = 2  \ii_t^{\infty} ( 1 - \rho^2) \theta '\, dx, 
\end{array}
\eeq
so that 
\beq
\label{lop6}
E^1(\psi_{t, 1} ) + E^1(\psi_{t, 2} ) = 2E^1( \psi ) \quad \mbox{ and } \quad 
p( \theta_{t, 1} , w_{t, 1}) +  p( \theta_{t, 2} , w_{t, 2}) = 2 p \left( \theta, ( \rho - 1) e^{ i\theta} \right).
\eeq

\medskip

{\it Step 2. For any $ p \in (0, \pi]$ satisfying $ \ds E_{\min}^1( p ) < 4 \ii_0^1 \sqrt{V( s^2) } \, ds $ and for any $ \e > 0 $ there exists 
$ \psi \in \Er $ such that $ E( \psi ) <  E_{\min}^1( p ) + \e$, $ \psi = \rho e^{ i \theta} $ with 
$ \rho \in H^1( \R)$, $ \theta \in \dot{H}^1( \R)$ and $ p \left( \theta, ( \rho - 1) e^{ i \theta} \right) = p$, and $ \psi $ is constant on  $( - \infty, -A]$ and on $[A, \infty)$ for some $ A >0$. }

Let $ p $ be as above. Let $ 0 < \e < \ds 4 \ii_0^1 \sqrt{V( s^2) } \, ds  - p$. 
By the definition of $E_{\min}^1 $ and by Lemma \ref{approx} (i), there exists $ \psi \in \Er \cap C^{\infty} ( \R)$ such that $\PR ( \psi ) = \pr{p} $, 
$ E^1( \psi ) < E_{\min}^1 + \e <\ds  4 \ii_0^1 \sqrt{V( s^2) } \, ds $ and
there exists $ A > 0 $ such that $ \psi $ is constant on $ ( - \infty, - A]$ and on $[A, \infty)$. 
By Corollary \ref{C4.5} we have $ |\psi | > 0 $ on $ \R$, hence there exist $ \theta, w \in C^{\infty}$ such that $ \psi = \rho e^{ i \theta}$, 
and $\mbox{supp}(w), \mbox{supp} ( \theta ') \subset [-A, A]$. 
We have $ p ( \theta, w ) = p + 2 k \pi$, where  $ k \in \Z$. 
If $ k = 0 $, the functions $ \psi, \; \theta $ and $ \rho $ satisfy all requirements of Step 2. 

Otherwise we construct $ \psi_{t, 1} = \rho_{t, 1} e^{ i \theta_{t, 1}}$ and $ \psi_{t, 2}=\rho_{t, 2} e^{ i \theta_{t, 2}}$ as in (\ref{lop1}) - (\ref{lop3}). Let $ w_{t, j } =  (\rho_{t, j}-1)  e^{ i \theta_{t, j}}$ for $ j = 1,2.$
By (\ref{lop5}) we see that the mappings $ m_j(t) :=  p(   \theta_{t, j}, w_{t,j}) $ are continuous, 
$ m_1 = 0 $ on $( - \infty, - A]$ and $ m_1 = 2p + 4 k \pi$ on $[A, \infty)$, 
$ m_2 = 2p + 4 k \pi $ on $( - \infty, - A]$ and $ m_2 =0 $ on $[A, \infty)$.

If $ k \ges 1 $ we may choose $ t_1 < t_2 $ such that $ m_1 ( t_1 ) = m_2( t_2) = p $. 
We have 
$$
E( \psi_{t_1, 1} ) + E( \psi _{t_2, 2 } ) = 2 \ii_{ ( -\infty, t_1 ] \cup [t_2, \infty )} |\psi '|^2 + V( |\psi |^2) \, dx \les 2 E( \psi), 
$$
hence $ E( \psi_{t_1, 1} ) \les E(\psi)$ or $ E( \psi_{t_2, 2} ) \les E(\psi)$. 
In the former case we replace $ \psi $ by $ \psi _{t_1, 1}$ and in the latter case we replace $ \psi $ by  $ \psi _{t_2, 2}$.

If $ k \les - 1 $ we have $ 2 k \pi + p \les - \pi \les - p $ because $ p \in (0, \pi]$. 
In this case by (\ref{lop5}) 
there exist $ t _1 \les t_2 $ such that $ m_1 ( t_1 ) = - p $ and $ m_2 ( t_2 ) = - p$. 
As above, we have $ E( \psi_{t_1, 1} ) + E( \psi _{t_2, 2 } ) \les 2 E( \psi)$ 
and we replace $ \psi $ by $ \ov{\psi_{t_1, 1}}$ if $ E( \psi_{t_1, 1} ) \les E(\psi)$, respectively by 
$ \ov{\psi_{t_2, 2}}$ if $ E( \psi_{t_2, 2} ) \les E(\psi)$. 

\medskip

{\it Step 3. $E_{\min}^1$ is nondecreasing on $[0, \pi]$. } 

 Let $ 0 \les p_1 < p_2 \les \pi$. 
 If $E_{\min}( p_2 ) = \ds  4 \ii_0^1 \sqrt{V( s^2) } \, ds $ we have $ E_{\min}^1 ( p_1 ) \les  E_{\min}^1 ( p_2 )$ by (vii).   
 Otherwise, consider any $ \e $ such that $ 0 < \e <  \ds  4 \ii_0^1 \sqrt{V( s^2) } -  E_{\min}^1 ( p_2 )$, then 
choose $ \psi = \rho e^{ i \theta}  = e^{ i \theta} + w$ as in Step 2.
Define $ \psi _{t, j }  = \rho_{t, j }  e^{ i \theta_{t, j }}  = e^{ i \theta_{t, j }} + w_{t, j }$ for $ j = 1, 2 $ as in (\ref{lop1}) - (\ref{lop3}). 
Using (\ref{lop5}) we see that there exist $ t_1 < t_2 $ such that $ p( \theta_{t_1, 1}, w_{t_1, 1} ) = p( \theta_{t_2, 2}, w_{t_2, 2} ) = p_1$.
By (\ref{lop4}) we have $ E^1( \psi _{t_1, 1 } ) + E^1( \psi _{t_2, 2 } ) \les 2 E^1( \psi )$. 
Then we have 
$$ E_{\min}^1 ( p_1 ) \les \min (  E^1( \psi _{t_1, 1 } ) , \,  E^1( \psi _{t_2, 2 } )) \les E^1( \psi ) \les E_{\min}^1( p_2) + \e .$$
Since $ \e $ was arbitrary we get $ E_{\min}^1 ( p_1 ) \les  E_{\min}^1( p_2)$, as desired. 

\medskip

{\it Step 4. $E_{\min}^1$ is concave. } 

Let $ 0 < p_1 < p_2 \les \pi $ and let $ p = \frac{ p_1 + p_2 }{2}$. 
We will prove that 
\beq
\label{concav1}
E_{\min}^1 \left( \frac{ p_1 + p_2}{2} \right) \ges \frac 12 E_{\min}^1 ( p_1 ) +  \frac 12 E_{\min}^1 ( p_2 ). 
\eeq
If $ E_{\min}^1 \left( \frac{ p_1 + p_2}{2} \right) = \ds  4 \ii_0^1 \sqrt{V( s^2) } \, ds $, (\ref{concav1}) obviously holds. 
Otherwise, take any $ \e $ such that  $ 0 < \e <  \ds  4 \ii_0^1 \sqrt{V( s^2) } \, ds -  E_{\min}^1 ( p)$, 
then choose  $ \psi = \rho e^{ i \theta}  = e^{ i \theta} + w$ as in Step 2.
There exists $ t_0 \in \R$ such that $ \ds \ii_{ - \infty } ^{t_0} ( 1 - \rho ^2 ) \theta '\, dx = \frac{p_1}{2}$, 
and then we have necessarily $\ds \ii_{t_0}^{ \infty }  ( 1 - \rho ^2 ) \theta '\, dx = \frac{p_2}{2}$.
Let $ \psi_{t_0, j} =   \rho_{t_0, j }  e^{ i \theta_{t_0, j }}  = e^{ i \theta_{t_0, j }} + w_{t_0, j }$ for $ j = 1,2$  be as in Step 1. 
It is clear that $p( \theta_{t_0, j } , w_{t_0, j } ) = p_j $ for $ j = 1, 2 $ and using (\ref{lop4}) we get
$$
E_{\min}^1 ( p_1 ) +   E_{\min}^1 ( p_2 ) \les E^1( \psi_{t_0, 1}) +  E^1( \psi_{t_0, 2 }) = 2 E( \psi ) \les 2 E_{\min}^1( p ) + 2 \e. 
$$
Since $ \e $ was arbitrary, we infer that (\ref{concav1}) holds. 

We have shown that $E_{\min}^1 $ is continuous and satisfies (\ref{concav1}) for any $ 0 < p_1 < p_2 \les \pi$. 
It is then standard to prove that $E_{\min}^1 $ is concave on $[0, \pi]$. 
We have $ E_{\min}^1 ( 2 \pi - p ) =  E_{\min}^1 ( - p ) =    E_{\min}^1 ( p )$ for any $ p$ and we infer that 
$ E_{\min}^1 $ is non-increasing and concave on $[\pi, 2 \pi]$, and then it follows  that it is concave on $[0, 2 \pi]$. 
\hfill
$\Box $

\begin{Theorem}
\label{T5.2}
Assume that conditions {\bf (A1)} and {\bf (B1)} in the introduction hold and $ p \in (0, \pi]$ satisfies $ E_{\min}^1 ( p ) < \sqrt{2} p $. 
Let $ ( \psi_n )_{n \ges 1} \subset \Er $ be a sequence satisfying 
\beq
\label{conv0}
\PR ( \psi _n ) \lra \pr{p} \qquad \mbox{ and } \qquad E ( \psi _n ) \lra E_{\min}^1 ( p) \qquad \mbox{ as } n \lra \infty. 
\eeq
Then there exist a subsequence $(\psi_{n_k})_{k \ges 1}$, a sequence $ (x_k )_{k \ges 1} \subset \R $ and $ \psi \in \Er$ satisfying 
$\PR ( \psi ) = \pr{p}$, $E^1( \psi ) = E_{\min}( p )$, and
$$
\begin{array}{l}
\psi_{n_k} ( \cdot + x_k ) \lra \psi \quad \mbox{ uniformly on } [ - R, R] \mbox{ for any } R > 0, 
\\
\psi_{n_k}' ( \cdot + x_k ) \lra \psi ' \quad \mbox{ in } L^2 ( \R), 
\\
| \psi_{n_k} ( \cdot + x_k ) | - |\psi | \lra 0  \quad \mbox{ in } L^p ( \R) \mbox{ for } 2 \les p < \infty.  
\end{array}
$$

\end{Theorem}

{\it Proof. } 
Let $ p \in (0, \pi] $ such that $ E_{\min}^1 ( p ) < \sqrt{2} p $ and let $ ( \psi _n ) _{n \ges 1 } \subset \Er $ be a sequence 
satisfying (\ref{conv0}). 
We denote by $ p_n \in [0, 2 \pi)$ the canonical valuation of the momentum of $ \psi _n $. Then we have $ p_n \lra p $ as $ n \lra \infty$. 

Denoting $ M := \sup_{n \ges 1} \|\psi _n ' \|_{L^2 ( \R ) } < \infty $, 
by  (\ref{holder}) we have 
\beq
\label{holder1}
| \psi_n ( x ) - \psi _n ( y )| \les M | x - y |^{\frac12} \qquad \mbox{  for all } x, y \in \R \mbox{  and all } n \in \N^*. 
\eeq

We  use the  concentration-compactness principle (see \cite{lions}). 
The sequence $ f_n := |\psi _n '|^2 + V( |\psi _n |^2) $ is bounded in $ L^1( \R)$. We denote by $ \Lambda _n $ the concentration function of $ f_n$, namely 
$$
\Lambda _n ( t ) = \sup_{x \in \R} \ii_{x - t }^{ x + t }  |\psi _n ' ( y )|^2 + V( |\psi _n (y)|^2) \, dy. 
$$
Obviously, $ \Lambda _n $ is a non-decreasing function on $[0, \infty)$, 
$ \Lambda _n ( 0 ) = 0 $ and $ \Lambda _n ( t ) \lra E^1 ( \psi _n ) $ as $ t \lra \infty $ . 
Proceeding as in \cite{lions} we see that there exists a subsequence of $( \psi_n , \Lambda _n )_{n \ges 1}$, still denoted $( \psi_n , \Lambda _n )_{n \ges 1}$,
and there is a non-decreasing function $ \Lambda : [0, \infty ) \lra \infty $ satisfying 
\beq
\label{cc1}
\Lambda _n (t) \lra \Lambda (t) \mbox{ a.e on } [0, \infty) \mbox{ as } n \lra \infty .
\eeq
For any fixed $ t $ we have $ \Lambda _n ( t ) \les E^1( \psi _n ) $ and letting $ n \lra \infty $ we get $ \Lambda (t ) \les E_{\min}^1 ( p ) $. 

\smallskip

Let $ \al = \ds \lim_{ t \ra  \infty } \Lambda ( t)$. 
It is clear that $ \ds 0 \les \al \les  
E_{\min}^1 ( p ) $.
We will show that $ \al = E_{\min}^1 ( p ) $. 

We prove first that $ \al > 0 $. We argue by contradiction and we assume that $ \al = 0$. 
This implies that $ \Lambda ( t ) = 0 $ for any $ t > 0$, which means that  for any fixed $ t > 0 $ we have 
\beq
\label{van1}
 \sup_{x \in \R} \ii_{x - t }^{ x + t }  |\psi _n ' ( y )|^2 + V( |\psi _n (y)|^2) \, dy \lra 0 \qquad \mbox{ as } n \lra \infty. 
\eeq
We claim that if (\ref{van1}) occurs  then necessarily   $ |\psi _n | \lra 1 $ uniformly on $ \R$. 
By {\bf (A1)} there is $ \eta _0 > 0 $ such that 
$$
 \frac 14 ( 1 - s^2 )^2 < V( s^2 ) < ( 1 - s^2 )^2  \qquad \mbox{ for all } s \in [1 - \eta_0, 1 + \eta _0]. 
$$
Fix $ \eta \in (0, \frac{\eta_0}{2} ]$. 
Assume that there exists $ x_n \in \R$ such that $ \big| \, |\psi _n ( x_n ) | - 1 \big|   =  \eta $. 
From (\ref{holder1}) it follows that there  exists $ r > 0 $, independent of $n$, such that 
$  \big| \, |\psi _n ( x_n ) | - 1 \big|  \in [ \frac{\eta}{2},  \frac{3\eta}{2}]$ 
for $ y \in [ x_n - r, x_n + r]$ and therefore
$$
\Lambda_n ( r ) \ges \ii_{x_n - r } ^{ x_n + r } V ( |\psi _n (y) |^2 ) \, dy \ges \ii_{x_n - r } ^{ x_n + r } \frac{1}{4}  \big| \, |\psi _n ( x_n ) | - 1 \big|^2  \, dy 
\ges \frac 18 \eta ^2 r. 
$$
By  (\ref{van1})  there exists $ n_{\eta } \in \N$ such that for all $ n \ges n_{\eta} $ we  have 
$ \Lambda_n ( r ) < \frac 18 \eta^2 r$, and the above inequalities imply that for any  $n\ges n_{\eta}$ we must have  
$ 1 - \eta < |\psi _n | < 1 + \eta$. Thus $ |\psi _n | \lra 1 $ uniformly on $ \R$. 

Choose $ \de \in (0, 1) $ such that $ E_{\min}^1 ( p ) < ( 1 - \de )^2 \sqrt{2} p $ (this is possible because $ E_{\min}^1 ( p ) <  \sqrt{2} p $).
By {\bf (A1)} there is $ \eta_{\de}> 0 $ such that $ V( s^2 ) \ges \frac 12 ( 1 - \de )^2 ( 1 - s ^2) ^2 $ for any $ s \in [1 - \eta_{\de}, 1 + \eta_{\de}]$. 
For all $ n $  sufficiently large  we have $ \big| \, |\psi _n | - 1 \big| < \min ( \de, \eta_{\de} )$ on $ \R$. 
For any such $n$ we may write $ \psi _n = \rho _n e^{ i \theta_n}$ where $ 1 - \rho _n \in H^1( \R)$ 
and $ \theta_n \in \dot{H}^1( \R, \R)$,  and we have 
\beq
\label{van2}
\begin{array}{l}
\ds E^1( \psi _n) \ges \ii_{\R} \rho _n ^2 |\theta_n '|^2 + V( \rho _n ^2 ) \, dx 
\ges ( 1 - \de )^2 \ii_{\R} |\theta_n '|^2 + \frac 12 ( 1 - \rho _n ^2 ) ^2 \, dx 
\\
\\
\ds \ges  ( 1 - \de )^2 \sqrt{2} \bigg| \ii_{\R}  ( 1 - \rho _n ^2 ) \theta_n '\, dx \bigg| .
\end{array}
\eeq
Recall that  $\ds  \ii_{\R}  ( 1 - \rho _n ^2 ) \theta_n '\, dx $ is a valuation of the momentum of $ \psi _n$ by Remark \ref{lift1}, hence 
there exists $ \ell_n \in \Z$ such that $ \ds  \ii_{\R}  ( 1 - \rho _n ^2 ) \theta_n '\, dx  = p _n + 2 \ell_n \pi$,  and consequently 
$$
 \bigg| \ii_{\R}  ( 1 - \rho _n ^2 ) \theta_n '\, dx \bigg| \ges \mbox{dist}( p_n, 2 \pi \Z). 
$$
Letting $ n \lra \infty $ in (\ref{van2}) we get 
$$
E_{\min}^1 ( p ) \ges ( 1 - \de)^2 \sqrt{2} \lim_{n \ra \infty }  \mbox{dist}( p_n, 2 \pi \Z) = ( 1 - \de)^2 \sqrt{2} p, 
$$
contradicting the choice of $ \de$. We have thus proved that $ \al \neq 0 $. 

\medskip

Assume that $ 0 < \al < E_{\min}^1 ( p )$. 
Arguing as in the proof of Theorem 5.3 in \cite{M10} (see (5.12) p. 156 there) we infer that there is a nondecreasing sequence $ R_n \lra \infty $ such that 
$$
\lim_{n \ra \infty } \Lambda_n ( 2 R_n ) = \lim_{n \ra \infty } \Lambda_n (  R_n ) =  \al . 
$$
For each $ n $ choose $ x_n \in \R$ such that 
$$
 \ii_{x_n  - R_n }^{ x_n + R_n }  |\psi _n ' ( y )|^2 + V( |\psi _n (y)|^2) \, dy > \Lambda_n ( R_n) - \frac 1n. 
$$
Then we have 
$$
\Lambda_n ( R_n) - \frac 1n 
< \ii_{x_n  - 2R_n }^{ x_n + 2R_n }  |\psi _n ' ( y )|^2 + V( |\psi _n (y)|^2) \, dy \les \Lambda_n ( 2 R_n), 
$$
and we infer that 
\beq
\label{dic1}
 \ii_{x_n  - R_n }^{ x_n + R_n }  |\psi _n ' ( y )|^2 + V( |\psi _n (y)|^2) \, dy \lra \al, 
\eeq
\beq
\label{dic2}
  \ii_{x_n  - 2R_n }^{ x_n - R_n}  |\psi _n ' ( y )|^2 + V( |\psi _n (y)|^2) \, dy
  + \ii_{ x_n + R_n}^{x_n + 2 R_n}   |\psi _n ' ( y )|^2 + V( |\psi _n (y)|^2) \, dy \lra 0 , \qquad \mbox{ and } 
\eeq
\beq
\label{dic3}
\ii_{- \infty }^{x_n - 2 R_n }  \! \!  |\psi _n ' ( y )|^2 \! + V( |\psi _n (y)|^2) \, dy  +  \ii_{x_n + 2 R_n }^{\infty}   \! \!  |\psi _n ' ( y )|^2 \! + V( |\psi _n (y)|^2) \, dy \lra E_{\min}^1( p ) - \al . 
\eeq
Assume that for infinitely many $n$'s we have 
$$\ds \ii_{- \infty }^{x_n - 2 R_n }  \! \!  |\psi _n ' ( y )|^2 \! + V( |\psi _n (y)|^2) \, dy  \ges   \ii_{x_n + 2 R_n }^{\infty}   \! \!  |\psi _n ' ( y )|^2 \! + V( |\psi _n (y)|^2) \, dy . $$
(A similar argument will work  if the opposite inequality holds true for  infinitely many $n$'s.)
Passing to a subsequence, still denoted the same, we may assume that 
\beq
\label{dic4}
\ii_{- \infty }^{x_n - 2 R_n }  \! \!  |\psi _n ' ( y )|^2 \! + V( |\psi _n (y)|^2) \, dy \lra \beta > 0 \quad \mbox{ as } n \lra \infty, 
\quad \mbox{ where  } 0 < \beta < E_{\min}^1 ( p ) . 
\eeq
Let $ z_n  = x_n - \frac 32 R_n$. 
For $ n $ sufficiently large we have $ [z_n - 2, z_n + 2 ] \subset [ x_n - 2R_n, x_n - R_n]$. Then using (\ref{holder1}) and (\ref{dic2}) and arguing as in the proof of the fact that $ \al > 0 $ we infer that $ \ds \sup_{y \in [z_n - 1, z_n + 1 ]} \big| \, |\psi _n ( y ) | - 1 \big| \lra 0 $ as $ n \lra \infty$. 
For $ n $ large enough we have $ \frac 12 < |\psi _n | < \frac 32 $ on $  [z_n - 1, z_n + 1 ]$, thus we have a lifting 
$ \psi _n = \rho _n e^{ i \theta_n }$ on that interval. Let $ r_n = \rho _ n ( z_n)$ and $ \al _ n = \theta _ n ( z_n )$. 
It is clear that $ r_n \lra 1$. 
Define 
$$
\psi _{n , 1}( x ) = \left\{
\begin{array}{l}
\psi _n ( x ) \mbox{ if } x \les z_n, 
\\
\big( ( 1 - r_n ) ( x - z_n ) + r_n \big) e^{ i \al _n } \mbox{ if }  x \in (z_n,  z_n + 1), 
\\
e^{ i \al _n }  \mbox{ if }  x \ges z_n + 1,  
\end{array}
\right. 
$$
$$
\psi _{n , 2}( x ) = \left\{
\begin{array}{l}
e^{ i \al _n }  \mbox{ if } x \les z_n -1, 
\\
\big( (  r_n -1 ) ( x - z_n ) + r_n \big) e^{ i \al _n } \mbox{ if }  x \in (z_n-1 ,  z_n ), 
\\
\psi _n ( x )  \mbox{ if }  x \ges z_n . 
\end{array}
\right. 
$$
It is  easy to see that $ \psi_{n, 1} , \psi_{n , 2} \in \Er $. It follows form (\ref{dic2}) and (\ref{dic4}) that  $ E^1( \psi_{n , 1}) \lra \beta$ 
and  $  E^1( \psi_{n , 2}) \lra E_{\min}^1( p ) - \beta$ as $ n \lra \infty. $
By Lemma \ref{E1min} (vii) we have $ 0 < \beta <  E_{\min}^1( p ) \les  4 \ii_0^1 \sqrt{V( s^2) } \, ds $, hence 
 $  E^1( \psi_{n , j }) <  4 \ii_0^1 \sqrt{V( s^2) } \, ds $  for all $ n $ sufficiently large and  $ j = 1, 2 $. 
Then using Corollary \ref{C4.5} we infer that $ \psi_{n, 1}$ and $ \psi_{n, 2}$ do not vanish and consequently these functions admit 
liftings, thus we may write 
$ \psi_{n , j} = \rho_{n , j } e^{ i \theta _{n , j }} $ on $ \R$. 
Replacing $ \theta_{n, j }$ by  $ \theta_{n, j } + 2 k \pi $ for some $ k \in \Z$, we may assume that 
$ \theta_{n , 1 } = \al _ n $ on $[z_n, \infty)$ and that $ \theta_{n , 2 } = \al _ n $ on $( - \infty, z_n ]$. 
Let 
$$
\rho _ n (x) = \left\{ 
\begin{array}{ll}
\rho_{n , 1 } (x )& \mbox{ if } x \les z_n, \\
\rho_{n , 2 } (x )& \mbox{ if } x > z_n,
\end{array}
\right.
\qquad
\theta _ n (x) = \left\{ 
\begin{array}{ll}
\theta_{n , 1 } (x )& \mbox{ if } x \les z_n, \\
\theta_{n , 2 } (x )& \mbox{ if } x > z_n.
\end{array}
\right.
$$
Then $ 1 - \rho _ n \in H^1( \R)$, $ \theta_n \in \dot{H}^1( \R)$ and $ \psi _ n = \rho_n e^{ i \theta _ n }$ on $ \R$. 
By Remark \ref{lift1} and the fact that $ \theta_{n, 1}' = 0 $ on $( z_n, \infty)$, $ \theta_{n, 2}' = 0 $ on $(- \infty,  z_n)$, 
 valuations of momenta of $ \psi_n $, $ \psi_{n, 1 }$ and $ \psi_{n, 2 }$ are, respectively, 
 \beq
\label{p-val}
{p}_n = \ii_{\R} ( 1 - \rho _n ^2 ) \theta_n ' \, dx, 
\quad
{p}_{n, 1} = \ii_{- \infty }^{ z_n } ( 1 - \rho _n ^2 ) \theta_n ' \, dx, 
\quad
{p}_{n, 2} = \ii_{z_n}^{ \infty } ( 1 - \rho _n ^2 ) \theta_n ' \, dx.
\eeq
Passing to a subsequence, we may assume that $ | \pr{p_{n,1}} | \lra p_1 ^* \in [0, \pi] $ and 
 $ | \pr{p_{n,2}} | \lra p_2 ^* \in [0,\pi] $  as $ n \lra \infty$. 
 From (\ref{p-val}) we have $ p_{n, 1 } + p_{n,2} = p_n$. 
 Letting $ n \lra \infty $ and using the inequality $|\pr{a + b } | \les |\pr{a}| + |\pr{b}|$, we get 
$$
p_1 ^* + p_2 ^* \ges p . 
$$
We have 
$ E_{\min}^1 ( |\pr{p_{n, j }} |) =  E_{\min}^1 ( p_{n, j } )  \les E^1( \psi _{n, j})$. 
Passing to the limit,  we find 
$$
 E_{\min }^1 ( p_1^* ) \les \beta <  E_{\min }^1(p)  \qquad \mbox{ and } \qquad  E_{\min }^1 ( p_2^* ) \les p - \beta <  E_{\min }^1 ( p) .
 $$
Since $ E_{\min }^1$ is nondecreasing on $[0, \pi]$,  we have $ p_1^* < p $ and $ p_ 2^* < p $. Then the inequality $
p_1 ^* + p_2 ^* \ges p 
$ implies that $ 0 < p_j ^* < p $ for $ j = 1,2$. 
We have
\beq
\label{subaddit1}
E_{\min}^1 ( p_{n, 1} ) + E_{\min}^1 ( p_{n, 2} ) \les E^1 ( \psi_{n, 1} ) +  E^1 ( \psi_{n, 2} )  = E^1 ( \psi _n ) + o(1)
\eeq
and letting $ n \lra \infty $ we discover $ E_{\min}^1 ( p_1 ^* ) + E_{\min}^1 ( p_2 ^* ) \les E_{\min}^1 ( p )$. 
The concavity of $E_{\min}^1 $ on $[0, 2 \pi]$ implies that 
$ E_{\min}^1 ( p_ j  ^* ) \ges \frac{p_j ^*}{p } E_{\min}^1 ( p )$, and equality may occur if and only if $E_{\min}^1 $ is linear on $[0, p]$. 
Summing up these inequalities for $ j = 1, 2 $ and comparing to the previous inequality we infer that we must have equality, and 
consequently $  E_{\min}^1 $ must be  linear on $[0, p ]$. 
Taking into account the behaviour of $E_{\min}^1 $ at the origin (see Lemma \ref{E1min} (ii) and (iii)) we infer that 
$ E_{\min}^1 ( s ) = \sqrt{2} s $ for all $ s \in [0, p]$,  contradicting the fact that $ E_{\min}^1 ( p ) <  \sqrt{2} p $.  

We conclude that we cannot have $ 0 < \al  < E_{\min}^1 ( p )$, thus necessarily $ \al =  E_{\min}^1 ( p )$.

\medskip

It is then standard to prove that there exists a sequence $ ( x_n ) _{n \ges 1 } \subset \R^N$ such that for any $ \e > 0 $ 
there exists $ R_{\e } > 0 $ satisfying 
$$ \ds \ii_{ - \infty }^{ x_n - R_{\e}} |\psi _n '|^2  + V( |\psi _n |^2 ) \, dx +  \ii_{ x_n + R_{\e}}^{  \infty } |\psi _n '|^2  + V( |\psi _n |^2 ) \, dx < \e $$
 for all sufficiently large $ n $. 
Let $ \tilde{\psi }_n = \psi _n ( \cdot + x_n )$. Then for any $ \e > 0 $ there exist $ R_{\e } > 0 $ and $ n_{\e } \in \N$ such that 
\beq
\label{conc1}
 \ii_{ - \infty }^{  - R_{\e}} |\tilde{\psi} _n '|^2  + V( |\tilde{\psi } _n |^2 ) \, dx 
+  \ii_{ R_{\e}}^{  \infty } |\tilde{\psi }_n '|^2  + V( |\tilde{\psi } _n |^2 ) \, dx < \e \qquad \mbox{ for all } n \ges n_{\e}. 
\eeq

Take $ \eta _0 > 0 $ such that 
$ \frac 14 ( 1 - s^2 )^2 < V( s^2 ) < ( 1 - s^2 )^2 $ for any $ s \in [1 - \eta_0, 1 + \eta _0].$ 
Let $H $ be as in assumption {\bf (B1)}. 
Let $ \e _1 = 2 \min (\big| H( \frac 12 ) \big|, H( \frac 32 ), \big| H( 1 - \eta_0 ) \big| , H( 1 + \eta_0))$. Take $ R_1 = R_{\e _1} > 0 $ 
such that (\ref{conc1}) holds with $ \e_1 $ instead of $ \e$. 
Using (\ref{escape}) we see that for all $ n $ sufficiently large and for all $ x \in ( - \infty, - R_1 ] \cup [R_1, \infty)$ we have 
\beq
\label{conc1bis}
|\tilde{\psi}_n ( x )| \in \left[\frac 12, \frac 32 \right] \cap [1 - \eta_0, 1 + \eta_0] \;   \mbox{ and } \; 
\frac 14 ( 1 - |\tilde{\psi}_n ( x )|^2 )^2 < V( |\tilde{\psi}_n ( x )|^2 ) < ( 1 -|\tilde{\psi}_n ( x )|^2 )^2.
\eeq
Then using (\ref{holder1}) we infer  that $ \tilde{\psi}_n $ are uniformly bounded on $ [ - R_1, R_1]$, hence on $ \R$. 
Since $ (\tilde{\psi}_n ')_{n \ges 1} $ is bounded in $ L^2( \R)$, 
it is standard to prove that there exists a function $ \psi \in H_{loc}^1 ( \R ) $ such that $ \psi ' \in L^2 ( \R)$ and there is a subsequence 
$(\tilde{\psi}_{n _k})_{k \ges 1}$  of $ (\tilde{\psi}_n)_{n \ges 1}$ satisfying 
\beq
\label{conc2}
\begin{array}{l}
\tilde{\psi}_{n_k} ' \rightharpoonup \psi ' \qquad \mbox{ weakly in } L^2 ( \R), 
\\
\tilde{\psi}_{n_k}  \lra \psi  \qquad \mbox{ strongly in } L^p ([-R, R] ) \mbox{ for any } R > 0 \mbox{ and any } 1 \les p \les \infty, 
\\
\tilde{\psi}_{n_k}  \lra \psi  \qquad \mbox{uniformly on  } [-R, R] \mbox{ for any } R > 0 .
\end{array}
\eeq
We may use (\ref{holder1}) and  Arzelà-Ascoli's  Theorem to get uniform convergence on compact intervals.
The weak convergence $ \tilde{\psi}_{n_k} ' \rightharpoonup \psi ' $ in $ L^2( \R)$ implies that for any interval $ I \subset \R$ we have
\beq
\label{conc3}
\ii_I |\psi '|^2 \, dx \les \liminf_{n \ra \infty } \ii_{I} |\tilde{\psi}_{n_k } ' |^2 \, dx. 
\eeq
Fatou's Lemma gives 
\beq
\label{conc4}
\ii_I V( |\psi |^2 ) \, dx \les  \liminf_{n \ra \infty } \ii_{I} V \left( |\tilde{\psi}_{n_k }  |^2  \right) \, dx. 
\eeq

Fix $ \e > 0 $. Take $  R_{\e}$ as in (\ref{conc1}). 
Since $ ( \tilde{\psi}_{n_k} )_{k \ges 1}$ is bounded in $ L^{\infty} ([-R_{\e}, R_{\e}])$, converges to $ \psi $ almost everywhere and 
$V$ is continuous, the dominated convergence theorem gives 
$
 V ( |\tilde{\psi}_{n_k}|^2)  \lra  V ( |\psi |^2 )  
$
in $ L^1 ( [-R_{\e}, R_{\e} ])$, 
hence  there is $ k_{\e}$ such that $ n_{k_{\e}} \ges n_{\e }$ and
$$
\ds \ii_{- R_{\e}}^{R_{\e}} \big| V ( |\tilde{\psi}_{n_k}|^2) -  V ( |\psi |^2 ) \big| \, dx < \e \qquad \mbox{ for all } k \ges k_{\e }.
$$
Then using (\ref{conc1}) and the similar estimate for $ \psi $  we have for all $ k \ges k_{\e}$, 
$$
\begin{array}{l}
\ds \ii_{\R} \big| V ( |\tilde{\psi}_{n_k}|^2) -  V ( |\psi |^2 ) \big| \, dx 
\\
\\
\les 
\ds \ii_{- R_{\e}}^{R_{\e}} \big| V ( |\tilde{\psi}_{n_k}|^2) -  V ( |\psi |^2 ) \big| \, dx
+ \ii_{( - \infty, - R_{\e}] \cup [ R_{\e}, \infty )}  V ( |\tilde{\psi}_{n_k}|^2)  +  V ( |\psi |^2 )\, dx \les 3\e . 
\end{array}
$$
Since $ \e $ was arbitrary we get $ V ( |\tilde{\psi}_{n_k}|^2) \lra   V ( |\psi |^2 )$ in $ L^1( \R)$, and in particular 
\beq
\label{conc4bis}
\ds \ii_{\R}  V ( |\tilde{\psi}_{n_k}|^2)\, dx  \lra  \ii_{\R}  V ( |\psi |^2 )\, dx \qquad \mbox{ as }  k\lra \infty.
\eeq
Similarly we show that $ \left(1 - |\tilde{\psi}_{n_k}|^2  \right) \lra \left(1 - |\psi |^2\right)  $ strongly in $ L^p ( \R)$ for any $ p \in [2, \infty)$. 

\medskip

Next we show that the momenta of $ \tilde{\psi}_{n_k}$ converge to the momentum of $ \psi $. To do this we find convenient representations of $ \tilde{\psi}_{n_k}$ in the form $ \tilde{\psi}_{n_k} = e^{ i \theta _k } + w_k$, and we do the same for $ \psi $. 

On $ ( - \infty, - R_1] $ and on $  [R_1, \infty )$ we have liftings, that is we may write $  \tilde{\psi}_{n_k} = \rho_k e^{ i \theta _k}$ and 
$ \psi = \rho e^{ i \theta}$. 
Given $ \e > 0 $, take  $ R_{\e } > R_1$ such that (\ref{conc1}) holds. 
Since $ \frac{\tilde{\psi}_{n_k}}{\psi }( \pm R_{\e}) \lra 1 $ as $ k \lra \infty$, 
we may replace if necessary $ \theta_k $ by $ \theta _k + 2 \ell_k \pi $ for some $ \ell_k \in \Z$ on $ ( - \infty, - R_1] $ and/or on $  [R_1, \infty )$ 
in such a way that $ \theta _k ( - R_{\e} ) = \theta ( - R_{\e}) + \al_k ^{ - }$ 
and $ \theta _k (  R_{\e} ) = \theta (  R_{\e}) + \al_k ^{ + }$, where $ \al_k ^{ \pm } \lra 0 $ as $ k \lra \infty$. 
 Then we extend $ \theta $ and $ \theta_k $ as affine functions on $ [ - R_{\e}, R_{\e}]$. 
 It is easily seen that $ e^{ i \theta _k } \lra e^{ i \theta} $ uniformly on  $ [ - R_{\e}, R_{\e}]$ and in $H^1(  [ - R_{\e}, R_{\e}] )$, 
 and that $ \theta_k, \theta \in \dot{H}^1( \R)$. 
 Denote $ w_k = \tilde{\psi}_{n_k} - e^{ i \theta_k }  $ and $ w = \psi - e^{ i \theta}$. 
From (\ref{conc2}) it follows  that $ w_k, w \in H^1( \R)$, $ w_k ' \rightharpoonup w'$ weakly in $ L^2 (  [ - R_{\e}, R_{\e}] )$, 
and $ w_k \lra w $ strongly in $L^2  (  [ - R_{\e}, R_{\e}] )$.
Then we infer that 
\beq
\label{conc6}
\ii_{- R_{\e}}^{ R_{\e}} -2 \langle \theta_k ' e^{i \theta _k } , w_k \rangle + \langle i w_k ', w_k \rangle \, dx 
\lra \ii_{- R_{\e}}^{ R_{\e}} -2 \langle \theta ' e^{i \theta } , w \rangle + \langle i w ', w \rangle \, dx 
\quad \mbox{ as } n \lra \infty. 
\eeq
On the other hand, on $ ( - \infty,- R_1 ] \cup [R_1 , \infty ) $ we have 
$  -2 \langle \theta_k ' e^{i \theta _k } , w_k \rangle + \langle i w_k ', w_k \rangle = ( 1 - \rho_k ^2 ) \theta _ k ' $.
Using (\ref{conc1bis}), proceeding as in (\ref{van2}), then using (\ref{conc1}) we get 
$$
\bigg| \ii_{- \infty } ^{-R_{\e}}  ( 1 - \rho_k ^2 ) \theta _ k ' \, dx \bigg|
\les 2\ii_{- \infty } ^{-R_{\e}} \frac 14 |\theta_k ' |^2 +  \frac 14 ( 1 - \rho_k ^2 )^2 \, dx 
\les  2\ii_{- \infty } ^{-R_{\e}} \rho _k ^2 |\theta_k ' |^2 +  V( \rho _k ^2 ) \, dx 
\les 2 \e. 
$$
A similar estimate holds true on $ [R_{\e} , \infty)$, 
as well as  for the function $ \psi $. 
Thus we get 
\beq
\label{conc8}
\begin{array}{l}
\ds \bigg| \ii_{\R}  -2 \langle \theta_k ' e^{i \theta _k } , w_k \rangle + \langle i w_k ', w_k \rangle \, dx - 
 \ii_{\R}  -2 \langle \theta ' e^{i \theta } , w \rangle + \langle i w ', w \rangle \, dx 
\bigg|
\\
\\
\ds \les 8 \e + \bigg| \ii_{- R_{\e}}^{R_{\e}}  -2 \langle \theta_k ' e^{i \theta _k } , w_k \rangle + \langle i w_k ', w_k \rangle \, dx - 
 \ii_{\R}  -2 \langle \theta ' e^{i \theta } , w \rangle + \langle i w ', w \rangle \, dx \bigg|
 \end{array}
\eeq
for all $k $ sufficiently large. 
Then using (\ref{conc6})   we see that the right-hand side of (\ref{conc8}) is smaller than $ 9 \e $ if $k  $ is large enough. 
Since $ \e $ was arbitrary, we have proved that 
\beq
\label{conc9}
p( \theta_k, w_k ) \lra p( \theta, w ) \qquad \mbox{ as } k \lra \infty. 
\eeq
We have $ \pr{p( \theta_k, w_k )} = \PR ( \tilde{\psi}_{n_k}) = \pr{p_{n_k}} \lra \pr{p}$, hence $ \pr{p ( \theta, w )} = \pr{p}$
and consequently  
$$ E^1( \psi ) \ges E_{\min}^1 ( p( \theta, w)) 
= 
 E_{\min}^1 ( p ) = \lim_{k \ra \infty } E^1 ( \tilde{\psi}_{n_k}). $$
On the other hand, we have $ \ds \ii_{\R} |\psi '|^2 \, dx \les \liminf_{k \ra \infty } \ii_{\R} | \tilde{\psi}_{n_k} ' |^2 \, dx $ because 
$ \tilde{\psi}_{n_k} '  \rightharpoonup \psi '$ in $ L^2( \R)$. 
Taking into account (\ref{conc4bis}), we infer that necessarily 
$ \| \tilde{\psi}_{n_k} ' \|_{L^2 ( \R)} ^2 \lra \| \psi ' \|_{L^2}^2 $ as $ k \lra \infty$. 
The weak convergence and the convergence of norms imply that  $ \tilde{\psi}_{n_k} ' \lra \psi '$ strongly in $ L^2( \R)$. 
\hfill
$\Box$

\begin{Proposition}
\label{P5.3}
Let $ p \in (0, \pi]$ and let $ \psi \in \Er $ be a solution of  the minimization problem considered above, that is 
$ \PR( \psi ) = \pr{p}$ and $ E^1 ( \psi ) = E_{\min}^1 ( p )$. 
Then there exists $c \in \left[ (E_{\min}^1)_r ' ( p ) ,  (E_{\min}^1)_{\ell} '(p)\right]$, where $  (E_{\min}^1)_{\ell} '(p)$ 
and $(E_{\min}^1)_r '(p)$ are the 
left and right derivatives of $ E_{\min}^1$ at $p$, respectively,  such that 
$$
i c \psi ' + \psi '' + F( |\psi |^2 ) \psi = 0 
\qquad \mbox{ in } \R. 
$$
In other words, $ \psi $ is a one-dimensional traveling wave of speed $c$ for (\ref{1.1}), and $ \psi \in C^2 ( \R)$. 

Moreover, for any $ p \in (0, \pi]$ such that $ E_{\min}^1 ( p ) < \sqrt{2} p $ 
and $ (E_{\min}^1 )_{\ell } ' ( p ) > (E_{\min}^1 )_{r } ' ( p )$, 
there exist two minimizers $ \psi _1 , \psi _2 \in \Er $ for $ E_{\min}^1 ( p ) $ that solve (\ref{4.1}) with speeds 
$ c _ 1 = (E_{\min}^1 )_{\ell } ' ( p )$ and $ c _ 2 = (E_{\min}^1 )_{r } ' ( p )$, respectively.

\end{Proposition}

The proof of Proposition \ref{P5.3} is standard and is similar to the proof of Proposition 4.14 in \cite{CM}, so we omit it.

\begin{remark} 
\label{R5.4}
\rm 
i) Let 
$$ \mathscr{C} = \left\{ c \in [0, \sqrt{2}) \; \Big| \; \begin{array}{l} \mbox{ there exist } p \in (0, \pi] \mbox{ and a minimizer } \psi \in \Er \mbox{ of } E_{\min}^1( p ) \\ \mbox{ that solves (\ref{4.1}) with speed } c \end{array} \right\}.
$$
The set $ \mathscr{C}$ is relatively closed in $ [0, \sqrt{2})$. 
Indeed, assume that $ c_n \in \mathscr{C}$ and  $ c_n \lra c_* \in [0, \sqrt{2})$. 
For each $n$ choose $ p_n \in  (0, \pi]$ and $ \psi _n \in \Er$ such that $ \psi _n $ is a minimizer for $ E_{\min}^1 ( p_n)$ and solves (\ref{4.1}) with speed $  c_n $. 
Since $E_{\min}^1 $ is concave and  $ {\ds \lim_{p \ra 0 }} \frac{E_{\min}^1 ( p)}{p } = \sqrt{2}$, we have $  {\ds \lim_{p \ra 0 }} \left(E_{\min}^1 \right)_{\ell}' ( p) =  {\ds \lim_{p \ra 0 }} \left(E_{\min}^1 \right)_{r}' ( p) = \sqrt{2}$. 
Since $c_n \lra c_* < \sqrt{2}$, we  have necessarily $ \ds \liminf_{n \ra \infty } p_n > 0 $. 
Passing to a subsequence we may assume that $ p_n \lra p _* \in (0, \pi]$. 
Then there exist a subsequence $(\psi_{n_k})_{k \ges 1}$, a sequence of points $ (x_k ) \subset \R$ and $ \psi \in \Er $ satisfying the conclusion of Theorem \ref{T5.2}.
The function $ \psi $ is a minimizer for $E_{\min}^1 ( p ) $ and solves (\ref{4.1}) for some $ c \in \left[ (E_{\min}^1)_r ' ( p ) ,  (E_{\min}^1)_{\ell} '(p)\right]$. 
On the other hand, writing the equation (\ref{4.1}) satisfied by each $ \psi _{n_k}$ and letting $ k \lra \infty$ we see that $ \psi $ satisfies (\ref{4.1}) in $ \Do '( \R)$ with speed $  c_*$. Hence $ c_* = c \in \mathscr{C}$. 

\smallskip

ii) If $F$ is locally Lipschitz, $ E_{\min}^1 $ cannot be affine on any interval $ [p_1, p_2] \subset [0, \pi]$ such that $ p_1 < p_2 $ and $ E_{\min}^1 ( p_2 ) < \sqrt{2} p_2$. 
To see this we argue by contradiction and we assume that $ E_{\min}^1 $ is affine on an interval as above. 
Then we have $  E_{\min}^1 ( p ) < \sqrt{2} p $ for any $ p \in ( p_1, p_2]$, and Theorem \ref{T5.2} implies that there exist minimizers for $ E_{\min}^1 (p) $ for all 
$ p \in (p_1, p_2]$. 
All these minimizers satisfy equation (\ref{4.1}) with the same speed $c$, namely the slope of $ E_{\min}^1 $ on the interval $ [p_1, p_2]$, but have different momenta. 
This contradicts the uniqueness of solutions of  (\ref{4.1}) modulo translations and multiplication by complex numbers of modulus one (see Proposition \ref{P4.1} (v)). 

If  there is some $ p_0 \in (0, \pi]$ such that $ E_{\min}^1 ( p_0 ) < \sqrt{2} p$, 
 the continuity of $ E_{\min}^1$ implies that this inequality holds for all $ p $ in an interval $I \subset (0, \pi]$ containing $p_0$. Then Theorem \ref{T5.2} gives the existence of minimizers for 
 $ E_{\min}^1( p)$ for all $ \tilde{p} \in I$. 
If $ p_1, p_2 \in I$, $ p_1 < p_2$ and $ \psi_1$, $ \psi _2$ are minimizers for 
$ E_{\min}^1( p_1)$ and $ E_{\min}^1( p_2)$ with associate Lagrange multipliers $ c( \psi_1)$ and $ c( \psi _2),$  respectively, Proposition \ref{P5.3} and the fact that $  E_{\min}^1$ is not affine on $ [p_1, p_2]$ imply that 
\beq
\label{Lag-monotone}
c( \psi _1 ) \ges  (E_{\min}^1 )_{r } ' ( p _1) >  (E_{\min}^1 )_{\ell } ' ( p_2 ) \ges c ( \psi _2). 
\eeq
In particular, the set $ \mathscr{C}$ is uncountable.

\end{remark}

\begin{example} 
\label{Example1}
\rm
Consider $ V \in C^{\infty }([0, \infty))$ such that $ V( s ) = \frac 12 ( 1 - s )^2 $ on $ [1 - \de , \infty)$ for some $ \de > 0 $, 
$V$ is decreasing on $ [0, 1)$ and $4 \ii_0 ^1 \sqrt{ V( s^2 ) } \, ds > \sqrt 2 \pi$. 
Then all solutions of (\ref{4.1}) in $ \Er $ are given by Proposition \ref{P4.1} (i). 
Let $ \zeta ( c ) $ be as in Proposition \ref{P4.1}. 
If $\psi _c $ is a solution of (\ref{4.1}), the infimum of $ |\psi |^2 $ is $ \zeta ( c)$. 
We have $ \zeta ( c ) \lra 0 $ as $ c \lra 0 $, and using (\ref{escape}) we see that 
$$
E^1( \psi _c)  \ges 4 \big| H (\sqrt{\zeta(c)} ) \big| \lra  4 \ii_0 ^1 \sqrt{ V( s^2 ) } \, ds > \sqrt 2 \pi \qquad \mbox{ as } c \lra 0 . 
$$
By  Lemma \ref{E1min} (vi) we have $E_{\min}^1 ( p ) < \sqrt{2} p $ for any $ p > 0 $, thus 
$ E_{\min}^1 ( p ) \les  E_{\min}^1 ( \pi ) < \sqrt{2} \pi $ for any $ p \in (0, \pi]$. 

There is $ c_ 0 > 0 $ such that for $ c \in [0, c_0)$ we have $ E^1( \psi _c)  > E_{\min}^1 ( \pi ) \ges  E_{\min}^1 ( p )$ for any $ p \in (0, \pi]$. 
Thus $ \psi _c$ cannot be a minimizer for $E_{\min}^1 (p)$ if $ 0 \les c < c_0$. 

On the other hand,
Theorem \ref{T5.2} implies that there exist minimizers for $E_{\min}^1 (p)$  for all $ p \in (0, \pi]$. By Proposition \ref{P5.3},  the minimizers are necessarily solutions of (\ref{4.1}), but they must have speeds $ c \ges c_0$. 
We infer that $ \ds \inf_{p \les \pi}  \left(E_{\min}^1\right) _{\ell}' ( p ) = \left(E_{\min}^1\right) _{\ell}' ( \pi )\ges c_0$, and  $ E_{\min}^1$ has a cusp at $p = \pi$.

\end{example}

\begin{example} 
\label{Example2}
\rm
Consider $ V \in C^{\infty}(\R, \R)$ having all properties in Example \ref{E4.3}. 
We assume, in addition, that $ \de_1$, $ \de _2$ are small enough, and 

$ \bullet $ $ c_0 \in \left( \frac{\sqrt{2}}{2} , \sqrt{ 2( 1 - \de_1)} \right), $

$ \bullet$ $ V > 0 $ on $[0, 1)$ and $ \ds \int_0 ^1 \sqrt{ V ( s^2 ) } \, ds < \frac{\pi \sqrt{2}}{12}. $

By Lemma \ref{E1min} (v) we have $ E_{\min}^1 ( p ) < \sqrt{2} p $ for any $ p > 0 $, and then Theorem \ref{T5.2} implies that for any $ p \in (0, \pi]$ there exists a minimizer $ \psi \in \Er $ for $ E_{\min}^1 ( p ) $ which satisfies (\ref{4.1}) for some
$c = c( \psi ) \in \left[ (E_{\min}^1)_r ' ( p ) ,  (E_{\min}^1)_{\ell} '(p)\right]$. 

Let $ \mathscr{C} $ be the set of all speeds achieved by minimizers, as in Remark \ref{R5.4}.  We have $ c_0 \not\in \mathscr{C}$ (see Example \ref{E4.3}). The set $ \mathscr{C}$ contains  speeds smaller than $ c_0 $ and speeds  larger than $ c_0$. 
Indeed, we have $ c( \psi ) \lra \sqrt{2} > c_0 $ as $ p \lra 0 $. 
The concavity of  $ E_{\min}^1 $ implies that for any minimizer $ \psi $ of $ E_{\min}^1 ( p ) $ we have 
$ 0 = E_{\min}^1 ( 0 ) \les E_{\min}^1 ( p ) - \left( E_{\min}^1\right)_{\ell}' ( \psi ) p \les E_{\min}^1 ( p ) - c ( \psi ) p $, and this gives 
$ c( \psi ) \les \frac 1p E_{\min}^1 ( p )$. 
In particular,   using Lemma \ref{E1min} (vii) we see  that for any $ p \in \left( \frac{ 2 \pi}{3}, \pi \right)$ we have
$ c ( \psi ) \les \frac 1p  \cdot 4 \int_0^1 \sqrt{V( s^2 )} \, ds < \frac{3}{2 \pi} \cdot 4 \frac{\pi \sqrt{2}}{12} = \frac{\sqrt{2}}{2} < c_0. $

Let $ c_1 = \sup \{ c \in \mathscr{C} \; \mid \; c < c_0 \}$ and $ c_2 = \inf \{ c \in \mathscr{C} \; \mid \; c > c_0 \}$. 
Since $ \mathscr{C} $ is closed in $ [ 0, \sqrt{2})$, we have 
$ c_1, c_2 \in   \mathscr{C} $, $ c_1 < c _0 < c_2 $, and $ ( c_1, c_2 ) \cap  \mathscr{C} = \emptyset$.   
There are $ p_1, p_2 \in (0, \pi]$ and there are minimizers $ \psi_1, \psi _2 $ for 
$ E_{\min}^1( p_1)$ and for $ E_{\min}^1( p_2)$ such that $ c( \psi _1 ) = c_1$  and $ c( \psi _2 ) = c_2$. 
In view of (\ref{Lag-monotone}), we have necessarily $ p_1 \ges p_2$. 
If $ p_1 > p_2$, take any $ p \in ( p_2, p_1)$ and any minimizer $ \psi $ of 
$ E_{\min}^1( p)$. Then we have  $ c( \psi ) \in \mathscr{C}$ and (\ref{Lag-monotone}) gives $ c_1 = c( \psi_1) < c( \psi ) < c( \psi _2) = c_2$, a contradiction. 
Thus $ p_1 = p_2$. 
Then $ E_{\min}^1( p_1)$ admits at least two minimizers $ \psi _1$, $ \psi _2$ having Lagrange multipliers $ c_1 < c_2$, and does not admit any minimizer having Lagrange multiplier $ c \in ( c_1, c_2)$. 
Obviously, the map $ p \longmapsto E_{\min}^1( p)$ is not differentiable at $ p_1$.

Notice that by Proposition \ref{P4.1} (vi), equation (\ref{4.1}) admits solutions in $ \Er $ for almost all $ c \in ( c_1, c_2)$. 
However, these solutions are not global minimizers of the energy at fixed momentum. (We suspect they are local minimizers, at least under some additional assumptions on $V$.)

\end{example}

\begin{example} 
\label{Example3}
\rm
Assume that condition {\bf (A1)} in the introduction is satisfied, that $F$ is locally Lipschitz on $[0, \infty)$, 
$ V > 0 $ on $(1, \infty )$ and there exists $ s \in (0,1)$ such that $ V( s ) = 0 $. 
With the notation at the beginning of section \ref{sect1D}, we have $ g ( s, 0 ) = 4 s V(s)$ and 
$$
\zeta( 0 ) = \sup \{ \zeta \in [0, 1 ) \; \mid \; g ( s, 0 ) = 0 \} = \sup \{ s \in (0, 1 ) \; \mid \; V( s ) = 0 \} \in (0, 1 ) . 
$$
Let $ G( \cdot, 0 ) $ be a primitive of $ \frac{1}{\sqrt{ g ( \cdot, 0 ) }} $ on $ ( \zeta(0), 1 ) $, and let  $ L(0) =\ds  \lim_{s \searrow \zeta ( 0 ) } G( s, 0 )$.  
If $ L(0 ) = - \infty$, equation (\ref{4.1}) with $ c = 0 $ does not admit nontrivial solutions $ \psi $ such that 
$ \ds \lim_{x  \ra \pm \infty }  |\psi (x )  | = 1 $. 
If $L(0)$ is finite, let $ \psi _0 (x ) = \sqrt{\varrho _0 (x)}$, where $\varrho _0$ is as in Proposition \ref{P4.1} (i); 
then all solutions $ \psi $ of (\ref{4.1}) with $ c = 0 $ satisfying  $ \ds \lim_{x  \ra \pm \infty }  |\psi (x )  | = 1 $ are either constant, or are of the form $ e^{ i \al } \psi _0 ( \cdot - x_0)$ for some $ \al , x_0 \in \R$ by Proposition \ref{P4.1} (v), and  a valuation of their  momentum is $0$ (cf. Remark \ref{R4.2}).

Since $ V( 1 ) = V( \zeta(0 )) = 0$,  there exist functions $ \psi \in C^{\infty}( \R)$ such that $ \psi ' \in L^2( \R)$, $ V( |\psi |^2) \in L^1( \R)$ and 
$|\psi ( x ) |^2  \lra \zeta (0)$ as $x  \lra \pm \infty $
(or  $|\psi ( x ) | ^2 \lra  \zeta (0)$ as $x  \lra - \infty $ and  $|\psi ( x ) |^2 \lra 1$ as $x  \lra  \infty $, or vice-versa). 
Since we want to avoid this situation, we need to slightly modify the definition of $ \Er$ as follows: 
$$
\Er  =  \{ \psi \in L_{loc}^1 ( \R) \; \mid \; \psi '\in L^2( \R), \, V( |\psi |^2) \in L^1( \R) \; \mbox{ and } 
\lim_{x  \ra \pm \infty } |\psi ( x ) | = 1  \}. 
$$
The energy $E^1$ does no longer control the Ginzburg-Landau energy $ E_{GL}^1$. 
Indeed, let $ \chi \in C^{\infty } ( \R) $ such that $ \chi = 1 $ on $ ( - \infty, 0 ]$ and $ \chi = \sqrt{ \zeta (0) } $ on $[1, \infty)$. 
Let $ \psi _R( x ) = \chi ( x + R ) $ if $ x \les 0 $, and $ \psi _R ( x ) = \chi ( R - x ) $ if $ x > 0 $. 
Then $ \psi _R \in \Er$, $ E^1 ( \psi _R)$ is independent of $R$  for $ R \ges 1$, and   $ E_{GL}^1 ( \psi _R) \lra \infty $ as $ R \lra \infty$ because  $ \psi _R = \sqrt{ \zeta (0) } $ on $[1 - R, R-1]$.
However, for any $ \eta > 0 $ we may find bounds of $ E_{GL}^1 $ in terms of $ E^1 $ in the set 
$ \{ \psi \in \Er \ \mid \; \inf_{x \in \R} |\psi ( x ) | \geq \zeta(0)  + \eta \}. $

Let $E_{\min}^1 $ be as in (\ref{E1-min}). Then we have $\ds E_{\min}^1 ( p ) \les 4 \int_{\sqrt{\zeta ( 0 ) }} ^1 \sqrt{V ( \tau ^2 ) } \, d \tau $ for any $ p \in \R$. 
To prove this claim, assume for simplicity that $ L(0) $ is finite. Let $ \varrho _0$ be as in Proposition \ref{P4.1} (i). Fix $ p \in \R$.   For $ R> 1$, let 
$$
\rho_R(  x) = \left\{ \begin{array}{ll} \sqrt{ \varrho _0 ( x + R )} & \mbox{if } x < - R, \\ \\
 \sqrt{\zeta( 0 ) } & \mbox{if } x \in [ - R, R],  \\  \\
 \sqrt{ \varrho _0 ( x - R )} & \mbox{if } x >  R,
\end{array}
\qquad
\right. 
\theta_R(  x) = \left\{ \begin{array}{ll} \ds - \frac{p }{2 ( 1 - \zeta ( 0) ) } & \mbox{if } x < - R, \\ \\
\ds \frac{p x}{2R ( 1 - \zeta ( 0) ) } & \mbox{if } x \in [ - R, R],  \\ \\
\ds \frac{p }{2 ( 1 - \zeta ( 0) ) } & \mbox{if } x >  R.
\end{array}
\right. 
$$
Let $ \psi^R  ( x ) = \rho_R ( x ) e^{ i \theta _R (x)}$. 
It is easily seen that $ \psi ^R \in \Er$, 
$$
E^1 ( \psi ^R ) = \int _{\R} |\rho_R '|^2 + \rho_R^2 |\theta_R'|^2 + V ( \rho _R^2) \, dx = 4  \int_{\sqrt{\zeta ( 0 ) }} ^1 \sqrt{V ( \tau ^2 ) } \, d \tau + \frac{ p^2 \zeta ( 0 )}{2 R ( 1 - \zeta(0))^2} , 
$$
and a valuation of the momentum of $ \psi ^R $ is 
$ \ds
\int_{\R} ( 1 - \rho_R^2 ) \theta_R ' \, dx = p .
$
Thus we have $ E_{\min}^1 ( p ) \les E^1 ( \psi^R)$ and letting $ R \lra \infty $ we get the desired conclusion. 
If $ L( 0) = - \infty$, it is easy to see that for any $ \e > 0 $ there exists  a real-valued function $ \rho \in \Er $ such that 
$\ds  E^1 ( \rho ) < 4 \int_{\sqrt{\zeta ( 0 ) }} ^1 \sqrt{V ( \tau ^2 ) } \, d \tau  + \e$ and $ \rho ( 0 ) = \sqrt{\zeta (0) }$. 
We may construct the functions $ \rho_R, \; \theta _R $ and $ \psi ^R $ as above with $ \rho $ instead of $ \sqrt{\varrho_0}$, and we infer that  the same conclusion holds. 
 
Lemma \ref{E1min} still holds with  $ \ds 4 \int_{\sqrt{\zeta ( 0 ) }} ^1 \sqrt{V ( \tau ^2 ) } \, d \tau   $ as upper bound for $E_{\min}^1 (p) $ in part (vii). 
Theorem \ref{T5.2} is still true if we require, in addition, that $ E_{\min}^1 ( p ) < \ds 4 \int_{\sqrt{\zeta ( 0 ) }} ^1 \sqrt{V ( \tau ^2 ) } \, d \tau  $. 
If  $ E_{\min}^1 ( p ) = \ds 4 \int_{\sqrt{\zeta ( 0 ) }} ^1 \sqrt{V ( \tau ^2 ) } \, d \tau  $, the  example here  above suggests that there exist minimizing sequences for $ E_{\min}^1 ( p )$ that do not converge in any reasonable way. 

The set $ \mathscr{C} $ in Remark \ref{R5.4} (i) does not contain $0$. 
Indeed, either $L(0) = - \infty $ and then equation (\ref{4.1}) with $ c = 0 $ does not have solutions in $ \Er$, 
or $ L(0 ) $ is finite and then 
all its solutions have momentum zero. 
If 
$ \ds E_{\min}^1 ( \pi ) < \min \left( \sqrt{2} \pi, 4 \int_{\sqrt{\zeta ( 0 ) }} ^1 \sqrt{V ( \tau ^2 ) } \, d \tau  \right)$, 
there exist minimizers for $  E_{\min}^1 ( \pi ) $, but there exists $ c_0 > 0 $ such that none of these minimizers has speed $ c \in (- c_0, c_0)$. In particular, $(E_{\min}^1)_{\ell}' ( \pi ) \ges c_0$, and then $(E_{\min}^1)_{\ell}' ( p ) \ges c_0$ for any 
$ p \in (0, \pi)$.  For almost every $ c \in (0, c_0)$ we have $ L( c ) > - \infty$ and the functions $ \psi _c $ constructed in Proposition \ref{P4.1} (i) are solutions of (\ref{4.1}), but these solutions cannot be minimizers of $ E_{\min}^1 ( p ) $ for some $ p \in (0, \pi]$.

\end{example}

\section{Minimizing the energy at fixed momentum in $ \Ep $}

Throughout  this section we suppose that the assumptions {\bf (A1), (A2), (B2)} in the introduction hold. We define
$$
E_{\la, \min} ( q ) = \inf \{ E_{\la}( \psi ) \; \big| \; \psi \in \Ep \mbox{ and } \QR (\psi ) = \pr{q} \}.
$$
The next lemma collects the main properties of the function $E_{\la, \min}$. 

\begin{Lemma}
\label{Emin}
Assume that $V$ satisfies the assumptions {\bf (A1)}, {\bf (A2)} and {\bf (B2)} in the introduction.
The function $ E_{\la , \min}$ has the following properties: 

\medskip

i) $ E_{\la, \min} $ is non-negative,  $2 \pi-$periodic, 
$
E_{\la, \min}( p ) = E_{\la, \min} ( | \pr{p} |) 
$
for any $ p \in \R $,  and
$$
E_{\la, \min} ( p ) = \inf \{ E_{\la}( e^{ i \ph } + w ) \; \big| \; \ph \in \dot{H}^1( \R, \R), 
w \in H_{per}^1 ( \R, \C) 
 \mbox{ and } q (\ph , w ) = p \}.
$$

ii) For any $ p \in \R$ and any $ \la > 0 $ we have $ E_{\la,  \min}( p ) \les E_{\min}^1 ( p )$. 
Consequently we have $ E_{\la,  \min}( p )   \les \sqrt{2} |p| $  for all $ p $.

 \medskip

iii) $ E_{\la, \min}$ is sub-additive: for any $ p_1, p_2 \in \R$ there holds 
 $$ E_{\la, \min}( p_1 + p_2 ) \les E_{\la,  \min}( p_1 ) + E_{\la, \min} ( p_2).$$
 
iv) $ E_{\la, \min} $  is $ \sqrt{2}-$Lipschitz on $ \R$. 
\medskip

v) For any $ \de > 0 $ there exists $ p_{\de } > 0 $ such that $ E_{\la, \min}( p ) \ges (1 - \de) \sqrt{2} p $
 for any $ p \in (0, p_{\de})$. 

\medskip

vi) For any fixed $p $, the mapping $ \la \longmapsto E_{\la, \min} ( p ) $ is non-decreasing.
Assume that  $ p \in (0,  2\pi ) $ and there exist $ p_1, p_2 \in (0, 2 \pi)$ satisfying $ p  = \frac{ p_1 + p_2 }{2 } $ and 
$ E_{\min}^1 ( p ) > \frac 12 \left(  E_{\min}^1 ( p_1 )  + E_{\min}^1 ( p _2) \right). $
Then there exists $ \la _ * ( p ) > 0 $ such that 
 $ E_{\la, \min} ( p ) < E_{\min}^1 ( p ) $ for any $ \la \in (0, \la _*( p))$. 

\medskip

vii) Assume that $ p_0 \in (0, 2 \pi) $ satisfies 
$$\ds \liminf_{ h \ra 0 } \frac{E_{\min}^1 ( p_0 + h )  + E_{\min}^1( p_0 - h ) - 2 E_{\min}^1 ( p_0 ) }{ h^2} > - \infty . $$ 
Then there exists $ \la ^* ( p_0 ) > 0 $ such that $ E_{\la, \min} ( p _0 ) = E_{\min}^1 ( p _0 ) $ for any $ \la \ges \la ^*( p_0 )$.

\medskip

viii) $ E_{\min}^1 $  is concave on $[0, 2 \pi]$. 

\medskip

ix) Let $ p \in (0, \pi]$ and let $ \la _ s ( p ) = \ds \sup \{ \la > 0 \; \mid \; E_{\la, \min}( p ) < E_{\min}^1 (p)\} $.  
The mapping $  \la \longmapsto E_{\la, \min} ( p ) $ is strictly increasing on $ (0, \la_s ( p ) )$ and $ E_{\la, \min} ( p ) = E_{\min}^1 (p)$ if $ \la \ges \la_s ( p )$. 

\end{Lemma}

\begin{remark} 
\rm
i)  If the assumption in Lemma \ref{Emin} (vi) does not hold for some $ p \in (0, \pi]$, the concavity of $E_{\min}^1 $ implies that 
$  E_{\min}^1 ( p ) = \frac 12 \left(  E_{\min}^1 ( p - \de ) +  E_{\min}^1 ( p + \de) \right) $ for any $ \de \in (0, p )$, and then we infer that $ E_{\min}^1$ must be affine on $ [0, 2 p ]$. 
This is impossible if $ p \in ( \frac{\pi}{2}, \pi ] $ because $ E_{\min}^1 ( 0 ) = 0 $ and 
$ E_{\min}^1$ achieves its positive maximum at $ \pi$. 
If $ V > 0 $ on  $[0, 1 )$ and $ V( s ) \les \frac 12 ( 1 - s )^2 + \frac 38 (1 - s )^3 $ on some interval $ ( 1 - \eta, 1]$, 
Lemma \ref{E1min} (ii) and (vi) implies that $E_{\min}^1$ cannot be linear on $[0, p ]$ for any $ p \in (0, \pi]$. 
Remember that  $E_{\la, \min} ( p ) =E_{\la, \min} ( 2 \pi - p )$.
Therefore  the conclusion of Lemma \ref{Emin} (vi) holds without any additional assumption if $ p \in ( \frac{\pi}{2}, \frac{ 3 \pi}{2} ) $, and it holds for any $p\in (0, 2 \pi) $ under the assumption of Lemma \ref{E1min}  (vi).

\medskip

ii) It is well-known that a concave function is twice differentiable almost everywhere. 
The limit in Lemma \ref{Emin} (vii) exists and is equal to $(E_{\min}^1 ) '' ( p )$ for almost every $ p \in (0, 2 \pi)$, and for any such $p$ the conclusion of Lemma \ref{Emin} (vii) holds true. Then parts (vi) and (vii) of the above lemma show that 
there is some critical value $ \la^*(p)$ such that $ E_{\la , \min} (p)< E_{\min}^1 (p)$ for $ \la < \la^*(p)$, and 
$ E_{\la , \min} (p) = E_{\min}^1 (p)$ for $ \la \ges \la^*(p)$.
The proof shows that one can give upper bounds for the critical value $ \la^*(p)$ 
if a lower bound for $(E_{\min}^1 ) '' ( p )$ is known (see (\ref{LAM}) below).

Using an argument in the proof of part (viii) one can show that 
$ \la \longmapsto E_{\la, \min} ( p)$ is strictly increasing on $(0, \la^* (p)]$. 

The results in Lemma \ref{Emin} (vi) and (vii) are not surprising. 
Recall that $E_{\la }$ is a rescaled energy that  comes from minimizing the energy $E$ on  $ \R \times \T_{\frac{1}{\la}}$ where $ \T_{\tau}$ is the $1-$dimensional torus of length $\tau$.
When $ \la $ is large the torus $\T_{\frac{1}{\la}}$ is too narrow  and variations with respect to the variable $y$ would be energetically too costly. 
On the contrary, on large  tori one can find better competitors than the $1-$dimensional minimizers of $E_{\min}^1$. 

In the case of the Gross-Pitaevskii nonlinearity $F(s ) = 1- s$, the function  $E_{\min}^1 $ is known explicitly 
(see \cite{BGS-survey} or Example \ref{GP}) and it turns out that it is $C^2 $ on $(0, 2 \pi)$, hence in this particular case 
the conclusion of Lemma \ref{Emin} (vii) holds for any value of $p$.

\end{remark}

{\it Proof. } 
The proof of (i) is similar to the proof of Lemma \ref{E1min} (i). 

\medskip

(ii) Consider any $ \psi \in \Er $ satisfying $ \PR ( \psi ) = p $. Let $ \psi^{\sharp} ( x, y ) = \psi ( x ) $. 
It is obvious that $ \psi^{\sharp}  \in \Ep $ and 
$ E_{\la } ( \psi ^{\sharp} ) = E^1 ( \psi )$ for any $ \la > 0 $. 
If $ \ph \in \dot{H}^1( \R)$ and $ w \in H^1 ( \R)$ are such that $ \psi = e^{i \ph } + w$ in $ \R$, we have 
$ \psi^{\sharp } = e^{ i \ph ( x)} +  w^{\sharp}$, where $ w^{\sharp} ( x, y ) = w( x)$. 
It is obvious that $ q(\ph, w^{\sharp} ) = p ( \ph , w)$, thus $ \QR( \psi^{\sharp} ) = \PR( \psi ) = \pr{p}$. 
We infer that $E_{\la, \min} ( p ) \les E_{\la } ( \psi ^{\sharp}) = E^1 ( \psi)$. 
Since the last inequality holds for all  $ \psi \in \Er $ such that $ \PR ( \psi ) = p $, the conclusion follows. 

\medskip

The proofs of (iii) and (iv) are very similar to the proofs of Lemma \ref{E1min} (iv) and (v), respectively, and we omit them. 
For the proof of part (v) we need some results from \cite{CM} and \cite{M10}. 
These results are stated in Lemma \ref{LemmaVeche} below. We postpone the proof of part (v) after the proof of Lemma \ref{L6.5}.

\medskip

vi)  If $ 0 < \la _1 < \la _ 2 $ it is obvious that $ E_{\la_1} ( \psi ) \les E_{\la _ 2} ( \psi ) $ for any $ \psi \in \Ep$, and this trivially implies that $ E_{\la_1, \min}( p ) \les E_{\la _ 2, \min} ( p )$ for all $p$. 

It suffices to consider the case $ p \in (0, \pi]$. Fix $ p_1, p_2 \in (0, 2 \pi)$ and $ \e > 0 $ such that  $ p = \frac{ p_1  + p_2 }{2} $ and 
$$
10 \e < E_{\min}^1 ( p ) - \frac 12 \left( E_{\min}^1 ( p _1) + E_{\min}^1 ( p _2) \right) . 
$$

By Lemma \ref{approx} (i), there exist functions $ \psi _ j = e^{ i \ph _j } + w _ j \in \Er$ 
such that $ \ph _j \in \dot{H}^1 \cap C^{\infty} ( \R, \R )$, $ w_j \in C_c^{\infty } ( \R, \C)$ and $ \tilde{A} > 0 $ such that 
$ \mbox{supp} ( w_j ) \subset [-\tilde{A}, \tilde{A} ]$, $ \ph _j $ are constant on $ ( - \infty, - \tilde{A} ]$ and on $ [\tilde{A} , \infty ) $ and 
$$
p ( \ph _ j , w_j ) = p_j \qquad \mbox{ and } \qquad E^1 ( \psi _ j ) < E_{\min}^1 ( p_j ) + \e 
\qquad \mbox{ for } j = 1, 2. 
$$
If $ \ph _1 $ and  $ \ph _2$  take different values near $ \pm \infty$, we modify $ \ph _j $ on 
$ ( - \infty, - \tilde{A} ]$ and on $ [\tilde{A}, \infty ) $ 
in such a way that $ \ph _1 = \ph _2 = constant $ on $ ( - \infty, - A ] \cup [ A, \infty ) $ (where $ A$ may eventually be much larger than $\tilde{A}$) 
and $ \ds \ii_{- A }^{- \tilde{A}} |\ph _j ' | ^2 \, dx +  \ii_{\tilde{A}} ^{A} |\ph _j ' | ^2 \, dx < \e $ for $ j = 1, 2$. 
We still denote $ \ph _j $ the modified functions. 
After this modification we have $ p ( \ph _ j , w_j ) = p_j  $ and $  E^1 ( \psi _ j ) < E_{\min}^1 ( p_j ) +2 \e $ for $ j = 1, 2$. 

Let $ 0 < \eta < \frac 18$ (the value of $\eta  $ will be chosen later). 
Take $ \chi \in C^{\infty }( \R) $ such that $\chi$ is $  1- $periodic, $  \quad 0 \les \chi \les 1$, 
$$   
\chi = 1 \mbox{ on } \left[ 0 , \frac 14 - \eta \right] \cup  \left[ \frac 34 + \eta, 1 \right], \quad
\chi = 0 \mbox{ on } \left[ \frac 14 + \eta, \frac 34 - \eta \right] , \quad 
\mbox{ and } \ii_0 ^1 \chi ( y ) \, dy = \frac 12. 
$$
Let 
$$
\begin{array}{c}
\ph ( x ) = \frac 12 ( \ph _1 ( x ) + \ph _2 ( x ) ),  
\\
\\
w( x, y  )= \chi ( y ) \left[ e ^{ i \ph _1 ( x ) } + w_1 ( x ) \right] +  ( 1 - \chi ( y ) )  \left[ e ^{ i \ph _2 ( x ) } + w_2 ( x ) \right] - e^{ i \ph ( x ) }, \quad \mbox{ and } 
\\
\\
\psi (x, y ) = e^{ i \ph ( x ) } + w ( x, y ) =  \chi ( y ) \psi _1 ( x ) +  ( 1 - \chi ( y ) ) \psi _2 ( x ) .
\end{array}
$$
Obviously, $ \ph \in \dot{H}^1 \cap C^{ \infty} ( \R, \R)  $, 
$ w \in C^{\infty} (\R ^2 ) $ and $ w $ is $1-$periodic with respect to the second variable, and $ w = 0 $ on 
$\left(   ( - \infty, - A ] \cup [ A, \infty ) \right) \times \R$. 
Let $ p( \ph, w)$ be as in (\ref{p}) and let  $d[\ph, w]$ be as at the beginning of section \ref{defmomp}. 
Since $ \ph = \ph _1 = \ph _ 2 = constant $ on  $ ( - \infty, - A ] \cup [ A, \infty ) $, 
using (\ref{quantif}) we get 
$$
p ( \ph _j , w_j ) = p ( \ph, w_j + e^{ i \ph _j } - e^{ i \ph } ) \qquad \mbox{ for } j = 1,2.
$$
Using this simple observation, after a straightforward computation we obtain 
$$
\begin{array}{l}
\ds \ii_{\R} d [ \ph , w ] ( x, y ) \, dx = \chi ( y ) p ( \ph _1 , w_1 ) + ( 1 -  \chi ( y ))  p ( \ph _2 , w_2) 
\\
\\
\qquad \qquad 
\ds - \chi ( y ) ( 1 - \chi ( y ) ) \ii_{\R} \langle  i \left( e^{ i \ph _1} + w_1 \right) '  - i \left( e^{ i \ph _2} + w_2 \right) ' , 
e^{ i \ph _1} + w_1 - ( e^{ i \ph _2} + w_2 ) \rangle \, dx . 
\end{array}
$$
 Integrating on $[0, 1 ] $ we find
$$
q ( \ph, w ) = \frac 12 p ( \ph _1, w_1 ) + \frac 12 p ( \ph _2, w_2 ) - \ii_0^1 \chi ( y ) ( 1 - \chi ( y ) )  \, dy 
\cdot \ii_{\R} \langle i \psi _1 '- i \psi _2 ' , \psi _1 - \psi _2 \rangle \, dx . 
$$
Notice that $ 0 \les  \ii_0^1 \chi ( y ) ( 1 - \chi ( y ) )  \, dy  < \eta $ because $ 0 \les  \chi ( y ) ( 1 - \chi ( y ) )\les \frac 14$ and 
$ \chi ( y ) ( 1 - \chi ( y ) ) = 0 $ on 
$ [0, \frac 14 - \eta ] \cup [ \frac 14 + \eta, \frac 34 - \eta ] \cup [ \frac 34 + \eta, 1 ] $. 
Denoting 
$ K = \ds \Big|  \ii_{\R} \langle i \psi _1 '- i \psi _2 ' , \psi _1 - \psi _2 \rangle \, dx \Big| $, we have shown that 
$$
\big| q ( \ph, w ) - \frac{ p_1 + p_2}{2} \big| < \eta K. 
$$
We know that for any $ \la > 0$, the function $ E_{\la, \min} $ is $ \sqrt{2} -$Lipschitz. 
If $ \eta K < \frac{\e}{\sqrt{2}} $ we have 
$$
\Big| E_{\la, \min } (  q ( \ph, w ) )-  E_{\la, \min } \left( \frac{ p_1 + p_2}{2} \right) \Big|  < \e 
\qquad \mbox{ for any } \la  > 0 . 
$$

Let $ M = \sup_{t \in [0, 1 ]} E^1 \left( t \psi _ 1 + ( 1 - t ) \psi _2  - e^{ i \ph } \right). $
It is straightforward to see that $ M $ is finite. 
Since 
$ \psi ( \cdot, y ) = \psi _1 $ for $ y \in [0, \frac 14 - \eta ] \cup [ \frac 34 + \eta, 1 ]$ and 
$ \psi ( \cdot, y ) = \psi _2 $ for $ y \in  [ \frac 14 + \eta , \frac 34 - \eta ] $, 
we infer that 
$$
\ii_{\R \times [0, 1 ] } \Big| \frac { \p \psi }{\p x } \Big| ^2  +  V ( |\psi |^2 ) \, dx \, dy
\les \left( \frac 12 - 2 \eta \right) E^1 ( \psi _1 ) + \left( \frac 12 - 2 \eta \right) E^1 ( \psi _2 )  + 4 \eta M . 
$$
Now choose $ \eta $ such that $ 0 < \eta < \min \left( \frac{ \e}{\sqrt{2} K}, \frac{\e}{4 M} \right).$
Then for any $ \la > 0 $ we have 
$$
\begin{array}{l}
E_{\la, \min} ( p ) = E_{\la, \min } \left( \frac{ p_1 + p_2}{2} \right) \les E_{\la , \min} ( q ( \ph, w )) + \e 
\les E_{\la } (\psi ) + \e 
\\
\\
\ds = \ii_{\R \times [0, 1 ] } \Big| \frac { \p \psi }{\p x } \Big| ^2  +  V ( |\psi |^2 ) \, dx \, dy + \la ^2 \Big\| \frac{ \p \psi }{\p y } \Big\|_{L^2 ( \R \times [0,1])} ^2  + \e
\\
\\
\ds \les \frac 12 E^1  ( \psi _1 ) +  \frac 12 E^1  ( \psi _2 )  + 2 \e + \la ^2 \Big\| \frac{ \p \psi }{\p y } \Big\|_{L^2 ( \R \times [0,1])} ^2
\\
\\
\ds \les \frac 12 E_{\min}^1 ( p_1 ) + \frac 12 E_{\min}^1 ( p_2 )  + 4 \e + \la ^2 \Big\| \frac{ \p \psi }{\p y } \Big\|_{L^2 ( \R \times [0,1])} ^2.
\end{array}
$$
For $ \la $ sufficiently small, so that $ \la ^2 \Big\| \frac{ \p \psi }{\p y } \Big\|_{L^2 ( \R \times [0,1])} ^2 < \e$, we get 
$$
E_{\la, \min} ( p ) < \frac 12 E_{\min}^1 ( p_1 ) + \frac 12 E_{\min}^1 ( p_2 )  + 5 \e < E_{\min}^1 ( p ), 
$$
as desired.

\medskip

vii) 
Choose $ \de > 0 $ sufficiently small and $ L > - \infty $  such that $ ( p_0 - \de , p_0 + \de ) \subset ( 0 , 2 \pi)$ and 
\beq
\label{6.1}
\frac{E_{\min}^1 ( p_0 + h )  + E_{\min}^1( p_0 - h ) - 2 E_{\min}^1 ( p_0 ) }{ h^2} \ges L  > - \infty 
\qquad \mbox{ for all } h \in (0, \de ].
\eeq

Let $ M =  E_{\min}^1  ( p_0 ) + 1 $. 
By Lemma \ref{approx} (ii) and the discussion preceding  Definition \ref{mom}, there exist mappings $ \psi = e^{ i \ph } + w \in \Ep$ such that $ \ph \in \dot{H}^1 ( \R ) \cap C^{\infty}( \R)$, 
$ w$ is the $1-$periodic extension with respect to the second variable of a function in $ C_{c}^{\infty} (\R \times ( 0, 1 ) )$, 
$ q ( \ph, w ) = p_0 $ and $ E_{\la } ( \psi )$ is arbitrarily close to $E_{\la, \min} ( p_0)$, in particular $E_{\la, \min} ( \psi) < M$.
For any such $ \psi$, the  upper bound 
$$
E_{\la }( \psi ) = \ii_{\R \times [0,1]  } \bigg| \frac{ \p \psi }{\p x } \bigg|^2  + \la ^2 \bigg| \frac{ \p \psi }{\p y } \bigg|^2 + V( |\psi |^2) \, dx \, dy \les M 
$$
gives $ \Big\| \frac{ \p \psi }{\p x } \Big\|_{L^2 ( \R \times [0,1] )} \les \sqrt{M}$ and  
$ \Big\| \frac{ \p \psi }{\p y } \Big\|_{L^2 ( \R \times [0,1] )} \les \frac{\sqrt{M}}{\la}$. 
By Lemma \ref{continuous}, for any $y_1, y_2 \in [0, 1] $ we have 
\beq
\label{strans}
|  p ( \ph , w ( \cdot, y_2)) - p ( \ph , w ( \cdot, y_1))  | \les 2 \Big\| \frac{ \p \psi }{\p x } \Big\|_{L^2((\R \times [y_1, y_2])}
\Big\| \frac{ \p \psi }{\p y } \Big\|_{L^2((\R \times [y_1, y_2])} \les \frac{2M}{\la}.
\eeq
The mapping $ y \longmapsto p ( \ph, w ( \cdot, y ))$ is continuous by (\ref{strans})  
and Lebesgue's dominated convergence theorem.  
Using Remark \ref{basic} (i), a valuation of the momentum of $ \psi $ is 
\beq
\label{bun}
p_0 = q ( \ph, w ) = \int_0 ^1 p ( \ph, w( \cdot, y ) ) \, dy  . 
\eeq
By the continuity of $ y \longmapsto p ( \ph, w ( \cdot, y ))$, there exists $ y_0 \in [0,1]$ such that $ p ( \ph, w ( \cdot, y_0 )) = p_0$.

Choose $ \la^* ( p_0 ) $ sufficiently large, such that 
\beq
\label{LAM}
 \frac{2M}{\la^*(p_0)} < \frac{ \de}{2} \qquad \mbox{ and } \qquad \frac L2 + \frac{\left(\la^* (p_0)\right)^2}{4M } > 1. 
 \eeq
From now on we will assume that $ \la > \la^*(p_0 ) $. 
For $ \psi $ as above, denote $ p( y ) =  p ( \ph , w ( \cdot, y)) $ and let $ \de _{\psi } = \sup_{y \in [0,1] } | p( y ) - p_0 |$. 
By (\ref{strans}) and (\ref{LAM}) we have $ \de_{\psi } \les \frac{\de}{2} $ and $ p(y ) \in [p_0  -\de_{\psi}, p_0 + \de_{\psi } ] \subset ( p_0 - \de , p_0 + \de ) \subset ( 0 , 2 \pi)$ for all $ y$.

Obviously, for any $ y \in [0, 1 ] $ we have $ \psi ( \cdot, y ) = e^{ i \ph } + w( \cdot, y )  \in \Er $ and consequently
\beq
\label{partial1}
E^1( \psi ( \cdot, y ) ) \ges E_{\min}^1 ( p ( \ph, w ( \cdot, y )))  = E_{\min}^1 ( p (  y)). 
\eeq
If $ f : [a, b ] \lra \R$ is concave and continuous, for any $ t \in [a, b ] $ we have 
$$ f( t ) \ges f( a ) + \frac{ t - a }{ b - a } ( f( b ) - f ( a)) . 
$$
Since $E_{\min}^1$ is concave on $(0, 2 \pi )$ and $ p(y ) \in [p_0 - \de_{\psi}, p_0 + \de _{\psi}]$,  we  get 
\beq
\label{lower}
E_{\min}^1( p(y)  ) \ges E_{\min}^1 ( p_0 - \de_{\psi } ) + \frac{ p(y) - ( p_0 - \de_{\psi})}{2 \de_{\psi } } \left[ E_{\min}^1 ( p_0 + \de_{\psi} ) -  E_{\min}^1 ( p_0 - \de_{\psi} ) \right].
\eeq
Using Fubini's Theorem, (\ref{partial1}), (\ref{lower}) and the fact that $ \ds \ii_0^1 p(y ) \, dy = p_0 $ (see (\ref{bun}))  we obtain 
\beq
\label{lower2}
\begin{array}{l}
\ds \ii_{\R \times [0, 1 ] } \bigg| \frac{ \p \psi}{\p x } \bigg|^2 + V( |\psi |^2 ) \, dx 
 =  \ds \ii_0 ^1 E^1 ( \psi ( \cdot , y ) ) \, dy 
\ges \ii_0^1 E_{\min}^1 ( p( y ) ) \, dy 
\\
\\
  \ges   \ds \frac 12 \left[ E_{\min}^1 ( p_0 + \de_{\psi} ) +  E_{\min}^1 ( p_0 - \de_{\psi} ) \right].
 \end{array}
\eeq
From (\ref{strans}) we get 
\beq
\label{lower3}
\Big\| \frac{ \p \psi }{\p y } \Big\|_{L^2((\R \times [0,1])}\ges \frac{ \de_{\psi}}{2 \Big\| \frac{ \p \psi }{\p x } \Big\|_{L^2((\R \times [0,1])} } \ges \frac{ \de_{\psi}}{ 2 \sqrt{M}}. 
\eeq 
From (\ref{lower2}), (\ref{lower3}),  (\ref{6.1}) and (\ref{LAM}) we obtain 
$$
\begin{array}{l}
\ds E_{\la}( \psi ) - E_{\min}^1 ( p_0 ) \ges  \ds \frac 12 \left[ E_{\min}^1 ( p_0 + \de_{\psi} ) +  E_{\min}^1 ( p_0 - \de_{\psi} ) \right] - 
 E_{\min}^1 ( p_0 )
 + \frac{ \la^2  \de_{\psi}^2 }{ 4 M } 
 \\
 \\
 = \ds \de_{\psi}^2 \left( \frac 12 \frac{ E_{\min}^1 ( p_0 + \de_{\psi} ) +  E_{\min}^1 ( p_0 - \de_{\psi} )  - 2 E_{\min}^1 ( p_0 )}{\de_{\psi }^2 } + \frac{\la^2 }{4M} \right)  >\de_{\psi}^2  .
 \end{array}
 $$
Since the last estimate holds for any $ \psi$ as considered above, the conclusion follows.

\medskip

(viii) We   prove the concavity of $E_{\la, \min}$ in several steps. 

\medskip

{\it Step 0. Functional setting. } 
We consider the space 
$$
\begin{array}{rcl}
\Xo & = & \{ w :\R^2 \lra \C \; \mid \; w \mbox{ is } 1-\mbox{periodic with respect to the second variable, }
\\ & &  w \mbox{ is piecewise } C^2 \mbox{ and there exists } a > 0 \mbox{ such that supp}(w) \subset [ - a, a ] \times \R \}.
\end{array}
$$
By "piecewise $C^2$" we mean that $ w$ is continuous and there exist finitely many points $ 0 = y_0 < y_1 < \dots < y_n = 1 $ such that 
for each $ j \in \{ 1, \dots, n \}$ there exists a mapping $ \tilde{w}^j $ that  is $C^2 $ on some larger strip 
$ \R \times (y_{j - 1} - \de, y _j + \de )$, and  
$w_{|  \R \times [ y_{j- 1}, y _j] } = \tilde{w}_{|  \R \times [ y_{j - 1}, y _j] } ^j$.
We consider the space $ \Xo $ for the following reasons: we need a function space $ \Xo \subset H_{per}^1 $ 
such that any function $ \psi \in \Ep $ can be approximated by functions 
of the form $ e^{ i \ph ( x ) } + w( x, y )$, where $ \ph \in \dot{H}^1( \R)$ and $ w \in \Xo$ 
(this can be done in view of Lemma \ref{approx} (ii)),
we need to use Lemma \ref{continuous} (which obviously extends to functions in $ \Xo$), and we need $\Xo $ to be stable 
under  a reflection procedure that we will describe below.

Given any $ w \in \Xo $ and $ y _ 0 \in \R$, we define the functions $ T_{1, y_0} w $ and $ T_{2, y_0} w $ 
on $ \R \times [ 0 , 1 ] $ as follows, then we extend them to $ \R^2$ as $1-$periodic functions with respect to the second variable:
$$
\begin{array}{c}
T_{1, y_0} w ( x, y ) = \left\{ \begin{array}{ll} 
w(x, y_0 + y ) & \mbox{ if } y \in [0, \frac 12], \\
w(x, y_0 +1 -  y ) & \mbox{ if } y \in [ \frac 12, 1],
\end{array}
\right.
\\
\\
T_{2, y_0} w ( x, y ) = \left\{ \begin{array}{ll} 
w(x, y_0 + 1 - y ) & \mbox{ if } y \in [0, \frac 12], \\
w(x, y_0 + y ) & \mbox{ if } y \in [ \frac 12, 1].
\end{array}
\right.
\end{array}
$$
It is obvious that $ T_{1, y_0} w $ and $ T_{2, y_0} w $ belong to $ \Xo $ for any $ w\in \Xo $ and any $ y_0 \in \R$, and we have 
 $  T_{j, y_0} w ( x, y ) = T_{j, y_0} w ( x, 1 - y ) $ for all $ y \in [0, 1]$ and $ j = 1, 2$. 

For $ \ph \in \dot{H}^1( \R)$ and $ w \in \Xo $, let $ \psi ( x , y ) = e^{ i \ph( x )} + w ( x, y )$ and let 
$ d [\ph, w ]$ and $ q ( \ph, w )$ be as in subsection \ref{defmomp}.
It is easily seen that 
\beq
\label{lopes}
\begin{array}{l}
\ds q ( \ph , T_{1, y_0}  w ) = \ii_0^1 \ii_{\R} d  [ \ph, T_{1, y_0}  w ]( x, y ) \, dx \, dy 
= 2 \ii_{y_0}^{ y_0 + \frac 12}  d  [ \ph,   w ]( x, y ) \, dx \, dy , 
\\
\\
\ds q ( \ph , T_{2, y_0}  w ) = \ii_0^1 \ii_{\R} d  [ \ph, T_{2, y_0}  w ]( x, y ) \, dx \, dy 
= 2 \ii_{y_0 - \frac 12}^{ y_0 }  d  [ \ph,   w ]( x, y ) \, dx \, dy , 
\\
\\
\ds E_{\la } ( e^{ i \ph } + T_{1, y_0}  w )  
= 2 \ii_{y_0}^{ y_0 + \frac 12}  \Big| \frac{ \p \psi }{\p x } \Big|^2 + \la ^2 \Big| \frac{ \p \psi }{\p y } \Big|^2 + V( |\psi |^2 )  \, dx \, dy , 
\\
\\
\ds E_{\la } ( e^{ i \ph } + T_{2, y_0}  w ) 
= 2 \ii_{y_0- \frac 12 }^{ y_0 }  \Big| \frac{ \p \psi }{\p x } \Big|^2 + \la ^2 \Big| \frac{ \p \psi }{\p y } \Big|^2 + V( |\psi |^2 )  \, dx \, dy ,
\qquad \mbox{ and } 
\\
\\
E_{\la } ( e^{ i \ph } + T_{1, y_0}  w )  + E_{\la } ( e^{ i \ph } + T_{2, y_0}  w )  = 2 E_{\la } ( \psi ) . 
\end{array}
\eeq
For any $ \ph $ and $ w $ as above, the function 
$$
\Upsilon_{\ph, w } ( t ) =  \ii_{t}^{ t  + \frac 12}  d  [ \ph,  w ]( x, y ) \, dx \, dy - 
 \ii_{ t - \frac 12}^{t}  d  [ \ph,  w ]( x, y ) \, dx \, dy
$$
is continuous, $ 1-$periodic on $ \R$ and $ \Upsilon_{\ph, w } ( t + \frac 12 ) = - \Upsilon_{\ph, w } ( t )$ for any $ t $.
We denote 
$$
\omega ( \ph, w ) = \sup_{t \in [0,1]} \mid \Upsilon_{\ph, w } ( t ) \mid . 
$$ 
For any $ q \in (0, 2 \pi)$ and any $ \la > 0 $, we denote
$$
\de_{\la }( q ) = \inf_{\e > 0 } \left( \sup \left\{  \omega ( \ph, w ) \; \mid \; 
\ph \in \dot{H}^1( \R), w \in \Xo, q ( \ph , w ) = q \mbox{ and } E_{\la } ( e^{ i \ph } + w ) < E_{\la, \min}( q ) + \e \right\} \right) . 
$$

{\it Step 1.  Assume that $ \de_{\la} ( q ) > 0 $. 
Then for any $ \eta \in (0, \de_{\la}(q)) $ we have }
\beq
\label{round1}
E_{\la, \min}( q - \eta ) + E_{\la, \min}( q + \eta )  \les 2 E_{\la, \min}( q) . 
\eeq
To see this, fix $ \eta \in (0, \de_{\la}(q)) $. Fix $ \e > 0 $. 
By the definition of $  \de_{\la }( q ) $, there exist $ \ph \in \dot{H}^1( \R)$ and   $w \in \Xo$ such that  $ q ( \ph , w ) = q $, 
$ E ( e^{ i \ph } + w ) < E_{\la, \min}( q ) + \e $ and $  \omega ( \ph, w ) > \eta$. 
Since $ \Upsilon_{\ph, w}$ is continuous, $1-$periodic and $ \Upsilon_{\ph, w} ( t + \frac 12 ) = -  \Upsilon_{\ph, w}(t)$ for any $t$,
there exists $ t_0 \in [0, 1 ] $ such that $ \Upsilon_{\ph, w} ( t_0 ) = \eta$. 
Let $ w_1 = T_{1, t_0} w $ and $ w_2 = T_{2, t_0} w$. 
From the first two equalities in (\ref{lopes}) we get $q(\ph, w_1 ) = q + \eta $ and $ q( \ph, w_2 ) = q - \eta$, 
and we infer that $ E_{\la } ( e^{ i \ph } + w_1 ) \ges E_{\la, \min} ( q + \eta )$ and $ E_{\la } ( e^{ i \ph } + w_2 ) \ges E_{\la, \min} ( q - \eta )$.
Then using the last equality in (\ref{lopes}) we find
$$
 E_{\la, \min}( q - \eta ) + E_{\la, \min}( q + \eta ) \les E_{\la } ( e^{ i \ph } + w_2 ) + E_{\la } ( e^{ i \ph } + w_1 )
= 2E_{\la } ( e^{ i \ph } + w ) < 2 E_{\la, \min}( q ) + 2\e. 
$$
Since the above inequality  holds for any $ \e > 0 $, (\ref{round1}) is proven. 

\medskip

{\it Step 2. If $ q \in (0, 2 \pi) $ is such that $ \de_{\la}(q ) = 0 ,$  then $ E_{\la, \min}( q ) = E_{\min}^1 ( q)$. }

Let  $ M = E_{\la, \min}( q ) + 1$.  Choose  $ r > 0 $ such that $ ( q - r, q + r ) \subset (0, 2 \pi)$. 
Then choose $ \e_0 \in (0, 1 ) $ such that $ \e _ 0 < r $ and $   \frac{ 4 \sqrt{M \e _0} }{\la }  < r $. 

Fix $ \e \in (0, \e_0) $. 
By the definition of $\de_{\la}(q)$, there exist $ \eta_{\e} > 0 $ such that for any 
$\ph \in \dot{H}^1( \R) $ and any $ w \in \Xo $ satisfying $ q ( \ph , w ) = q $ and $ E_{\la } ( e^{ i \ph } + w ) < E_{\la, \min}( q ) + 2 \eta_{\e } $ we have $\omega( \ph, w ) < \e$. 
We may assume that $ \eta_{\e } \les \e$. 

Consider  $ \ph \in \dot{H}^1( \R)$ and   $w \in \Xo$ such that  $ q ( \ph , w ) = q $ and 
$ E ( e^{ i \ph } + w ) < E_{\la, \min}( q ) + \eta_{\e } $. 
Proceeding as in step 1, we see that  there exists $ y_0 \in [0,1]$ such that  $ \Upsilon_{\ph, w} ( y_0 ) = 0.$
Denote  $ w_1 = T_{1, y_0} w $, $ w_2 = T_{2, y_0} w$,  $ \psi = e^{ i \ph } + w $,  and $ \psi_j ( x, y ) = e^{ i \ph (x ) } + w_j ( x, y)$ for $j=1,2$. 
By (\ref{lopes}) we have
$$
q( \ph, w_1 ) = q( \ph, w_2 ) = q( \ph, w) = q \quad \mbox{ and } \quad 
E_{\la} ( \psi _1 ) + E_{\la}( \psi _2 ) = 2 E_{\la } (\psi) < 2E_{\la, \min}( q ) + 2 \eta_{\e} 
$$
and we infer that 
$$
E_{\la, \min}( q )  \les E_{\la} ( \psi _j  ) < E_{\la, \min}( q )  + 2\eta_{\e}  \qquad \mbox{ for } j = 1,2. 
$$
Taking into account how $ \eta_{\e}$ was chosen,  we infer that $\omega( \ph, w _j) < \e$ for $ j = 1,2$. 

Let $ w_{j, 1 } = T_{1, \frac 14} w_j $ and $ w_{j, 2 } = T_{2, \frac 14} w_j $ for $ j = 1,2 $, then let 
$ \psi _{j, \ell } (x , y ) = e^{ i \ph(x)} + w_{j, \ell}( x , y )$ for $ j, \ell \in \{ 1, 2 \}$. 
Using the first equality in (\ref{lopes}) we find 
$$
\big| q ( \ph, w_{1,1} ) - q ( \ph, w_1 ) \big| = \big| \Upsilon_{\ph, w_1} \left( \frac 14 \right) \big| \les \omega( \ph, w _1) < \e. 
$$
Similarly we get $ \big| q ( \ph, w_{j,\ell} ) - q ( \ph, w_j ) \big| < \e $ for $ j, \ell \in \{ 1, 2 \}$, and this gives 
$  q ( \ph, w_{j,\ell} )  \in ( q - \e , q + \e )$. 
Since $ E_{\la, \min}$ is $ \sqrt{2}-$Lipschitz, we get
\beq
\label{round2}
E_{\la, \min} ( q ( \ph, w_{j, \ell})) \ges E_{\la, \min} (q ) - \sqrt{2 } \e
\qquad \mbox{ for }  j, \ell \in \{ 1, 2 \}.
\eeq

Now we observe that  by construction,  the functions $ w_{ j, \ell}$ are $\frac 12-$periodic with respect to the second variable $y$. 
Let $ \tilde{w}_{j, \ell} (x , y ) = w_{ j, \ell}( x, \frac y2)$. 
Then we have $  \tilde{w}_{j, \ell} \in \Xo $ and 
$$
q( \ph,  \tilde{w}_{j, \ell} ) = \ii_0^1 \ii_{\R} d [\ph,  \tilde{w}_{j, \ell}](x, y ) \, dx \, dy 
= 2 \ii_0^{\frac 12} \ii_{\R} d [\ph, {w}_{j, \ell}](x, y ) \, dx \, dy = q ( \ph, w_{j, \ell}) 
$$
because 
$ \ds \ii_0^{\frac 12} \ii_{\R} d [\ph, {w}_{j, \ell}](x, y ) \, dx \, dy = \ii_{\frac 12} ^1 \ii_{\R} d [\ph, {w}_{j, \ell}](x, y ) \, dx \, dy$
(the last equality is a consequence of the fact that $w_{j, \ell}$ is $\frac 12-$periodic with respect to $y$). 
Denoting $ \tilde{\psi}_{j, \ell } ( x, y ) = e^{ i \ph ( x)} + \tilde{w}_{j, \ell}( x, y ) = \psi_{j, \ell}( x, \frac y2 ) $ we have 
$ \tilde{\psi}_{j, \ell } \in \Ep $  and  using (\ref{round2}) we get 
$$
 E_{\la } ( \tilde{\psi}_{j, \ell } )  \ges E_{\la, \min} ( q( \ph,  \tilde{w}_{j, \ell} )) = E_{\la, \min} ( q( \ph,  {w}_{j, \ell} )) 
 \ges E_{\la, \min} (q ) - \sqrt{2 } \e.
 $$

A  simple computation gives
$$
E_{\la } ( \tilde{\psi}_{j, \ell } ) = E_{\la} ( {\psi}_{j, \ell } ) - \frac 34 \la^2 \Big\| \frac{ \p \psi _{j, \ell} }{\p y } \Big\|_{L^2 ( \R \times[0,1])} ^2. 
$$
From the last sequence of inequalities we obtain 
\beq
\label{round3}
\frac 34 \la ^2 \Big\| \frac{ \p \psi _{j, \ell} }{\p y } \Big\|_{L^2 ( \R \times[0,1])} ^2 \les E_{\la} ( {\psi}_{j, \ell } )  -   E_{\la, \min} (q ) + \sqrt{2 } \e. 
\eeq
Summing up the inequalities (\ref{round3}) for $ j, \ell \in \{ 1, 2 \}$ we get 
\beq
\label{round4}
3 \la ^2 \Big\| \frac{\p  \psi  }{\p y } \Big\|_{L^2 ( \R \times[0,1])} ^2 \les 4  E_{\la} ( {\psi} )  -  4 E_{\la, \min} (q ) + 4\sqrt{2 } \e
<  4 \eta_{\e} + 4 \sqrt{2 } \e < 12 \e. 
\eeq
Using (\ref{round4}), the fact that $ \Big\| \frac{ \p \psi }{\p x}\Big\|_{L^2 ( \R \times [0,1])}^2 \les E_{\la } ( \psi ) \les M $ and Lemma \ref{continuous} we infer that 
$$
\big| p( \ph, w ( \cdot, y_2 )) -  p( \ph, w ( \cdot, y_1 ))  \big| \les 2  \Big\| \frac{\p  \psi  }{\p x } \Big\|_{L^2 ( \R \times[0,1])}  \Big\| \frac{\p  \psi  }{\p y } \Big\|_{L^2 ( \R \times[0,1])} 
\les \frac{4 \sqrt{ \e M}}{\la } 
$$
for any $ y_1, y_2 \in [0,1].  $
Since $ y \longmapsto  p( \ph, w ( \cdot, y )) $ is continuous and 
$ \ds \ii_0 ^1  p( \ph, w ( \cdot, y )) \, dy = q$, there exists $ y_* \in [0,1]$ such that $  p( \ph, w ( \cdot, y_* ))  = q$ 
and consequently 
\beq
 p( \ph, w ( \cdot, y ))  \subset \left( q - \frac{ 4 \sqrt{ \e M}}{\la }  , q + \frac{ 4 \sqrt{ \e M}}{\la }  \right)
 \qquad \mbox{ for any } y \in [0,1].
 \eeq
Using the fact that $ E_{\min}^1 $ is $\sqrt{2}-$Lipschitz, we deduce that 
\beq
\label{round5}
E^1 ( \psi (\cdot, y )) = \ii_{\R}  \Big| \frac{\p  \psi  }{\p x } (x, y ) \Big| ^2  + V( |\psi |^2 )(x, y ) \, dx \ges E_{\min}^1 ( q ) - 
\frac{4 \sqrt{ 2 \e M}}{\la } . 
\eeq
Integrating (\ref{round5}) over $[0,1]$ we discover
\beq
\label{round6}
E_{\la, \min}( q ) + \e > E_{\la } ( \psi ) \ges \ii_0^1 E^1 ( \psi ( \cdot, y )) \, dy \ges E_{\min}^1 ( q ) - 
\frac{4 \sqrt{ 2 \e M}}{\la }. 
\eeq
Since (\ref{round6}) holds for any $ \e \in (0, \e_0) $,   
we infer that $E_{\la, \min}( q ) \ges E_{\min}^1 ( q ) .$
Thus necessarily $E_{\la, \min}( q ) = E_{\min}^1 ( q ) $ (see part (ii)) and the proof of step 2 is completed.

 \medskip
 
{\it Step 3. Conclusion. } 
The concavity of $E_{\la, \min}$ on $[0, 2 \pi]$ follows from steps 1 and 2 and from Lemma \ref{concave-HighSchool} below 
with $[a, b ] = [0, 2\pi]$, $ f = E_{\la, \min}$ and $g = E_{\min}^1$. 

\begin{Lemma} 
\label{concave-HighSchool}
Let $f, g : [a, b ] \lra \R $ be two continuous functions. 
Assume that 

a) $g$ is concave and $ f \les g $ on $[a, b]$, and

b) for any $ x \in (a, b ) $ we have either $ f(x ) = g(x )$, or there exists $ \de _x > 0 $ such that 
\beq
\label{middle}
f( x ) \ges \frac 12 ( f( x - \eta ) + f( x + \eta)) \qquad \mbox{  for any  } 0 < \eta < \de_x. 
\eeq

Then $f$ is concave on $[a, b]$. 
\end{Lemma}

{\it Proof. } 
For any $ x \in (a, b ) $ there exists $ \de _x > 0 $ such that (\ref{middle}) holds. 
If $f( x ) < g( x )$, this follows from assumption (b). 
If $ f( x ) = g( x ) $ we may take $ \de _x = \min( x-a, b-x) $. 
Indeed, if $ 0 < \eta < \min(x-a, b-x)$ we have $ x-\eta, x + \eta \in [a, b]$. By assumption (a) we get
$$ f( x ) = g(x ) \ges \frac 12 (  g( x - \eta ) + g( x + \eta)) \ges  \frac 12 ( f( x - \eta ) + f( x + \eta)). $$

Next we see that for any $ \al, \beta \in \R$, the function $ f_{\al, \beta }( t) = f( t ) - \al t - \beta $ cannot achieve a  minimum on
an interval $ ( x_1, x_2 ) \subset  (a, b) $ unless it is constant on $ [x_1, x_2]$. 
Indeed, assume that $ f_{\al, \beta }$ reaches a minimum on $ (x_1, x_2 ) $ at some point $y \in ( x_1, x_2)$. 
It is obvious that $ f_{\al, \beta }$ also satisfies (\ref{middle}). 
Let $ S = \{ z \in ( x_1, x_2 ) \; \mid  \; f_{\al, \beta } ( z ) =  f_{\al, \beta } ( y ) \}$. 
By (\ref{middle}) we see that $ z - \eta, z + \eta \in S$ for any $ 0 \les \eta < \min ( \de_{y}, y - x_1, x_2 - y )$ and we infer that 
$ S $ is open in $ ( x_1, x_2 )$. By the continuity of $ f$, the set $S$ is also relatively closed in $(x_1, x_2)$. Hence $ S = (x_1, x_2)$.
 
Let $ x_1, x_2 \in [a, b ] $, $ x _ 1 < x_ 2$. 
The function $ t \longmapsto f( t ) - \frac{x_2 - t}{x_2 - x_1}f ( x_1 ) - \frac{t - x_1}{x_2 - x_1}f ( x_2 )$ takes the value $0 $ at $ t = x_1 $ and  at $ t = x_2$, hence it must be nonnegative on $[x_1, x_2]$. This means that 
$$
f( t ) \ges \frac{x_2 - t}{x_2 - x_1}f ( x_1 ) + \frac{t - x_1}{x_2 - x_1}f ( x_2 ) \qquad \mbox{ for any } t \in [x_1, x_2]. 
$$
Since $ x_1 $ and $ x_2 $ were arbitrary, $ f$ is concave on $[a, b]$.
\hfill
$\Box$

\medskip

(ix) We already know that $  \la \longmapsto E_{\la, \min} ( p ) $ is nondecreasing and  $ E_{\la, \min} ( p ) \les E_{\min}^1 ( p ) $ for all $ \la > 0 $. 
It suffices to show that, if $ 0 < \la_ 1 < \la _2$ and $ E_{\la_1, \min}( p ) = E_{\la_2, \min}( p )$ then necessarily $ E_{\la_1, \min}( p ) = E_{\min}^1( p )$. 

Let $ M = E_{\la_2, \min}( p ) + 1$. Choose  $ \de \in ( 0, \frac p2) $. 
Let $ \e _ 0 = \min \left( 1, \frac{ \de^2 ( \la_2 ^2 - \la_1 ^2 ) }{4M} \right).$

Take $ \e \in ( 0 , \e_0)$. There exist $ \ph \in \dot{H}^1 ( \R, \R)$ and $ w \in C_c^{\infty} ( \R \times ( 0, 1))$ such that $ q( \ph, w) = p $ and 
$E _{\la _2} ( e^{i \ph } + w ) < E_{\la_2, \min}( p ) + \e$. Denote $ \psi = e^{i \ph } + w$.
It is obvious that $ \big\| \frac{\p \psi }{\p  x} \big\|_{L^2( \R \times (0, 1) )} ^2 \les M$  and 
$ E_{\la _1 } ( \psi  )\ges  E_{\la_1, \min}( p ) = E_{\la_2, \min}( p )$, therefore
$$
( \la_2 ^2 - \la_1 ^2 )  \Big\| \frac{\p \psi }{\p  y} \Big\|_{L^2( \R \times (0, 1) )} ^2
= E_{\la_2}( \psi ) - E_{\la_1}( \psi ) < \e. 
$$
Using (\ref{strans}) we see that the mapping $ y \longmapsto p( \ph, w( \cdot, y ))$ is continuous. We have $ \int_0^1  p( \ph, w( \cdot, y ))  \, dy = q( \ph, w) = p$, hence there exists $ y_0 \in [0,1]$ such that $ p( \ph, w( \cdot, y _0)) = p$. 
From (\ref{strans}) it follows that for any $ y \in [0,1]$, 
$$
\big| p ( \ph, w( \cdot, y )) - p ( \ph, w( \cdot, y _0))\big|
\les 2  \Big\| \frac{\p \psi }{\p  x} \Big\|_{L^2( \R \times (0, 1) )}  \Big\| \frac{\p \psi }{\p  y} \Big\|_{L^2( \R \times (0, 1) )}
\les \frac{ 2 \sqrt{ \e M }}{\sqrt{ \la_2 ^2 - \la_1^2}} < \de. 
$$
Hence $ p ( \ph, w( \cdot, y )) \in ( p - \de , p + \de ) \subset (0, 2 \pi )$ for all $ y $. Since $ E_{\min}^1 $ is $\sqrt{2}-$Lipschitz, we find
$$
E^1 ( \psi( \cdot, y ) ) \ges E_{\min}^1 ( p ( \ph, w ( \cdot, y ))) \ges E_{\min}^1 ( p  ) - \frac{ 2 \sqrt{ 2 \e M }}{\sqrt{ \la_2 ^2 - \la_1^2}} \qquad \mbox{for any } y \in [0,1].
$$
By the choice of $\psi $ and  Fubini's Theorem, we obtain 
$$
E_{\la_1, \min}(p ) + \e > E_{\la_1} ( \psi ) 
 \ges 
\int_0^1 E^1 ( \psi ( \cdot, y )) \, dy 
\ges E_{\min}^1 ( p  ) - \frac{ 2 \sqrt{ 2 \e M }}{\sqrt{ \la_2 ^2 - \la_1^2}} .
$$
Since the above inequality is valid for any $ \e \in (0, \e_0)$ we infer that 
$E_{\la_1, \min}(p ) = E_{\min}^1(p )$. 
\hfill
$\Box$

\medskip

To perform minimization of the energy at fixed momentum in $ \Eo$ we will use a "regularization by minimization procedure" that has been developed in \cite{M10} and \cite{CM}. 
It will enable us to get rid of small defects of Sobolev functions and to approximate functions in $ \Eo $ that have  small energy on every ball of fixed radius    by functions whose modulus is close to $1$. 

We consider a function $ \nu \in C^{\infty } (\R)$ such that $ \nu $ is odd, $ \nu (s) = s $ 
for $ s \in [0, 2  ]$, $ 0 \leq \nu ' \leq 1 $ on $ \R$ and 
$ \nu (s) = 3   $ for $ s \geq 4  $. 
Given  $ \psi \in \Ep $ and $ \la > 0 $, 
the modified  Ginzburg-Landau  energy of $\psi $ in $ \Om $ is 
\beq
E_{GLm , \la } (\psi ) =  \ii_{\R \times [0,1] } \Big| \frac{ \p \psi }{\p x } \Big|^2 + \la ^2 \Big| \frac{ \p \psi }{\p y } \Big|^2 
+ \frac 12  \left( \nu ^2(|\psi |) - 1 \right)^2 \, dx \, dy. 
\eeq
For any given $ \psi \in \Ep, $ and for  $ h , \la > 0 $ we consider the functional 
\beq
\label{ghla}
\begin{array}{rcl}
G_{h, \la }^{\psi } ( \zeta ) & = &  E_{GLm , \la } (\zeta ) + \frac{1}{h^2 } \big\| \zeta - \psi \big\|_{L^2( \R \times [0,1]) }^2 
\\
\\
 & = & \ds
\ii_{\R \times [0, 1]} \Big| \frac{ \p \zeta }{\p x } \Big|^2 + \la ^2 \Big| \frac{ \p \zeta }{\p y } \Big|^2 
+ \frac 12  \left( \nu ^2(|\zeta |) - 1 \right)^2 + \frac{1}{h^2} |\zeta - \psi |^2 \, dx \, dy. 
\end{array}
\eeq
Notice that $ G_{h, \la }^{\psi } ( \zeta ) < \infty $ for any $ \zeta \in \Ep $ satisfying $ \zeta - \psi \in H_{per}^1$. 
We will consider the problem of minimizing $ G_{h, \la }^{\psi } $ in the set 
$ \Ep _{\psi } := \{ \zeta \in \Ep \mid \zeta - \psi \in H_{per}^1 \}$. 

Let $ \Lambda = \frac{ 1}{\la }$. 
Denoting $ \tilde{\psi } ( x, y ) = \psi ( x, \la y )  $ and $ \tilde{\zeta } ( x, y ) = \zeta ( x, \la y )$ we see that 
$ \tilde{ \psi } $ and $ \tilde{ \zeta }$ are $ \Lambda-$periodic with respect to the second variable and 
$$
G_{h, \la }^{\psi } (\zeta )  = \la \ii_{\R \times [0, \Lambda]} 
\Big| \frac{ \p\tilde{ \zeta } }{\p x } \Big|^2 +  \Big| \frac{ \p \tilde{\zeta} }{\p y } \Big|^2 
+ \frac 12  \left( \nu ^2(| \tilde{\zeta} |) - 1 \right)^2 + \frac{1}{h^2} | \tilde{\zeta} - \psi |^2 \, dx \, dy
= \la \tilde{G} _{h}^{\tilde{\psi } } ( \tilde{\zeta}). 
$$
Therefore $ \zeta $ is a minimizer for $G_{h, \la }^{\psi }  $ among $1-$periodic  functions  with respect to the second variable 
if and only if $ \tilde{\zeta}$ is a minimizer for $ \tilde{G} _{h}^{\tilde{\psi } } $ among functions that are 
$\Lambda-$periodic   with respect to the second variable. 
This observation enables us to  use directly the results established in \cite{CM, M10}. 

Proceeding exactly as in Lemma 3.1 p. 160 and Lemma 3.2 p. 164 in \cite{CM} (see also Lemma 3.1 p. 126 and Lemma 3.2 p. 132 in \cite{M10}), we get:

\begin{Lemma}
\label{LemmaVeche}
i) The functional $ G_{h, \la }^{\psi} $ has a minimizer in the set $ \Ep _{\psi} = \{ \zeta \in \Ep \mid \zeta - \psi \in H_{per}^1 \}$. 

ii) Any minimizer $ \zeta_ h $ satisfies 
\beq
\label{6.20}
E_{GLm , \la } (\zeta _ h  ) \les E_{GLm , \la } (\psi ) ; 
\eeq
\vspace{-10pt}
\beq
\label{6.21}
\| \zeta_ h - \psi \|_{L^2( \R \times [ 0 , 1 ] ) }^2 \les h^2 E_{GLm , \la } (\psi ) ; 
\eeq
\vspace{-10pt}
\beq
\label{6.22}
\ds \ii_{ \R \times [ 0 , 1 ]  } \Big\vert \left( \nu ^2( |\zeta_h|) - 1 \right)^2 - 
  \left( \nu ^2( |\psi |) - 1 \right)^2 \Big\vert \, dx \, dy 
\les C   h E_{GLm , \la } (\psi )  . 
\eeq
If $ \psi = e^{ i \ph } + w $ with $ \ph \in \dot{H}^1 ( \R) $ and $ w \in H_{per}^1$,  we have  $ \QR ( \zeta_h ) = \pr{ q( \ph, w + ( \zeta_h - \psi ) )}$ and 
\beq
\label{6.23}
| q ( \ph, w + ( \zeta _ h - \psi )) - q ( \ph, w ) | \les 2 h E_{GLm , \la } (\psi ) .
\eeq

iii) Let $L(z) = \left( \nu ^2( |z |) - 1 \right) \nu ( |z |) \nu '( |z |) 
\frac{z }{ |z |}$  if $ z \in \C ^*$ and $L(0) = 0$. 
Then any minimizer $\zeta_h $ of $ G_{h, \la }^{\psi}  $ in $ \Ep_{\psi} $ 
satisfies   the equation
\beq
\label{6.24}
\ds - \frac{ \partial^2 \zeta _h}{\partial x^2 } - \lambda ^2  \frac{ \partial^2 \zeta _h}{\partial y^2 } +  L( \zeta_h ) 
+ \frac{ 1}{ h^2 } (\zeta _h - \psi ) = 0  \qquad \mbox{ in }  \Do ' ( \R^2). 
\eeq
Moreover, we have $ \zeta_h \in W_{loc}^{2, p } ( \R^2)$ for $ 1 \les p < \infty$,
and consequently $ \zeta_ h \in C_{loc}^{1, \al } ( \R^2)$ for any $ \al \in [0, 1)$. 

iv) For any $ \la > 0$, $ h > 0 $ and  $ \de > 0 $ there exists   a constant 
$K = K ( \la,   h, \de ) > 0 $ such that  for any $ \psi \in \Ep $ satisfying 
$E_{GLm , \la }  ( \psi) \leq K$ and for any minimizer 
$\zeta_h $ of $ G_{h, \la }^{\psi}  $ in $ \Ep_{\psi} $ ~we~have 
\beq
\label{6.25}
1 - \de < |\zeta _h(x) | < 1 + \de 
\qquad \mbox{ in  } \R^2.
\eeq

\medskip

v) Let $ (\psi_n)_{n \ges 1} \subset \Ep $ be a sequence of functions such that 
 $ (E_{GL, \la }(\psi_n))_{n \ges 1} $ is bounded and
\beq
\label{vanishGL}
 \ds \lim_{n \ra \infty } \Big( \sup_{y \in \R } E_{GL ,  \la }^{([y- 1, y + 1])} (\psi_n) \Big) =0.
\eeq

\smallskip

There exists a sequence $ h_n \lra 0 $ such that for any minimizer $ \zeta_n $ of 
$ G_{h_n, \la }^{\psi_n}  $ in $ \Ep_{\psi_n} $  we have 
$ \| \, | \zeta_n | - 1 \|_{L^{\infty }(\R ^2)} \lra 0 $ as $ n \lra \infty $.

\end{Lemma}

The proof of Lemma \ref{LemmaVeche} is the same as the poof of Lemmas 3.1 and 3.2 in \cite{CM}, and we refer the interested reader to that article. See also Lemmas 3.1 and 3.2 in \cite{M10} for higher-dimensional variants. 
The existence of a minimizer is straightforward using the direct method in calculus of variations. 
Estimates (\ref{6.20}) and (\ref{6.21}) folow immediately from the fact that 
$ G_{h, \la }^{\psi } ( \zeta_h ) \les G_{h, \la }^{\psi} ( \psi )  =  E_{GLm , \la } ( \psi) $, and 
(\ref{6.23}) comes from (\ref{diff-mom2}),   (\ref{6.20}) and (\ref{6.21}).
Part (v) is a version of  Lemma 3.2 p. 164 in \cite{CM} in the periodic setting; see also Lemma 3.2 p. 132 in \cite{M10}.

\begin{Lemma}
\label{L6.5}

Fix $\la > 0 $.
Assume that $( \psi_n )_{n \ges 1 } \subset \Ep $ is a sequence of functions satisfying:

(a) There exists $ M > 0 $ such that $ E_{\la} ( \psi _n ) \les M $ for all $n$, 

(b) $ \ds \lim_{n \ra \infty } \Big( \sup_{y \in \R } E_{ \la }^{([y- 1, y + 1])} (\psi_n) \Big) =0, $

(c) $ |\QR ( \psi _n ) | \lra q \in [0, \pi ]$ as $ n \lra \infty$. 

\medskip
\noindent
Then we have $ \ds \liminf_{ n \ra \infty } E_{\la } ( \psi _n ) \ges \sqrt{2} q$. 

\end{Lemma}

{\it Proof.  } 
Using Lemma \ref{MutualBounds} we see that $E_{GL, \la }( \psi_n)$ is bounded and (\ref{vanishGL}) holds. 
We denote 
$$ \tilde{M} = \sup_{n \ges 1 }E_{GL, \la }( \psi_n) \qquad \mbox{ and } \qquad  \e _n =  \sup_{x \in \R } E_{GL ,  \la }^{([x- 1, x + 1])} (\psi_n) \underset{n \ra \infty}{\lra }  0. $$

Let $ u_n = |\psi _n | -1 $. 
Then we have $ |\nabla u _n | \les |\nabla \psi _n | $ almost everywhere on $ \R^2$. 
We also have $ | u _n | \les |\, |\psi _n | - 1 | \cdot | \, |\psi _n  |  + 1 | = |\, |\psi _n |^2 - 1 |$. We infer that 
$ u_n \in H^1 ( \R \times (0, 1 ) ) $ and 
$$
 \| u _n \|_{  H^1 ( \R \times (0, 1 ) ) } ^2 \les \max \left( 1 , \frac{1}{\la ^2} \right) E_{GL, \la } ( \psi _n ) 
\les \max \left( 1 , \frac{1}{\la ^2} \right) \tilde{M}. 
$$
Similarly we get 
$ \| u _n \|_{  H^1 ( [x-1, x + 1 ]  \times (0, 1 ) ) } ^2 \les \max \left( 1 , \frac{1}{\la ^2} \right) \e _n $ for any $ x \in \R$. 
Let $ p \in (2, \infty)$. By the Sobolev embedding, there is $ C_p > 0 $ such that 
$ \| u \|_{L^p ( [ a-1, a + 1 ]\times [ 0, 1 ] ) }\les C_p \| u \|_{H^1( ( a-1, a + 1 ) \times ( 0, 1 ) )} $ for any 
$ u \in {H^1( ( a-1, a + 1 ) \times ( 0, 1 ) )} $ and we infer that 
$$
\ii_{ [ a-1, a + 1 ]\times [ 0, 1 ] } | u _n |^p \, dx \, dy \les C_p ^p  \| u \|_{H^1( ( a-1, a + 1 ) \times ( 0, 1 ) )}^p 
\les C( p, \la ) 
\e_n ^{ \frac{p-2}{2}} \| u _n \|_{H^1( ( a-1, a + 1 ) \times ( 0, 1 ) )}^2. 
$$
Take a sequence $ (a_k)_{k \ges 1 }$ such that $ \ds \R = \cup_{k \ges 1 } [ a_k - 1 , a _k + 1 ] $ and each point $ x \in \R$ belongs to at most two of the intervals $ [a_k - 1, a_k + 1 ] $. 
We write the above inequality for each $k$ and we sum over $k$  to get 
$$
\| u _n \|_{L^p (\R \times [0,1]) }^p \les C( p, \la ) \e_n ^{  \frac{p-2}{2}} \| u_n \|_{H^1 ( \R \times (0, 1) )}^2 \les 
C( p, \la ) \e_n ^{  \frac{p-2}{2}}  E_{GL, \la } ( \psi _n ) \les C( p, \la ) \e_n ^{  \frac{p-2}{2}}  \tilde{M}. 
$$
We have thus proved  that $ \| u _n \|_{L^p (\R \times [0,1]) } \lra 0 $ as $ n \lra \infty$ for any $ p \in (2, \infty).$
 
Next we show that 
\beq
\label{approx-nonlin}
\ii_{\R \times [ 0, 1 ] } \Big| V( |\psi _n |^2 ) - \frac 12 \left( |\psi _n |^2 - 1 \right)^2 \Big| \, dx \, dy \lra 0 \qquad \mbox{ as } n \lra \infty. 
\eeq

Fix $ \e > 0 $. 
By assumption {\bf (A1)} there exists $ \eta ( \e ) > 0 $ such that 
$$
 \Big| V( s^2 ) - \frac 12 \left( s^2 - 1 \right)^2 \Big| \les \frac{\e}{2} 
\left( s^2 - 1 \right)^2 \qquad \mbox{  for any } s \in [ 1 - \eta ( \e), 1 + \eta( \e)].
$$ 
Choose $ p \ges \max( 4, 2p_0 +2)$, where $p_0 $ is as in {\bf (A2)}.  
By  assumption {\bf (A2)} there is $ C(\e, p ) > 0 $ such that 
$$
 \Big| V( s^2 ) - \frac 12 \left( s^2 - 1 \right)^2 \Big| \les C( \e, p ) \big| \, |s | - 1 \big|^p \qquad \mbox{ for any } s \in [0, 1 - \eta(\e)] \cup [1 + \eta( \e ) , \infty).
$$ 
We find 
$$
 \Big| V( |\psi _n |^2 ) - \frac 12 \left( |\psi _n |^2 - 1 \right)^2 \Big|  \les 
 \frac{\e}{2}  \left( |\psi _n |^2 - 1 \right)^2 + C( \e, p ) | \, |\psi_n | - 1 |^p 
$$
and integrating we obtain 
\beq
\label{V-GL}
\begin{array}{l}
\ds \ii_{\R \times [ 0, 1 ] } \Big| V( |\psi _n |^2 ) - \frac 12 \left( |\psi _n |^2 - 1 \right)^2 \Big| \, dx \, dy  
\\
\\
\ds \les \frac{\e}{2} \ii_{\R \times [ 0, 1 ] } \left( |\psi _n |^2 - 1 \right)^2 \, dx \, dy  
+ C( \e, p ) \| u _n \|_{L^p (\R \times [0,1]) }^p 
\les \e \tilde{M} + C( \e, p, \la ) \e_n ^{  \frac{p-2}{2}}  \tilde{M}.
\end{array}
\eeq
Since $ \e _n \lra 0$, 
there exists $ n( \e) \in \N^* $ such that $ C( \e, p, \la ) \e_n ^{  \frac{p-2}{2}} < \e $ for all $ n \ges n( \e)$. 
Then the right-hand side in the above inequality is smaller than $ 2 \e \tilde{M}$ for all $ n \ges n( \e)$. 
Since $ \e $ was arbitrary, (\ref{approx-nonlin}) is proven. 

Assume that $ \zeta \in \Ep $ satisfies $ 1 - \de \les  | \zeta | \les 1 + \de$ for some $ \de \in (0, 1)$. 
According to Lemma \ref{lift}, $ \zeta$ admits a lifting $ \rho e^{ i \theta }$, a valuation of the momentum of $ \zeta $ is 
$  \ii_{\R \times [0, 1] } ( 1 - \rho ^2 ) \frac{ \p \theta}{\p x } \, dx dy $,
and we have
\beq
\label{standard}
\begin{array}{l}
\ds \sqrt{2} ( 1 - \de ) \bigg|   \ii_{\R \times [0, 1] } ( 1 - \rho ^2 ) \frac{ \p \theta}{\p x } \, dx \, dy \bigg| 
\les   \ii_{\R \times [0, 1] }( 1 - \de )^2 \Big|  \frac{ \p \theta}{\p x } \Big|^2 + \frac 12 ( 1 - \rho ^2 )^2  \, dx \, dy
\\
\\
\ds \les  \! \! \ii_{\R \times [0, 1] } \! \rho ^2\Big|  \frac{ \p \theta}{\p x } \Big|^2 \! \! + \frac{ ( 1\! - \rho ^2 )^2 }{2} \, dx \, dy
\les  \! \ii_{\R \times [0, 1] } \Big|  \frac{ \p \zeta}{\p x } \Big|^2 \! \! + \frac 12 ( 1\! - |\zeta | ^2 )^2  \, dx \, dy \les E_{GLm, \la } ( \zeta ).  
\end{array}
\eeq

We may use Lemma \ref{LemmaVeche} (v) for the sequence $(\psi_{n})_{n \ges 1}$. 
We infer that there exists a sequence $ h_n \lra 0 $ and for each $n$ there is a 
minimizer $ \zeta_n $ of 
$ G_{h_n, \la }^{\psi_n}  $ in $ \Ep_{\psi_n} = \{ \zeta \in \Ep \mid \zeta - \psi _n \in H_{per}^1 \}$  such that 
$ \| \, | \zeta_n | - 1 \|_{L^{\infty }(\R ^2)} \lra 0 $ as $ n \lra \infty $. Denote $ \de _n = \| \, | \zeta_n | - 1 \|_{L^{\infty }(\R ^2)}$, 
so that $ 1 - \de_n \les |\zeta_n | \les 1 + \de_n . $
For all $ n $ sufficiently large we have $ \de _n < 1 $,  and then $ \zeta _n $ admits a lifting $ \zeta _ n = \rho_n e^{ i \theta _n}$ 
and (\ref{standard}) holds for $ \zeta_n $. 

From (\ref{6.23}) we have 
$$
|\QR ( \zeta_n ) - \QR( \psi_n) | \les 2 h_n E_{GLm, \la } ( \psi _n ) \les 2 h_n \tilde{M} \lra 0 \quad \mbox{ as } n \lra \infty .
$$
Since $ |\QR ( \psi _n ) | \lra q$,  we infer that $ |\QR ( \zeta_n ) | \lra q$ as $ n \lra \infty$. 

We have $ E_{GL, \la } ( \psi _n ) \ges E_{GLm, \la } ( \psi _n ) \ges E_{GLm, \la } ( \zeta _n )$,  and using (\ref{standard}) we obtain
$$
\begin{array}{l}
E_{\la} (\psi _n) = E_{GL, \la } (\psi _n) + \ds  \int_{\R \times[0,1]} V(|\psi_n |^2) - \frac 12 ( |\psi _n|^2 - 1)^2\, dx \, dy 
\\
\\
\geq 
E_{GLm, \la} (\zeta _n) + \ds  \int_{\R \times [0,1]} V(|\psi_n |^2) - \frac 12 ( |\psi _n|^2 - 1)^2\, dx \, dy 
\\
\\
\geq 
\sqrt{2} ( 1 - \de _n ) |\QR (\zeta _n )| + \ds  \int_{\R \times [0,1]} V(|\psi_n |^2) - \frac 12 ( |\psi _n|^2 - 1)^2\, dx .
\end{array}
$$
Letting $ n \lra \infty $ in the above inequality and using (\ref{approx-nonlin}) we get
$ \ds \liminf_{n \ra \infty } E_{\la} (\psi _n)  \ges \sqrt{2}    q$ and Lemma \ref{L6.5} is proven.
\hfill
$\Box$

\medskip

{\it Proof of Lemma \ref{Emin} (v).}
Fix $ p \ges \max( 4, 2 p_0+2 )$, where $p_0 $ is as in {\bf (A2)}.
Coming back to (\ref{V-GL}) and using the fact that $ \| |\psi | -1 \|_{H^1 ( \R \times (0, 1))}^2 \les C( \la ) E_{GL, \la } ( \psi )$ and the Sobolev embedding, we see that for any $ \e > 0 $ there exists $ C (\e, p ) > 0 $ such that for any $ \psi \in \Ep $ we have
\beq
\label{V-GLbis}
\begin{array}{l}
\big| E_{\la } ( \psi ) - E_{GL, \la } ( \psi ) \big| \les \ds \ii_{\R \times [ 0, 1 ] } \big| V( |\psi _n |^2 ) - \frac 12 \left( |\psi _n |^2 - 1 \right)^2 \big| \, dx \, dy  
\\
\\
\ds \les \frac{\e}{2} \ii_{\R \times [ 0, 1 ] } \left( |\psi _n |^2 - 1 \right)^2 \, dx \, dy  
+ C( \e, p ) \big\|   |\psi _n  | - 1 \big\|_{L^p (\R \times [0,1]) }^p 
\\
\\
\les \e E_{GL, \la} ( \psi )  + C( \e, p, \la )  E_{GL, \la} ( \psi )  ^{  \frac{p}{2}}  .
\end{array}
\eeq
 Using  (\ref{V-GLbis}) 
we infer that for any $ \e > 0 $ and for any $ \la > 0 $  there exists $ \tilde{M} ( \e, \la ) > 0 $ such that 
\beq
\label{V-GLter}
\big| E_{\la } ( \psi ) - E_{GL, \la } ( \psi ) \big| \les  2 \e E_{GL, \la} ( \psi )   \qquad \mbox{ for any } \psi \in \Ep \mbox{ satisfying } 
E_{GL, \la } ( \psi ) \les \tilde{M} ( \e, \la ). 
\eeq

Fix $\de \in (0, 1)$. 
Given any $ \psi \in \Ep $ and any $ h > 0 $,  denoting by $ \zeta _ h $ a minimizer of $ G_{h, \la } ^{\psi } $ in $ \Ep _{\psi}$ 
we have $ E_{GL, \la } ( \psi ) \ges E_{GLm, \la }( \psi )\ges E_{GLm, \la } ( \zeta _h) $ and therefore
\beq
\label{MOAB}
\begin{array}{l}
\ds E_{\la}( \psi ) - \sqrt{2} ( 1 - \de ) \big| \QR ( \psi ) \big| 
\ges \left( E_{\la}( \psi ) - E_{GL, \la}( \psi ) \right) 
+ \frac{ \de}{2} E_{GL, \la}( \psi ) 
\\ \\
\ds + \left(\!  1 - \frac{ \de}{2} \right) \! \! \left( \! E_{GLm, \la } ( \zeta_h) - \frac{1 - \de}{1 - \frac{ \de}{2}} \sqrt{2} \big| \QR ( \zeta _h) \big| \right) + ( 1 - \de ) \sqrt{2} \left(  \big| \QR ( \zeta _h) \big| - \big| \QR ( \psi ) \big|  \right). 
\end{array}
\eeq

Choose $ h = \frac{\de}{16\sqrt{2} ( 1 - \de)}$.  
Using (\ref{6.23}) we have 
\beq
\label{approximQ}
( 1 - \de ) \sqrt{2} \Big|  \big| \QR ( \zeta _h) \big| - \big| \QR ( \psi ) \big|  \Big| \les ( 1 - \de ) \sqrt{2} \cdot ( 2 h E_{GLm, \la } ( \psi ) )
\les \frac{ \de}{8 } E_{GL, \la } ( \psi ) .
\eeq

If $ h $ is as above, by Lemma \ref{LemmaVeche} (iv) there exists a constant $ K(\la,  \de) > 0 $ such that for any $ \psi \in \Ep $ satisfying 
$ E_{GL, \la} ( \psi ) \les K(\la,  \de)$ and for any minimizer $ \zeta_h $ of $G_{h, \la }^{\psi} $ in $ \Ep_{\psi}$ we have
$ 1 - \frac{\de}{2} \les |\zeta_h | \les 1 + \frac{\de}{2}. $
Then using (\ref{standard}) for $ \zeta_h$ we see that 
\beq
\label{6.34}
\frac{1 - \de}{1 - \frac{ \de}{2}} \sqrt{2} \big| \QR ( \zeta _h) \big| \les \left(1 - \frac{ \de}{2} \right) \sqrt{2} \big| \QR ( \zeta _h) \big|
\les E_{GLm, \la } ( \zeta_h ). 
\eeq

If $ E_{GL, \la }( \psi ) \les \tilde{M}\left( \frac{\de}{16}, \la \right)$, using (\ref{V-GLter}) we get 
\beq
\label{6.35}
\big| E_{\la } ( \psi ) - E_{GL, \la } ( \psi ) \big| \les  \frac{\de}{8} E_{GL, \la} ( \psi )  .
\eeq

By Lemma \ref{MutualBounds} there exists $ m > 0 $ such that for any $ \psi \in \Ep $ satisfying $ E_{\la }( \psi ) \les m $ we have 
$ E_{GL, \la } ( \psi ) \les \min \left(  K(\la,  \de), \tilde{M}\left( \frac{\e}{16}, \la \right) \right)$. 
Then using (\ref{MOAB})-(\ref{6.35}) we infer that 
$$
\ds E_{\la}( \psi ) - \sqrt{2} ( 1 - \de ) \big| \QR ( \psi ) \big|  \ges \frac{ \de}{4} E_{GL, \la } ( \psi ) \ges 0 
\qquad \mbox{ for all }  \psi \in \Ep \mbox{ satisfying } E_{\la}( \psi ) \les m . 
$$
The above inequality and the fact that $ E_{\la, \min} ( p ) \lra 0 $ as $ p \lra 0 $ imply that there exists $ p_{\de } > 0 $ such that 
$ E_{\la, \min}( p ) \ges \sqrt{2} ( 1 - \de ) p $ for all $ p \in [0, p_{\de} ]$. 
\hfill
$\Box$

\begin{Theorem}
\label{T6.6}
Assume that $V$ satisfies {\bf (A1)}, {\bf (A2)}, and {\bf (B2)}. 
Let $ p \in (0, \pi]$ and $ \la > 0 $ such that  $ E_{\la, \min} ( p ) < \sqrt{2} p $. 
Let $ ( \psi_n )_{n \ges 1} \subset \Ep $ be a sequence satisfying 
\beq
\label{conv0-per}
\QR ( \psi _n ) \lra \pr{p} \qquad \mbox{ and } \qquad E_{\la} ( \psi _n ) \lra E_{\la , \min} ( p) \qquad \mbox{ as } n \lra \infty. 
\eeq
Then there exist a subsequence $(\psi_{n_k})_{k \ges 1}$, a sequence $ (x_k )_{k \ges 1} \subset \R $ and $ \psi \in \Ep$ satisfying 
$\QR ( \psi ) = \pr{p}$, $E_{\la}( \psi ) = E_{\la, \min}( p )$, and
$$
\begin{array}{l}
|\psi_{n_k} | ( \cdot + x_k ) - | \psi| \lra 0  \quad \mbox{ in } L^p ( \R \times[0,1] ) \; \mbox{ for } 2 \les p < \infty, 
\\
V ( |\psi_{n_k} ( \cdot + x_k )|^2 )  \lra V( |\psi |^2 )  \quad \mbox{ in } L^1 ( \R \times [0, 1]),  
\\
\nabla \psi_{n_k} ( \cdot + x_k ) \lra  \nabla \psi \quad \mbox{ in } L^2 ( \R \times [0, 1]). 
\end{array}
$$
\end{Theorem}

\begin{remark} 
\rm 
If $ \la \ges \la_s ( q )$ we have $ E_{\la, \min}( q) = E_{\min}^1 ( q ) $. In this case the existence of minimizers follows from 
Theorem \ref{T5.2} in the previous section. 
However, even in this case Theorem \ref{T6.6} above is interesting because it gives the stability of the minimizers under two-dimensional periodic perturbations, and this stability does not follow directly  from Theorem \ref{T5.2}. 

\end{remark}

{\it Proof of Theorem \ref{T6.6}. } 
Let $(\psi_{n})_{n \ges 1} \subset \Ep $ be a sequence satisfying (\ref{conv0-per}). 
Then $ E_{\la } ( \psi _n )$ is bounded.

As in the proof of Theorem \ref{T5.2}, we  use the concentration-compactness principle (\cite{lions}). 
We denote by $ \Lambda _n $ the concentration function of 
$ f_n := \ds \Big| \frac{ \p \psi _n }{\p x } \Big|^2 + \la ^2 \Big| \frac{ \p \psi _n }{\p y } \Big|^2 + V( |\psi _n |^2  )$, namely
$$
\Lambda _n ( t ) = \sup_{a \in \R} \ii_{[a-  t, a + t] \times [0, 1 ]  }   \Big| \frac{ \p \psi _n }{\p x } \Big|^2 + \la ^2 \Big| \frac{ \p \psi _n }{\p y } \Big|^2 + V( |\psi _n |^2 ) \, dx \, dy. 
$$
It is clear that  $ \Lambda _n $ is a non-decreasing function on $[0, \infty)$, 
$ \Lambda _n ( 0 ) = 0 $ and $ \ds \lim_{t \ra \infty} \Lambda _n ( t ) = E_{\la }  ( \psi _n ) $. 
By a standard argument (see \cite{lions}),  there exists a subsequence of $( \psi_n , \Lambda _n )_{n \ges 1}$, still denoted $( \psi_n , \Lambda _n )_{n \ges 1}$,
and there is a non-decreasing function $ \Lambda : [0, \infty ) \lra [0,\infty)  $  satisfying  
\beq
\label{cc2}
\Lambda _n (t) \lra \Lambda (t) \mbox{ a.e on } [0, \infty) \mbox{ as } n \lra \infty . 
\eeq
Let $ \al = \ds \lim_{t \ra \infty} \Lambda (t)$. We have $ \ds 0 \les \al \les \lim_{n \ra \infty } E_{\la }  ( \psi _n )=  E_{\la, \min} ( p ).$

We will show that $ \al = E_{\la, \min} ( p ) $. 

\medskip

We prove first that $ \al > 0 $. We argue by contradiction and we assume that $ \al = 0$. 
Then we have $ \Lambda ( t ) = 0 $ for all $ t \ges 0 $, and in particular $ \Lambda ( 1 ) = 0 $. 
Then Lemma \ref{L6.5} implies that $ \ds \liminf_{n \ra \infty } E_{\la } ( \psi _n ) \ges \sqrt{2}  \lim_ { n \ra \infty } \big| \QR ( \psi _n ) \big| $, that is $ E_{\la, \min} ( p ) \ges \sqrt{2} p$, contradictiong the assumption of Theorem \ref{T6.6}. 

\medskip

Assume that $ 0 < \al < E_{\la, \min} ( p ) $.
Proceeding as in the proof of Theorem \ref{T5.2} (see (\ref{dic1}) - (\ref{dic4}) there),  we see that there exist a sequence 
$ R_n \lra \infty$, a sequence $ ( a_n )_{n \ges 1 } \subset \R$, and  $ \beta \in (0, E_{\la, \min}( p ) ) $ such that, 
after possibly extracting a further subsequence, we have 
\beq
\label{dic11}
E_{\la }^{( - \infty, a_n  - {R_n }) }  ( \psi _n ) \lra \beta, \qquad 
E_{\la }^{(  a_n  + {R_n }, \infty ) }  ( \psi _n ) \lra E_{\la, \min}( p ) - \beta, \quad \mbox{ and }
\eeq
\beq
\label{dic12}
E_{\la }^{[ a_n  - {R_n }, a _ n +{R_n} ] }  ( \psi _n ) \lra 0.
\eeq
Let $ \breve{\psi}_n ( x  ) = \ii_0 ^1 \psi_n ( x, y ) \, dy$ and $ v_{\psi_n} ( x, y ) = \psi_n (x, y ) - \breve{\psi}_n ( x)$. 
From (\ref{dic12}) and  Lemma \ref{PG-bis} we infer that 
\beq
\label{dic13}
E^{1, \, [a_n - R_n, a_n + R_n ] }( \breve{\psi}_n ) \lra 0 \quad \mbox{ and } \qquad 
\|  v_{\psi_n} \|_{H^1( ([a_n - R_n, a_n + R_n ] ) }\lra 0 
\quad \mbox{ as } n \lra \infty. 
\eeq
Proceeding as in the proof of Theorem \ref{T5.2} we see that 
$\big\| |\breve{\psi} _n | -1 \big\|_{L^{\infty}( [a_n - R_n +1, a_n + R_n -1] } \lra 0 $. For $n $ sufficiently large, $\breve{\psi} _n$ admits a lifting $ \breve{\psi}_n (x) = \rho_n (x) e^{ i \theta _n (x) }$ on 
$  [a_n - R_n +1, a_n + R_n -1] $, where $ \rho_n = |\breve{\psi}_n|$. 
Since $E^{1, \, [a_n - R_n, a_n + R_n ] }( \breve{\psi}_n ) \lra 0$ we have 
$ \| \rho_n - 1 \|_{H^1( [a_n - R_n+1 , a_n + R_n -1] ) } \lra 0 $. 

By Lemma \ref{PG} there exist $ \ph_n \in \dot{H}^1( \R)$ and $ v _n \in H^1( \R)$ such that 
$ \breve{\psi}  _n = e^{ i \ph_n } + v_n $ on $ \R$.  
We may assume that $ \ph _n = \theta_n $ and $ v_n = ( \rho_n - 1) e^{ i \theta_n } $ on $ [a_n - R_n +2, a_n + R_n - 2 ]$, where $ \breve{\psi}_n = \rho_n e^{ i \theta _n}$ is the above lifting
(to see this, we take $ \chi_n \in C_c^{\infty } ( \R) $ such that $ \mbox{supp}(\chi_n ) \subset [a_n - R_n +1, a_n + R_n -1]$ and
$ \chi_n = 1 $ on $[a_n - R_n+2,  a_n + R_n - 2 ]$ and we replace 
$ \ph _n $ by $ \tilde{\ph}_n = ( 1 - \chi_n ) \ph _ n + \chi_n \theta _n$ and $ v_n $ by 
$ \tilde{v}_n  = \breve{\psi}_n - e^{ i \tilde{\ph}  _n } $). 
Then we may write 
$$
\psi _ n = e^{ i \ph _n (x) } + w_n ( x, y ) \qquad \mbox{ on } \R ^2, 
$$
where $ \ph_n ( x ) = \theta_n ( x ) $ and $ w_n ( x, y ) = ( \rho_n ( x) - 1 ) e^{ i \theta_n (x) } + v_{\psi_n }( x, y )$ on 
$ [a_n - R_n + 2 , a_n + R_n - 2 ] $. 
We have
\beq
\label{dic14}
\| \ph _n ' \|_{L^2( [a_n - R_n + 2 , a_n + R_n - 2 ] )} \lra 0 \quad \mbox{ and } \quad 
\| w_n \|_{H^1 ( [a_n - R_n + 2 , a_n + R_n - 2 ]  \times [0,1]) } \lra 0 
\eeq
as $ n \lra \infty$. 
Consider a nondecreasing function $ \chi \in C^{\infty } ( \R)$ such that $ \chi  = 0 $ on $( - \infty, -1 ]$ and $ \chi = 1 $ on $[0, \infty)$. 
Define
$$
\ph_{n, 1 }( x ) = \left\{ 
\begin{array}{ll} \ph_n ( x ) & \mbox{ if } x \les a_n, \\
\ph_n ( a_n  ) & \mbox{ if } x  > a_n, 
\end{array}
\right. 
\qquad 
\ph_{n, 2 }( x ) = \left\{ 
\begin{array}{ll} \ph_n ( a_n ) & \mbox{ if } x < a_n, \\
\ph_n ( a_n  ) & \mbox{ if } x  \ges a_n, 
\end{array}
\right. 
$$
$$
w_{n, 1} (x, y ) = \chi( a_n - x  ) w_n ( x, y ), \qquad w_{n, 2 } ( x, y ) = \chi( x - a_n )w_n ( x, y ), 
$$
and 
$\psi_{n, j} (x, y ) = e^{ i \ph _{n, j } ( x)} + w_{n, j } ( x, y )$ for $ j = 1, 2$. 
It is then clear that $ \psi_{n, j } \in \Ep$ and we have 
\beq
\label{dic15}
\begin{array}{ll}
\psi_{n, 1 } = \psi_n \mbox{ on } (- \infty, a_n],  &  \psi_{n, 1} = e^{ i \ph_n ( a_n) } = constant \mbox{ on } [a_n + 1, \infty), 
\\ 
\psi_{n, 2 } = \psi_n \mbox{ on } [a_n,  \infty),  & \psi_{n, 2} = e^{ i \ph_n ( a_n) } = constant \mbox{ on } (-\infty,  a_n -1 ].
\end{array}
\eeq
Using (\ref{dic11}), (\ref{dic12}) and (\ref{dic14}) it is easily seen that 
\beq
\label{dic16}
E_{\la } ( \psi_{n, 1 } ) \lra \beta \qquad \mbox{ and } \qquad E_{\la } ( \psi_{n, 2 } ) \lra E_{\la, \min} (p ) - \beta
\qquad \mbox{ as } n \lra \infty. 
\eeq
Taking into account (\ref{dic15}), valuations of the momenta of $ \psi_{n, 1}$ and of $ \psi_{n, 2}$ are, respectively, 
$$
\begin{array}{rcl}
q( \ph_{n, 1}, w_{n,1} )  & = & \ds \ii_{( - \infty, a_n ] \times [0,1]}
-2 \langle \ph_n ' e^{ i \ph_n}, w_n \rangle + \langle i \frac{ \p w_n}{\p x } , w_n \rangle \, dx \, dy 
\\ \\  
 & &  \quad \ds +  \ii_{[ a_n , a_n +1] \times [0,1]} \langle i \frac{ \p w_{n,1}}{\p x } , w_{n,1} \rangle \, dx \, dy, \quad \mbox{ and}
\end{array}
$$
$$
\begin{array}{rcl}
q( \ph_{n, 2}, w_{n,2} ) & = & \ds \ii_{[a_n,   \infty )\times [0,1]}
-2 \langle \ph_n ' e^{ i \ph_n}, w_n \rangle + \langle i \frac{ \p w_n}{\p x } , w_n \rangle \, dx \, dy  
\\
\\
& &   \quad \ds  + \ii_{[ a_n-1 , a_n ] \times [0,1]} \langle i \frac{ \p w_{n,2}}{\p x } , w_{n,2} \rangle \, dx \, dy.  
\end{array}
$$
Since $ R_n \lra \infty $ and $\| w_{n , j } \|_{H^1 ( [a_n - R_n + 2 , a_n + R_n - 2 ]  \times [0,1]) } \lra 0$, it follows from the above that 
\beq
\label{dic17}
q( \ph_{n, 1}, w_{n,1} )  + q( \ph_{n, 2}, w_{n,2} ) = q( \ph_n , w_n ) + o(1) \qquad \mbox{ as } n \lra \infty.
\eeq
After extracting  a further subsequence we may assume that $ \big| \QR ( \psi_{n, 1} ) \big| \lra p_1 \in [0, \pi] $ and 
$ \big| \QR ( \psi_{n, 2} ) \big| \lra p_2 \in [0, \pi] $. Then from (\ref{dic17}) and the fact that $ \big| \pr{a + b } \big| \les \big|\pr{a} \big| +  \big|\pr{b} \big|$, we  get 
\beq
\label{dic18}
p_1 + p_2 \ges p . 
\eeq
We have $ E_{\la} ( \psi_{n, j } ) \ges E_{\la, \min } \left(  \big| \QR ( \psi_{n, j} ) \big| \right) $. Passing to the limit and using (\ref{dic16}) we see that 
\beq
\label{dic19}
\beta \ges  E_{\la, \min } ( p_1) \qquad \mbox{ and } \qquad E_{\la, \min} ( p ) - \beta \ges  E_{\la, \min } ( p_2). 
\eeq
Since $E_{\la, \min} $ is nondecreasing on $[0, \pi]$,  
 (\ref{dic19}) implies that $ p_1 < p $ and $ p_2 < p$  and then from (\ref{dic18}) we infer that $ 0 < p_j < p $ for $ j = 1,2$. 

The concavity of $E_{\la, \min}$ implies that $ E_{\la, \min}( p_j ) \ges \frac{ p_j}{p} E_{\la, \min}(p)$ for $ j = 1,2$, and equality may occur if and only if 
$ E_{\la, \min} $ is linear on $[0, p]$. 
Summing up these two inequalities we find $ E_{\la, \min} ( p_1 ) + E_{\la, \min}( p_2 ) \ges \frac{p_1 + p_2}{p}E_{\la, \min}(p) \ges E_{\la, \min}(p)$. Comparing to (\ref{dic19}) we see that we must have equality in the two previous inequalities,  and this implies that $ p_1 + p_2 = p $ and $ E_{\la, \min}( p_j ) =  \frac{ p_j}{p} E_{\la, \min}(p)$, hence  $E_{\la, \min}$ must be linear on $[0, p]$. 
Then using  Lemma \ref{Emin} (v) we deduce that $ E_{\la, \min} ( p ) = \sqrt{ 2} p $,  contradicting the assumption that 
$ E_{\la, \min} ( p ) < \sqrt{ 2} p$. 
We conclude that we cannot have $ 0 < \al < E_{\la, \min}( p )$. 

\medskip

So far we have proved that $ \al =  E_{\la, \min}( p )$.
Then it is standard to prove that there  is a sequence $(x_n)_{n \geq 1} \subset \R^N$ such that 
for any $ \e > 0 $ there is  $ R_{\e} > 0 $  satisfying 
$E_{\la} ^{  ( - \infty, x_n -  R_{\e})} (\psi _n)  + E_{\la} ^{  (  x_n +  R_{\e}, \infty )} (\psi _n) < \e$ for all sufficiently large $ n$. 
Denoting $ \tilde{\psi}{_n} = \psi _n ( \cdot + x_n)$, we see that 
for any $ \e > 0 $ there exist   $ R_{\e} > 0 $ and $ n_{\e } \in \N$ such that 
\beq
\label{conv21} 
E_{\la} ^{( - \infty, - R_{\e})}(\tilde{\psi}{_n})  + E_{\la} ^{( R_{\e}, \infty )}(\tilde{\psi}{_n}) < \e \qquad \mbox{ for all } n \geq n_{\e}.
\eeq
Obviously,  $ (\nabla \tilde{\psi}_{n})_{n \ges 1} $ is bounded in $L^2(\R \times (0, 1)) $ and using the boundedness of 
$E_{GL, \la } ( \tilde{\psi}_n)$ it is easy  to see that
$( \tilde{\psi}{_n})_{n \ges 1} $ is bounded in $L^2 ((-R, R) \times (0, 1) )$ for any $R >0$.
By a standard argument,  there exist a function $ \psi \in H_{loc}^1(\R^2)$
which  is $1-$periodic with respect to the second variable, $ \nabla \psi \in L^2(\R \times (0, 1))$ and there is a subsequence 
$ (\tilde{\psi}_{n _k})_{k \geq 1} $
 satisfying 
 \beq
 \label{conv22}
 \begin{array}{ll}
 \nabla \tilde{\psi}_{n_k } \rightharpoonup \nabla \psi & \quad \mbox{ weakly in } L^2( \R \times (0, 1)), 
 \\
 \tilde{\psi}_{n_k } \rightharpoonup  \psi & \quad \mbox{ weakly in } H^1( ( - R, R ) \times (0,1) ) \mbox{ for any } R > 0, 
 \\
 \tilde{\psi}_{n_k } \lra  \psi \; & \quad \mbox{ strongly in } L^p(( - R, R ) \times (0,1) ) \mbox{ for any } R > 0 \mbox{ and }
  p \in [1, \infty) , 
 \\
  \tilde{\psi}_{n_k } \lra  \psi  & \quad   \mbox{  almost everywhere on } \R^2.
 \end{array}
 \eeq

The weak convergence implies that  for any interval $ I \subset \R$ we have 
\beq
\label{conv23}
\int_{ I \times (0, 1) } \Big| \frac{ \p  \psi }{\p x } \Big|^2  + \la ^2 \Big| \frac{ \p  \psi }{\p y } \Big|^2\, dx \, dy 
\les \liminf_{k \ra \infty } 
\int_{ I \times (0, 1)} \Big| \frac{ \p  \tilde{\psi}_{n_k } }{\p x } \Big|^2  + \la ^2 \Big| \frac{ \p  \tilde{\psi}_{n_k } }{\p y } \Big|^2 \, dx \,  dy .
\eeq
The almost everywhere  convergence and Fatou's Lemma give
\beq
\label{conv24}
\int_{ I \times [0, 1] }
\left( |\psi | ^2 - 1 \right)^2 \, dx \, dy 
\leq \liminf_{k \ra \infty } 
\int_{ I \times [0, 1] } \left(  |\tilde{\psi }_{n_k}| ^2 - 1 \right)^2 \, dx \, dy
\quad \mbox{ and } 
\eeq
\beq
\label{conv25}
\quad
\int_{ I \times [0, 1] }  V(|\psi |^2)  \, dx  \, dy \leq \liminf_{k \ra \infty } \int_{I \times [0, 1]}  V(|\tilde{\psi} _{n_k}|^2)  \, dx \, dy.
\eeq
From (\ref{conv23}) - (\ref{conv25}) we get
\beq
\label{conv26}
E_{GL, \la} ^ I ( \psi ) \leq \liminf_{k \ra \infty }  E_{GL, \la}^I( \tilde{\psi } _{n_k}) 
\quad
\mbox{ and } 
\quad
E_{\la}^I ( \psi ) \leq \liminf_{k \ra \infty }  E_{\la}^I( \tilde{\psi } _{n_k}) .
\eeq
We deduce   in particular that $ \psi \in \Ep $, 
$ E_{\la } ( \psi ) \les E_{\la, \min}( p ) $  and $ E_{\la }^{( - \infty, - R_{\e})} ( \psi ) +  E_{\la }^{(  R_{\e}, \infty )} ( \psi ) \les \e$.

Let $ \breve{\tilde{\psi}}_{n _k } = \int_0^1 \tilde{\psi}_{n _k } (x, y ) \, dy $. 
Using Lemma \ref{PG-bis} (iii)  as well as Lemma \ref{MutualBounds} and taking eventually a larger $ R_{\e}$,  we see that the sequence $ \breve{\tilde{\psi}}_{n _k } $ satisfies   
(\ref{conc1}). 
As in the proof of Theorem \ref{T5.2} we infer that there is a function $ \zeta \in H_{loc}^1( \R) $ such that $ \zeta '\in L^2 ( \R) $ and 
we may extract a subsequence, still denoted the same, 
such that (\ref{conc2}) holds for $  \breve{\tilde{\psi}}_{n _k }  $ and $ \zeta $. 
It is then clear that (\ref{conc3}) and (\ref{conc4}) hold for  $  \breve{\tilde{\psi}}_{n _k }  $ and $ \zeta $, too, and consequently 
$\zeta \in \Er$. 

Fix $ R_ 1 > 0 $ such that $ \breve{\tilde{\psi}}_{n _k } $ and $ \zeta $ satisfy 
(\ref{conc1bis}) on $ ( - \infty, - R_1] \cup [R_1, \infty). $
Then the functions $ \breve{\tilde{\psi}}_{n _k } $ and $ \zeta $ admit liftings on $ ( - \infty, - R_1] \cup [R_1, \infty) $, say 
$ \breve{\tilde{\psi}}_{n _k } = \rho_k e^{ i \theta _k } $ and $ \zeta = \rho e^{ i \theta}$, 
where $ \rho, \rho _ k \in H^1 (  ( - \infty, - R_1] \cup [R_1, \infty) )$, and $ \theta, \theta_k $ are continuous and 
$ \theta', \theta_k ' \in  L^2 (  ( - \infty, - R_1] \cup [R_1, \infty) )$.

Fix $ \e > 0 $. 
Then choose $ R_{\e} > 0 $  and $ k_{\e } \in \N$ such that  
$  \tilde{\psi}_{n_k }  $ and $ \psi $ satisfy (\ref{conv21}) and 
$ \breve{\tilde{\psi}}_{n _k } $ and $ \zeta $ satisfy (\ref{conc1}) and (\ref{conc1bis}) on $ ( - \infty, - R_{\e} ] \cup [R_{\e}, \infty )$ for all  $k \ges k_{\e}$. 
We proceed as in the proof of Theorem \ref{T5.2}. 
We have $ \breve{\tilde{\psi}}_{n _k }  ( \pm R_{\e} ) /  \zeta ( \pm R_{\e} ) \lra 1$ as $ k \lra \infty$, 
thus we may replace $ \theta_k $ by $ \theta _k + 2 \ell_k^{\pm}  \pi $ on $ ( - \infty, - R_{\e} ] $ and on $  [R_{\e}, \infty )$, where $ \ell_k^{\pm } \in \Z$, 
 so that $ \al_k ^{ \pm } := \theta_k ( \pm R_{\e}) - \theta  ( \pm R_{\e}) \lra 0 $ as $ k \lra \infty$. 
 We extend $ \theta_k $ and $ \theta $ as affine functions on $[ - R_{\e}, R_{\e}].$
Then we have $ \theta_k, \theta \in \dot{H}^1( \R)$, and $ \theta _k \lra \theta, $  $ e^{ i \theta_k } \lra e^{ i \theta }$ 
uniformly on $[ - R_{\e}, R_{\e}]$ as well as in $ H^1 ( [- R_{\e}, R_{\e}])$.

We write 
$$
 \tilde{\psi}_{n_k } = e^{ i \theta _k } + w_k \qquad \mbox{ and } \qquad \psi = e^{ i \theta } + w, 
 $$
where 
$$
w_k = \left( \breve{\tilde{\psi}}_{n _k }   - e^{ i \theta _k } \right) + \left( \tilde{\psi}_{n _k } - \breve{\tilde{\psi}}_{n _k } \right), 
\qquad \mbox{ and } \qquad
w =  ( \zeta - e^{ i \theta}) + ( \psi - \zeta). 
$$
It is clear that $ w_k \in H^1( \R \times (0, 1)) $
and $ \nabla w \in L^2( \R \times (0, 1)) $. 
On $ ( - \infty, - R_{\e} ] $ and on $  [R_{\e}, \infty )$ we have 
$  \zeta - e^{ i \theta} = (\rho - 1)e^{ i \theta } \in L^2 (( - \infty, - R_{\e} ] \cup [R_{\e}, \infty ))$ because $ \zeta \in \Er$, and clearly 
$\zeta - e^{ i \theta} \in L^{\infty} ([- R_{\e}, R_{\e}]) \subset L^2 ([- R_{\e}, R_{\e}]) $, hence
$ \zeta - e^{ i \theta} \in L^2( \R)$. 
We have $ \tilde{\psi}_{n _k } - \breve{\tilde{\psi}}_{n _k }  \lra \psi - \zeta $ almost everywhere on $\R ^2$, and
Lemma \ref{PG-bis} gives 
$$ \|  \tilde{\psi}_{n _k } - \breve{\tilde{\psi}}_{n _k } \|_{L^2( \R \times (0, 1))} ^2= \| v_{ \tilde{\psi}_{n _k }}  \|_{L^2( \R \times (0, 1))} ^2 \les C \Big\| \frac{\p  \tilde{\psi}_{n _k }}{\p y } \Big\| _{L^2( \R \times (0, 1))} ^2\les CM. $$
Then Fatou's Lemma implies that $  \psi - \zeta  \in L^2 ( \R \times (0, 1))$, hence $ w \in H^1 ( \R \times (0, 1)) $.

Denoting $ v = \psi - \zeta $, an easy computation shows that on 
$ (( - \infty, - R_{\e} ] \cup [R_{\e}, \infty )) \times (0, 1) $ we have 
$$
d( \theta, w ) = d ( \theta , ( \rho - 1) e^{ i \theta } + v ) 
= ( 1 - \rho ^2 ) \theta ' 
+ \Big\langle i \frac{\p }{\p x } ( \psi + e^{ i \theta }), v \Big\rangle + \langle i \frac{ \p v}{\p x } , (\rho - 1) e^{ i \theta}  \rangle. 
$$
We may estimate each term: 
$$
\big|  ( 1 - \rho ^2 ) \theta '\big| \les \frac 12 |\theta ' |^2 + \frac12 \left( 1 - \rho ^2 \right)^2 , 
$$
$$ 
\Big| \Big\langle i \frac{\p }{\p x } ( \psi + e^{ i \theta }), v \Big\rangle  \Big|
\les \left( \Big| \frac{ \p \psi }{\p x } \Big| + |\theta '| \right) |v | 
\les \frac 12 \Big| \frac{ \p \psi }{\p x } \Big|^2 + \frac 12 |\theta '|^2 + |v|^2, 
$$
$$ 
\big| \langle i \frac{ \p v}{\p x } , (\rho - 1) e^{ i \theta}  \rangle \big|
\les \frac 12 \Big| \frac{ \p v }{\p x } \Big|^2 + \frac 12 ( \rho - 1)^2 
\les \Big| \frac{ \p \psi }{\p x } \Big|^2 +  |\theta '|^2  + \frac 12 ( \rho^2 - 1)^2 .
$$
Since $ \psi $ and $ \zeta$ satisfy (\ref{conv21}) and (\ref{conc1}), respectively, we find 
\beq
\label{conv28}
\ii_{(( - \infty, - R_{\e} ] \cup [R_{\e}, \infty )) \times (0, 1)} \big| d( \theta, w ) \big| \, dx \, dy 
\les C \e, 
\eeq
where $C$ does not depend on $\e$. 

Using the fact that  $\tilde{\psi}_{n _k } $ and $  \breve{\tilde{\psi}}_{n _k }$ also satisfy (\ref{conv21}) and (\ref{conc1}), respectively, 
and proceeding similarly we  find that  $ \theta_k $ and $ w_k $  satisfy (\ref{conv28}), too,  with $C$ independent of $ \e $ and of $k$. 

On the other hand, we have 
\beq
\label{conv29}
\begin{array}{l}
\ds \ii_{[ - R_{\e}, R_{\e} ] \times (0, 1)} d( \theta_k, w_k ) \, dx \, dy 
= \ii_{[ - R_{\e}, R_{\e} ]  \times (0, 1)} - 2 \langle \theta _k ' e^{ i \theta_k} , w_k \rangle + \langle i \frac{\p w_k }{\p x} , w_k \rangle \, dx \, dy 
\\
\\
\lra 
\ds  \ii_{[ - R_{\e}, R_{\e} ]  \times (0, 1)} - 2 \langle \theta  ' e^{ i \theta} , w  \rangle + \langle i \frac{\p w }{\p x} , w \rangle \, dx \, dy 
= \ii_{[ - R_{\e}, R_{\e} ]  \times (0, 1)}  d ( \theta, w) \, dx \, dy 
\end{array}
\eeq
because $  \theta _k ' e^{ i \theta_k} \rightharpoonup \theta  ' e^{ i \theta}  $ and 
 $ \frac{\p w_k }{\p x} \rightharpoonup \frac{\p w }{\p x} $ weakly in $L^2( [ - R_{\e}, R_{\e}] \times [0, 1]$, while 
 $  w_k \lra w $ strongly in $L^2( [ - R_{\e}, R_{\e}] \times [0, 1])$. 
 
From (\ref{conv28}) and (\ref{conv29}) we deduce that for all $k$  sufficiently large  we have 
$$
| q( \theta_k, w_k ) - q( \theta, w ) | = \Big| \ii_{\R \times (0, 1)} d( \theta _k, w _k ) - d( \theta, w) \, dx \, dy \Big| \les ( 2 C + 1 ) \e. 
$$
Hence there exists $ k(\e) \in \N$ such that for all $ k \ges k( \e)$, 
\beq
\label{conv30}
\big| \QR( \tilde{\psi}_{n_k} ) - \QR ( \psi ) \big| \les \big|  q( \theta_k, w_k ) - q( \theta, w ) \big| \les  ( 2 C + 1 ) \e. 
\eeq
We have thus shown that $  \QR ( \psi ) = \ds \lim_ {k \ra \infty } \QR( \tilde{\psi}_{n_k} ) = \pr{p}$. 
Therefore we must have 
$$
E_{\la } ( \psi ) \ges E_{\la, \min} ( p ) = \lim_{k \ra \infty } E_{\la } \left(  \tilde{\psi}_{n_k}\right). 
$$
Comparing the above inequality to (\ref{conv23}), (\ref{conv25})
 and (\ref{conv26}) (with $ I = \R$), we infer that necessarily 
$$
E_{\la } ( \psi ) =  E_{\la, \min} ( p ) = \lim_{k \ra \infty } E_{\la } \left(  \tilde{\psi}_{n_k}\right) 
$$
and 
\beq
\label{conv31}
\int_{ \R \times (0, 1) } \Big| \frac{ \p  \psi }{\p x } \Big|^2  + \la ^2 \Big| \frac{ \p  \psi }{\p y } \Big|^2\, dx \, dy 
=  \lim_{k \ra \infty } 
\int_{ \R \times (0, 1)} \Big| \frac{ \p  \tilde{\psi}_{n_k } }{\p x } \Big|^2  + \la ^2 \Big| \frac{ \p  \tilde{\psi}_{n_k } }{\p y } \Big|^2 \, dx \,  dy , 
\eeq
\beq
\label{conv32}
\quad
\int_{ \R \times (0, 1) }  V(|\psi |^2)  \, dx  \, dy  =  \lim_{k \ra \infty } \int_{\R \times (0, 1)}  V(|\tilde{\psi} _{n_k}|^2)  \, dx \, dy.
\eeq
The weak convergence $ \nabla   \tilde{\psi}_{n_k } \rightharpoonup \nabla \psi $ in $ L^2 ( \R \times (0, 1) )  $ and (\ref{conv31}) 
give the strong convergence $ \nabla   \tilde{\psi}_{n_k } \lra \nabla \psi $ in $ L^2 ( \R \times (0, 1) )  $. 
The fact that $ V \ges 0$, $ V(|\tilde{\psi} _{n_k}|^2)   \lra  V(|\psi |^2) $ almost everywhere on $ \R^2$ and (\ref{conv32}) imply that 
$  V(|\tilde{\psi} _{n_k}|^2)   \lra  V(|\psi |^2) $ in $ L^1 ( \R \times (0, 1)). $

Fix $ \e > 0 $. Let $ R_{\e}$ be as in (\ref{conv21}). Using (\ref{conv21}) and Lemma \ref{MutualBounds} we find
\beq
\label{conv33}
\ii_{(( - \infty, - R_{\e}] \cup [R_{\e} , \infty)) \times [0, 1]} \big| |\tilde{\psi}_{n_k } | - 1 \big|^2 \, dx \, dy \les 
C E_{GL, \la } ^{( - \infty, - R_{\e}] \cup [R_{\e} , \infty) } ( \tilde{\psi}_{n_k } ) 
\les \de( \e), 
\eeq
where $ \de( \e ) \lra 0 $ as $ \e \lra 0 $. 
It is obvious that a similar estimate holds for $\psi$. 
Since $ \tilde{\psi}_{n_k } \lra \psi $ in $ L^2 ([ - R_{\e}, R_{\e} ] \times [0,1]) $, we have 
$ \| \, |\tilde{\psi}_{n_k }  | - |\psi | \|_{L^2 ( ([ - R_{\e}, R_{\e} ] \times [0,1] ) }^2 \les \e $ for all sufficiently large $ k $. 
It is obvious that $ \left( \, |\tilde{\psi}_{n_k }  | - |\psi | \right)^2 \les \frac 12  \big|\tilde{\psi}_{n_k } | - 1 \big|^2 + \frac 12 \big| \, |\psi | - 1\big|^2$ 
and using (\ref{conv33}) we get 
\beq
\label{conv34}
  \big\| \, |\tilde{\psi}_{n_k } | - |\psi |  \big\|_{L^2( \R \times [0,1])}^2 \les 2 \de( \e) + \e
  \eeq
for all sufficiently large $k$. 
We have thus shown that $ |\tilde{\psi}_{n_k } | - |\psi |   \lra 0 $ in $ L^2( \R \times [0,1]) $. 
Since $ \nabla \left(  |\tilde{\psi}_{n_k } | - |\psi | \right) $ is bounded in $ L^2( \R \times [0,1]) $, using the Sobolev embedding and interpolation we infer that $ |\tilde{\psi}_{n_k } | - |\psi | \lra 0 $ in $ L^p( \R \times [0,1]) $  for any $ p \in [2, \infty)$.
\hfill
$\Box$

\begin{Proposition}
\label{P6.8}
Let $ \la > 0 $ and $ p \in (0, \pi ]$. 
Assume that  $\psi \in \Ep $ is a minimizer of $E_{\la} $ in the set 
$ \{ \phi \in \Ep \; \big| \;  \QR( \phi )   = p \}$. Then: 

\smallskip

(i) There is $ c \in [\left( E_{\la, \min} \right)_r' (p),  ( E_{\la , \min})_{\ell}' (p) ] $ such that $ \psi $ satisfies
\beq
\label{6.60}
i c \psi _{x_1} + \Delta \psi + F(|\psi |^2) \psi = 0 \qquad \mbox{ in } \Do' ( \R^2).
\eeq

\smallskip

(ii) Any solution $\psi \in \Ep $ of (\ref{6.60}) satisfies 
$ \psi \in W_{loc}^{2, p } ( \R^2)$ 
for any $ p \in [2, \infty)$, $\psi $ and $ \nabla \psi $ are bounded 
and $ \psi \in C^{1, \al } (\R^N)$ for any $ \al \in [0, 1)$.

\smallskip

(iii) If $E_{\la, \min}(p )< \sqrt{2} p$ and $ c_1 = \left( E_{\la, \min} \right)_r' (p) <  ( E_{\la , \min})_{\ell}' (p) =  c_2$, there exist two minimizers 
$ \psi _1, \, \psi _2 $ 
of $E_{\la} $ in the set 
$ \{ \phi \in \Ep \; \big| \;  \QR( \phi )   = p \}$
that solve (\ref{6.60}) with $ c = c_1$ and $ c = c_2$, respectively. 
\smallskip

(iv) Assume,  in addition, that $F$ is $C^1$.    Then $\psi $ is symmetric in the following sense: there exists $ y_0 \in [0,1]$ such that $ \psi ( x, y_0 + y  )= \psi (x, y_0 - y)$ for any $ x \in \R$ and any $ y \in [0, \frac 12]$. 

\end{Proposition} 

The proof of Proposition \ref{P6.8} is similar to the proof of Proposition 4.14 p. 187 in \cite{CM}, so we omit it. 
We only sketch the proof of (iv). 
Assume that $ \psi = e^{ i \ph } + w $ is a minimizer for $ E_{\la , \min}(p)$.    
Let $ T_{1, y} w$, $ T_{2, y} w$ and $ \Upsilon _{\ph, w}$ be as in Step 0 of  the proof of Proposition \ref{Emin} (viii). 
There exists $ y_ 0 \in [0,1] $ such that  $ \Upsilon _{\ph, w} ( y _0 ) = 0$.
It follows from (\ref{lopes}) that $ q( \ph, T_{1, y_0} w ) = q( \ph, T_{2, y_0} w ) = q( \ph, w)$. 
Denoting $ \psi _j = e^{ i \ph } + T_{j, y_0} w$ for $ j = 1,2$, we see that 
$ \QR ( \psi _1 ) = \QR( \psi _2) = \QR ( \psi )$, hence $ E_{\la } ( \psi _1 ) \ges E_{\la, \min} ( p )$ and 
$ E_{\la } ( \psi _2 ) \ges E_{\la, \min} ( p )$. 
On the other hand, by (\ref{lopes}) we have 
$$
 E_{\la } ( \psi _1 )  +  E_{\la } ( \psi _2)  = 2 E_{\la} (\psi ) =  2E_{\la, \min} ( p )
$$
and we infer that necessarily $ E_{\la } ( \psi _1 )  =  E_{\la } ( \psi _2)  = E_{\la} (\psi )$ and $ \psi _1, \psi _2$ are also minimizers. 
Then arguing exactly as in the proof of Proposition 4.14 (iii) p. 187 in \cite{CM} we infer that $ \psi = \psi _1 = \psi _2$.

\medskip

\noindent 
{\bf Acknowledgements. } MM gratefully acknowledges partial financial support from Institut Universitaire de France (IUF). 
AM is very grateful to Labex CIMI for financial support. 

\bigskip

At the editorial office's request, we add the following statements.

\medskip

\noindent
{\bf Conflict of interest. } There is no conflict of interest. 

\medskip

\noindent
{\bf Data availability. } 
This article is self-contained and does not use any external data, except for the  references listed below. 
All references are published books or research articles. 
The manuscript will be posted on the preprint server arxiv.org.

\end{document}